\documentclass{amsart}

\usepackage{enumerate} %

\title[Harish-Chandra's theorem]
{Generalization of Harish-Chandra's basic theorem for Riemannian symmetric spaces of non-compact type}
\author{Hiroshi Oda}
\address{Faculty of Engineering, Takushoku University, %
815-1, Tatemachi, Hachioji-shi, Tokyo 193-0985, Japan}
\email{hoda@la.takushoku-u.ac.jp}

\usepackage{txfonts}%
\usepackage{mathrsfs}%
\numberwithin{equation}{section}
\numberwithin{figure}{section}
\numberwithin{table}{section}

\theoremstyle{plain}
 \newtheorem{thm}{Theorem}[section]
 \newtheorem{cor}[thm]{Corollary}
 \newtheorem{lem}[thm]{Lemma}
 \newtheorem{prop}[thm]{Proposition}
\theoremstyle{definition}
 \newtheorem{defn}[thm]{Definition}
 \newtheorem{exm}[thm]{Example}
\theoremstyle{remark}
 \newtheorem{rem}[thm]{Remark}

\allowdisplaybreaks

\makeatletter
\def\@ibid@@#1#2]{%
 [ibid., #2]}
\def\@ibid@{%
 \if\@@tmp@@[%
  \def\nxt@todo{\@ibid@@}%
 \else%
  [ibid.]%
  \let\nxt@todo\relax%
 \fi\nxt@todo}
\def\ibid{%
 \futurelet\@@tmp@@\@ibid@}
\makeatother
\def\norbra#1{\/{\normalfont(}#1\/{\normalfont)}}

\DeclareMathOperator{\Ad}{Ad}
\DeclareMathOperator{\ad}{ad}
\DeclareMathOperator{\gr}{gr}
\DeclareMathOperator{\Ann}{Ann}

\DeclareMathOperator{\Hom}{Hom}
\DeclareMathOperator{\End}{End}

\DeclareMathOperator{\sgn}{sgn}
\DeclareMathOperator{\sym}{symm}
\DeclareMathOperator{\Lie}{Lie}
\DeclareMathOperator{\Real}{Re}
\DeclareMathOperator{\diag}{diag}
\DeclareMathOperator{\id}{id}
\DeclareMathOperator{\refl}{ref}
\DeclareMathOperator{\coldet}{\underrightarrow{\text{\normalfont det}}}

\def\triv{{\rm{triv}}}
\def\single{{\rm{single}}}
\def\double{{\rm{double}}}
\def\ss{{\rm{ss}}}

\def\iu{\sqrt{-1}}
\def\Lieg{{\mathfrak g}}
\def\Lieh{{\mathfrak h}}
\def\Lieu{{\mathfrak u}}

\def\Lien{{\mathfrak n}}
\def\Liea{{\mathfrak a}}

\def\Liek{{\mathfrak k}}
\def\Liep{{\mathfrak p}}
\def\Liem{{\mathfrak m}}

\def\Liez{{\mathfrak z}}
\def\Liesl{{\mathfrak{sl}}}
\def\LieSL{{\mathit{SL}}}
\def\Lieso{{\mathfrak{so}}}
\def\CC{{\mathbb C}}
\def\RR{{\mathbb R}}
\def\ZZ{{\mathbb Z}}
\def\simarrow{\xrightarrow{\smash[b]{\lower 0.8ex\hbox{$\sim$}}}}
\def\bmat#1#2#3#4{\begin{pmatrix}#1&#2\\#3&#4\end{pmatrix}}
\def\ang#1#2{\langle{#1},{#2}\rangle}
\def\cang#1#2{\langle\!\langle{#1},{#2}\rangle\!\rangle}
\newcommand{\bigzerot}{%
	\smash{\lower2ex\hbox{\Huge $0$}}}
\newcommand{\bigzerob}{\smash{\hbox{\Huge $0$}}}

\begin{document}

\begin{abstract}
A basic exact sequence by Harish-Chandra
related to the invariant differential operators on a Riemannian symmetric space $G/K$
is generalized for each $K$-type in a certain class which we call `single-petaled'.
The argument also includes 
a further generalization of
Broer's generalization of the Chevalley restriction theorem.
\end{abstract}
\maketitle

\section{Introduction and main results}\label{sec:intro}
Let $\Lieg$ be a real semisimple Lie algebra
and $\theta$ a fixed Cartan involution of $\Lieg$.
In this paper the subscript $\CC$ is used for indicating the complexification of a real object.
Denote the universal enveloping algebra of the complex Lie algebra $\Lieg_\CC$ by $U(\Lieg_\CC)$,
the center of $U(\Lieg_\CC)$ by $Z(\Lieg_\CC)$,
and the symmetric algebra of $\Lieg_\CC$ by $S(\Lieg_\CC)$.
Similar notation is used for other complex Lie algebras or vector spaces.
Let $G_{\ad}$ be the adjoint group of $\Lieg$, 
$G_\theta$ the subgroup of the adjoint group of $\Lieg_\CC$
consisting of all the elements that leave $\Lieg$ stable,
and $G$ an arbitrary group such that $G_{\ad}\subset G\subset G_\theta$.
The adjoint action of 
$G_{\ad}$ or $G$ or $G_\theta$ (resp.~$\Lieg_\CC$)
is denoted by $\Ad$ (resp.~$\ad$).
Let $\Lieg=\Liek\oplus\Liep$ be the Cartan decomposition.
Take a maximal Abelian subspace $\Liea$ of $\Liep$
and fix a basis $\Pi$ of the restricted root system $\Sigma$ for $(\Lieg, \Liea)$.
$\Pi$ defines the system $\Sigma^+$ of positive roots.
Put $K=G^\theta$ (the subgroup of $G$ which commutes with $\theta$),
$M=Z_K(\Liea)$ (the centralizer of $\Liea$ in $K$),
and $\Liem=\Lie(M)$.
Using $N_K(\Liea)$ (the normalizer of $\Liea$ in $K$)
define the Weyl group $W=N_K(\Liea)/M$.
For each $\alpha \in \Sigma$,
let $H_\alpha$ be the corresponding element of $\Liea$
via the Killing form $B(\cdot,\cdot)$ of $\Lieg$
and put $|\alpha|=B(H_\alpha, H_\alpha)^{\frac12}$,
$\alpha^\vee=\frac{2H_\alpha}{|\alpha|^2}$ (the coroot of $\alpha$).
Denote the restricted root space for $\alpha\in\Sigma$ by $\Lieg_\alpha$
and put $\Lien=\sum_{\alpha\in\Sigma^+}\Lieg_\alpha$,
$\rho=\frac12\sum_{\alpha\in\Sigma^+}(\dim\Lieg_\alpha)\alpha$.
Let us now define the map $\gamma$ of $U(\Lieg_\CC)$ into
$S(\Liea_\CC)$ by the projection
\[
U(\Lieg_\CC)=U(\Liea_\CC)\,\oplus\,\bigr(\Lien_\CC U(\Lieg_\CC)+U(\Lieg_\CC)\Liek_\CC\bigl)
\rightarrow U(\Liea_\CC)\simeq S(\Liea_\CC)
\]
followed by the translation
\[
S(\Liea_\CC)\ni f(\lambda)\mapsto f(\lambda+\rho)\in S(\Liea_\CC).
\]
Here we identified $S(\Liea_\CC)$ with the space of holomorphic polynomials
on the dual space $\Liea_\CC^*$ of $\Liea_\CC$.
We call $\gamma$ the {\it Harish-Chandra homomorphism}.
Let $U(\Lieg_\CC)^K$ (resp.~$S(\Liea_\CC)^W$) be the subalgebra of invariants
in $U(\Lieg_\CC)$ (resp.~$S(\Liea_\CC)$)
under the action of $K$ (resp.~$W$).
In this paper the superscript of an operator domain generally indicates
the subspace of invariants.  
Harish-Chandra showed in \cite{HC}
the following exact sequence of algebra homomorphisms\footnote{Actually
his proof targets only the case of $G=G_{\ad}$,
but the general case follows from it
since \eqref{eq:Ch} is valid for any $G$ (cf.~\cite[Proposition 10]{KR}).}:
\begin{equation}\label{eq:HChom}
 0 \rightarrow U(\Lieg_\CC)^K\cap U(\Lieg_\CC)\Liek_\CC
   \rightarrow U(\Lieg_\CC)^K \xrightarrow{\gamma} S(\Liea_\CC)^W \rightarrow 0.
\end{equation}
On the other hand, let $\Liea^\perp$ be the orthogonal complement of $\Liea$ in $\Liep$
relative to the Killing form $B(\cdot,\cdot)$
and $\gamma_0$ the projection of $S(\Liep_\CC)$ onto $S(\Liea_\CC)$ defined by
\[
S(\Liep_\CC)=S(\Liea_\CC)\oplus S(\Liep_\CC)(\Liea^\perp)_\CC
\rightarrow S(\Liea_\CC).
\]
Then the restriction of $\gamma_0$ to $S(\Liep_\CC)^K$
gives the algebra isomorphism
\begin{equation}\label{eq:Ch}
	\gamma_0:S(\Liep_\CC)^K\simarrow S(\Liea_\CC)^W,
\end{equation}
which is known as the {\it Chevalley restriction theorem\/}.
Let $\sym:S(\Lieg_\CC)\rightarrow U(\Lieg_\CC)$ be the symmetrization map.
Then one has the $K$-module decomposition
\begin{equation}\label{eq:K-decomp}
U(\Lieg_\CC) = \sym(S(\Liep_\CC)) \oplus U(\Lieg_\CC)\Liek_\CC,
\end{equation}
so that \eqref{eq:HChom} is considered as a non-commutative counterpart of \eqref{eq:Ch}.
Hereafter we use the same symbol $\mathscr A$
for the three algebras
$U(\Lieg_\CC)^K/U(\Lieg_\CC)^K\cap U(\Lieg_\CC)\Liek_\CC, S(\Liep_\CC)^K$ and $S(\Liea_\CC)^W$
identified with one another.

Note that \eqref{eq:HChom} and \eqref{eq:Ch}
are rewritten as
\begin{equation}\label{eq:HC2}
 0 \rightarrow \Hom_K(\triv, U(\Lieg_\CC)\Liek_\CC)
   \rightarrow \Hom_K(\triv, U(\Lieg_\CC))
   \xrightarrow{\Gamma^\triv} \Hom_W(\triv, S(\Liea_\CC))\rightarrow 0
\end{equation}
and
\begin{equation}\label{eq:Ch2}
	\Gamma_0^\triv:\Hom_K(\triv, S(\Liep_\CC))
	\simarrow\Hom_W(\triv, S(\Liea_\CC)),
\end{equation}
respectively.
Here `$\triv$' denotes the trivial representation of $K$ or $W$ over $\CC$.
The definitions of $\Gamma^\triv$ and $\Gamma_0^\triv$ are clear.

First,
we generalize the Chevalley restriction theorem in the form \eqref{eq:Ch2}.
We say a $K$-type $(\sigma,V)$ is {\it quasi-spherical\/}
if the $\sigma$-isotypic component of $(\Ad|_K,S(\Liep_\CC))$ is not $0$.
From \cite{KR}
we know $(\sigma,V)$ is quasi-spherical
if and only if $V^M\ne0$.
Suppose $(\sigma,V)$ is quasi-spherical.
Then $W$ naturally acts on $V^M$.
Define the map
\begin{equation}\label{eq:gChmap}
	\Gamma_0^\sigma:\Hom_K(V, S(\Liep_\CC))\ni\Phi
	\mapsto \varphi \in \Hom_W(V^M, S(\Liea_\CC))
\end{equation}
so that the image $\varphi$ is given by the composition
\begin{equation}\label{eq:gChmap2}
	\varphi:V^M\hookrightarrow V\xrightarrow{\Phi}
	S(\Liep_\CC)\xrightarrow{\gamma_0}S(\Liea_\CC).
\end{equation}
It is obviously well-defined and if $(\sigma,V)=(\triv,\CC)$,
the map $\Gamma_0^{\triv}$ coincides with \eqref{eq:Ch2}.
Both $\Hom_K(V, S(\Liep_\CC))$ and $\Hom_W(V^M, S(\Liea_\CC))$
have natural $\mathscr A$-module structures coming from multiplication of images
by elements in $S(\Liep_\CC)^K$ or $S(\Liea_\CC)^W$.
Observe that $\Gamma^\sigma_0$ intertwines these $\mathscr A$-module structures.
Let us introduce a new class of $K$-types:
\begin{defn}\label{defn:single-petaled}
Put $\Sigma_1=\Sigma\setminus2\Sigma$.
Choose a subset $\mathscr R$ of $\Sigma_1$
so that $\mathscr R$ intersects each $W$-orbit of $\Sigma_1$.
For each $\alpha\in \mathscr R$
fix $X_\alpha\in\Lieg_\alpha\setminus\{0\}$.
Then we call a quasi-spherical $K$-type $(\sigma, V)$ is
{\it single-petaled} if and only if
\begin{equation}\label{eq:single-petaled}
 \sigma(X_\alpha+\theta X_\alpha)
 \bigl( \sigma(X_\alpha+\theta X_\alpha)^2-2|\alpha|^2 B(X_\alpha,\theta X_\alpha)\bigr)
 \varv=0
 \qquad\forall \varv\in V^M, \forall \alpha\in \mathscr R.
\end{equation}
\end{defn}
\begin{rem}\label{rem:single-petaled}
Suppose $\alpha \in \Sigma^+$.
Then $\Lieg_\alpha$ has the $M$-invariant inner product $-B(\cdot,\theta\,\cdot)$
and $M$ acts transitively on the unit sphere of $\Lieg_\alpha$
if $\dim \Lieg_\alpha>1$ (see \cite[Theorem~2.1.7]{Ko:Mat2}).
Hence multiplying \eqref{eq:single-petaled} by elements in $M$ or scalars,
we see the above definition does not depend on the choice of $\{X_\alpha\}$.
Furthermore, multiplying \eqref{eq:single-petaled} by elements in $N_K(\Liea)$,
we see $\mathscr R$ may be replaced with $\Sigma_1$.
Hence the definition is also independent of the choice of $\mathscr R$.
\end{rem}

Then we have
\begin{thm}\label{thm:Ch}
For any quasi-spherical $(\sigma,V)$,
$\Gamma_0^\sigma$ is injective.
On the other hand, $\Gamma_0^\sigma$ is surjective if and only if
$(\sigma,V)$ is single-petaled.
\end{thm}
This theorem gives
a generalization of Broer's theorem
for a complex semisimple Lie algebra (\cite{Br})
into the case of a Riemannian symmetric space of non-compact type.
As mentioned in a footnote of \cite{Br}, Broer's theorem can also be proved by
using the results of \cite{PRV}.
Similarly, using the results of \cite{Ko:Mat2},
which generalize the results of \cite{PRV} into the Riemannian symmetric case,
we can show Theorem~\ref{thm:Ch} in a purely algebraic manner.
However our proof in \S\ref{sec:Ch} employs an analytic method
modeled on \cite{Da}.
This method leads to further generalizations of Theorem~\ref{thm:Ch}
in some directions (Proposition~\ref{prop:Chinj}, Theorem~\ref{thm:Chsurj}, Theorem~\ref{thm:WhenSurj}).
In particular, 
for some wider class of $K$-types than `single-petaled',
which will be called {\it `quasi-single-petaled'\/},
a result close to Theorem~\ref{thm:Ch} holds.

The generalization of \eqref{eq:HC2}, which is the main theme of this paper,
requires the notion of
the {\it degenerate affine Hecke algebra\/} $\mathbf H$ associated naturally to the data $(\Lien, \Liea)$
(Definition~\ref{defn:H}).
We here state a few properties of $\mathbf H$.
$\mathbf H$ is an algebra over $\CC$ including $S(\Liea_\CC)$ and 
the group algebra $\CC[W]$ of $W$ as subalgebras with the same $1$
and the center of $\mathbf H$ is $S(\Liea_\CC)^W$.
The map $S(\Liea_\CC)\otimes\CC[W]\rightarrow\mathbf H$ defined by multiplication
gives a $\CC$-linear isomorphism.
Hence the left $\mathbf H$-module
\begin{equation}\label{eq:Sh}
S_{\mathbf H}(\Liea_\CC)\coloneqq\mathbf H\ \Bigm/\!\!\sum_{w\in W\setminus\{1\}}\! \mathbf H(w-1)
\simeq S(\Liea_\CC)\otimes\CC[W] \Bigm/ S(\Liea_\CC)\otimes \!\!\sum_{w\in W\setminus\{1\}}\! \CC[W](w-1)
\end{equation}
is naturally identified with $S(\Liea_\CC)$ as a left $S(\Liea_\CC)$-module.
It is notable that although the left $W$-action on $S_{\mathbf H}(\Liea_\CC)$
differs from the original $W$-action on $S(\Liea_\CC)$,
the space of $W$-fixed elements in $S_{\mathbf H}(\Liea_\CC)$
equals $S(\Liea_\CC)^W$ (Corollary~\ref{cor:ShW}).
Hence we may replace $\Hom_W(\triv, S(\Liea_\CC))$ in \eqref{eq:HC2}
with $\Hom_W(\triv, S_{\mathbf H}(\Liea_\CC))$.

Suppose $(\sigma,V)$ is a quasi-spherical $K$-type.
Define the map
\begin{equation}\label{eq:gHChom}
	\Gamma^\sigma:\Hom_K(V, U(\Lieg_\CC))\ni\Psi
	\mapsto \psi \in \Hom_\CC(V^M, S_{\mathbf H}(\Liea_\CC))
\end{equation}
so that the image $\psi$ is given by the composition
\begin{equation}\label{eq:gHChom2}
	\psi:V^M \hookrightarrow V\xrightarrow{\Psi}
	U(\Lieg_\CC)\xrightarrow{\gamma}S(\Liea_\CC)\simeq S_{\mathbf H}(\Liea_\CC).
\end{equation}
Note that the space in the right-hand side of \eqref{eq:gHChom} is not
$\Hom_W(V^M, S_{\mathbf H}(\Liea_\CC))$.
In fact the map $\psi$ defined by \eqref{eq:gHChom2}
does not always commute with the $W$-actions.
Now we state the main result of this paper:
\begin{thm}\label{thm:HChom}
For any $(\sigma,V)$,
the kernel of\/ $\Gamma^\sigma$ equals $\Hom_K(V, U(\Lieg_\CC)\Liek_\CC)$.
The image of\/ $\Gamma^\sigma$ is included
in $\Hom_W(V^M, S_{\mathbf H}(\Liea_\CC))$
if and only if $(\sigma,V)$ is single-petaled.
If this condition is satisfied,
the image equals $\Hom_W(V^M, S_{\mathbf H}(\Liea_\CC))$
and therefore we have
the exact sequence
\begin{equation}\label{eq:gHChom3}
 0 \rightarrow \Hom_K(V, U(\Lieg_\CC)\Liek_\CC)
   \rightarrow \Hom_K(V, U(\Lieg_\CC))
   \xrightarrow{\Gamma^\sigma} \Hom_W(V^M, S_{\mathbf H}(\Liea_\CC))
   \rightarrow 0.
\end{equation}
\end{thm}
\begin{rem}
If $D\in U(\Lieg_\CC)^K$ and $\Psi\in\Hom_K(V,U(\Lieg_\CC))$,
then the right multiplication of the image of $\Psi$ by $D$
gives a new element $\Psi\cdot D \in \Hom_K(V,U(\Lieg_\CC))$.
This $U(\Lieg_\CC)^K$-module structure of $\Hom_K(V,U(\Lieg_\CC))$
induces an $\mathscr A$-module structure of
$\Hom_K(V,U(\Lieg_\CC)/U(\Lieg_\CC)\Liek_\CC)
\simeq \Hom_K(V,U(\Lieg_\CC))/\Hom_K(V,U(\Lieg_\CC)\Liek_\CC)$.
Also, we naturally consider $\Hom_\CC(V^M, S_{\mathbf H}(\Liea_\CC))$ and
$\Hom_W(V^M, S_{\mathbf H}(\Liea_\CC))$ as $\mathscr A$-modules.
Since it is clear that
\begin{equation}\label{eq:Uk-inv}
\Gamma^\sigma(\Psi\cdot D)=\Gamma^\sigma(\Psi)\cdot\gamma(D)
\qquad \forall D\in U(\Lieg_\CC)^K, \forall \Psi\in\Hom_K(V,U(\Lieg_\CC)),
\end{equation}
$\Gamma^\sigma$ induces an $\mathscr A$-homomorphism
$\Hom_K(V,U(\Lieg_\CC)/U(\Lieg_\CC)\Liek_\CC)
\rightarrow \Hom_\CC(V^M, S_{\mathbf H}(\Liea_\CC))$.
Moreover, if $(\sigma,V)$ is single-petaled,
we get a natural $\mathscr A$-isomorphism
\[
\Hom_K(V,U(\Lieg_\CC)/U(\Lieg_\CC)\Liek_\CC)
\simarrow \Hom_W(V^M, S_{\mathbf H}(\Liea_\CC)).
\]
\end{rem}
The proof of Theorem~\ref{thm:HChom} is given in \S\ref{sec:HC}
with a related result on the quasi-single-petaled $K$-types (Theorem~\ref{thm:HCsurj}).
If $\Lieg$ is a complex semisimple Lie algebra,
then a quasi-spherical $K$-type is naturally identified with 
a finite-dimensional irreducible holomorphic representation of $G$.
Under this identification,
a single-petaled $K$-type is nothing but an irreducible {\it small\/} representation of $G$
in the sense of \cite{Br} (Corollary~\ref{cor:small}).
Moreover, in this case we can deduce from Theorem~\ref{thm:HChom}
a generalization of the celebrated {\it Harish-Chandra isomorphism} (Theorem~\ref{thm:HCisomo}).
In \S\ref{sec:complex}
we also study two topics related to the generalized Harish-Chandra
isomorphism---construction of new kinds of non-commutative determinants,
and a natural correspondence between the submodules of the Verma module $M(\lambda)$ of $U(\Lieg)$
and
the submodules of a certain basic module $A(\lambda)$
of the degenerate affine Hecke algebra $\Tilde{\mathbf H}$ associated to this complex case.
\subsection*{Acknowledgments}
The author would like to express his sincere appreciation to Professor~T.~Oshima and Professor~H.~Matsumoto
for their guidance and encouragement.
The author is also indebted to Dr~T.~Honda for many valuable discussions.
\section{Quasi-spherical K-types}\label{sec:Kostant}
We shall prepare some results on quasi-spherical K-types
which will be used in the subsequent sections.
Most of the results in this section are known.

Identify the $K$-module $S(\Liep_\CC)$
with the $K$-module $\mathscr{P}(\Liep)$ of $\CC$-valued polynomial functions on $\Liep$
via the Killing form.
Each $X\in\Liep$ defines the partial differential operator $\partial(X)$ on $\Liep$.
Extend the correspondence $\partial : X\mapsto\partial(X)$
to the algebra homomorphism from $S(\Liep_\CC)$ to the algebra of partial differential operators
on $\Liep$.

We say an element in $S(\Liep_\CC)\simeq\mathscr P(\Liep)$
is $K$-{\it harmonic\/} if it is killed by 
$\partial(F)$ for any $F\in S(\Liep_\CC)^K\cap S(\Liep_\CC)\Liep_\CC$.
Let $\mathscr H_K(\Liep)$ denote the set of $K$-harmonics.
Note that $\mathscr H_K(\Liep)$ %
is independent of the choice of $G$ ($G_{\ad}\subset G\subset G_\theta$).
The following is essentially due to \cite{KR}:
\begin{prop}\label{prop:KR}
The map
\[
S(\Liep_\CC)^K \otimes \mathscr H_K(\Liep)\rightarrow S(\Liep_\CC)
\]
defined by multiplication is a $K$-module isomorphism.
Moreover, for any finite-dimensional representation $(\sigma, V)$
of $K$ over $\CC$,
$\dim_\CC \Hom_K(V,\mathscr H_K(\Liep)) = \dim_\CC V^M$.
Hence 
\begin{equation}\label{eq:Ast1}
\Hom_K(V,S(\Liep_\CC)) \simeq \mathscr A \otimes \Hom_K(V,\mathscr H_K(\Liep))
\simeq \mathscr A^{\oplus m(\sigma)}\quad\text{with }m(\sigma)=\dim_\CC V^M.
\end{equation}
\end{prop}
\begin{cor}\label{cor:UgDecomp}
As $K$-modules,
$U(\Lieg_\CC)=U(\Lieg_\CC)\Liek_\CC\oplus \sym(\mathscr H_K(\Liep))
\otimes \sym(S(\Liep_\CC)^K).$
\end{cor}
Let $(\sigma,V)$ a quasi-spherical $K$-type
and put $m(\sigma)=\dim_\CC V^M$.
Let $\{\varv_1,\ldots,\varv_{m(\sigma)}\}$ be a basis of $V^M$ and
$\{\Phi_1,\ldots,\Phi_{m(\sigma)}\}$ a basis of
$\Hom_{K}(V,\mathscr H_{K}(\Liep))$.
We put $\Psi_j=\sym\circ\,\Phi_j\in\Hom_{K}(V,U(\Lieg))$
($j=1,\ldots, m(\sigma)$).
In \cite{Ko:Mat, Ko:Mat2},
Kostant studied
the $S(\Liea_\CC)$-valued $m(\sigma)\times m(\sigma)$-matrix
$P^\sigma=\bigl(\gamma\circ\Psi_j[\varv_i]\bigr)_{1\leq i,j\leq m(\sigma)}$,
which is closely related to the theme of the present paper.
In particular he determined the value of $\det P^\sigma$.
It is clear that $\det P^\sigma$, up to a scalar multiple,
does not differ for any choice of bases.
Let $\ang{\cdot}{\cdot}$ be the $\CC$-bilinear form on $\Liea_\CC^* \times \Liea_\CC^*$
induced from $B(\cdot,\cdot)$.
\begin{prop}\label{prop:Kdet}
Suppose $\lambda\in\Liea_\CC^*$ satisfies 
$\Real\ang{\lambda}{\alpha}\ge 0$ for any $\alpha\in\Sigma^+$. 
Then $(\det P^\sigma)(\lambda)\ne0$ for any $(\sigma,V)$.
\end{prop}
\begin{proof}
The proposition is a direct consequence of
\cite{Ko:Mat, Ko:Mat2} if $G=G_\theta$.
We shall translate this result to the general case.
Denote $K, M$ for $G_\theta$
by $K_\theta, M_\theta$.
Let $F$ be the subgroup of 
the adjoint group of $\Lieg_\CC$
consisting of all elements of $\exp \Liea_\CC$
with order not greater than $2$.
Then $K_\theta=KF$, $M_\theta=MF$,
and $F$ normalizes $K$ and $M$ (\cite[Proposition 1, Lemma 20]{KR}).
Since $F$ is isomorphic to a direct product of $\ZZ/2\ZZ$,
we can choose a subgroup $F_1$ so that $K_\theta=K\rtimes F_1$ and $M_\theta=M\times F_1$.
Suppose $(\sigma,V)$ is a quasi-spherical $K$-type.
Then $V_\theta\coloneqq\CC[F_1]\otimes V$
has the natural $K_\theta$-module structure $\sigma_\theta$ defined by
$\sigma_\theta(ka)(a'\otimes \varv)=aa'\otimes(a'akaa')\varv$ for $k\in K, a,a'\in F, \varv\in V$.
Observe that $(V_\theta)^M=\bigoplus_{a\in F_1}a\otimes V^M$ and
$(V_\theta)^{M_\theta}=\left(\frac1{\#F_1}\sum_{a\in F_1}a\right)\otimes V^M$.
Let $(\sigma_\theta,V_\theta)=(\sigma_1,V_1)\oplus\cdots\oplus (\sigma_{t'},V_{t'})$ be
an irreducible decomposition as a $K_\theta$-module
such that $(\sigma_1,V_1),\ldots,(\sigma_t,V_t)$ are just all of the quasi-spherical
components.
Since $\left(\frac1{\#F_1}\sum_{a\in F_1}a\right)\otimes V^M=(V_1)^{M_\theta}\oplus\cdots \oplus (V_t)^{M_\theta}$,
we have $m(\sigma)=m(\sigma_1)+\cdots+m(\sigma_t)$.
Take a basis $\{\varv_1^{(1)},\ldots,\varv_{m(\sigma_1)}^{(1)},\ldots,\varv_1^{(t)},\ldots,\varv_{m(\sigma_t)}^{(t)}\}$ of $V^M$
so that
$\Bigl\{ \left(\frac1{\#F_1}\sum_{a\in F_1}a\right)\otimes \varv_1^{(s)},\ldots,
 \left(\frac1{\#F_1}\sum_{a\in F_1}a\right)\otimes \varv_{m(\sigma_s)}^{(s)} \Bigr\}$
forms a basis of $(V_s)^{M_\theta}$ for each $s=1,\ldots,t$.
On the other hand,
since $(\sigma_\theta,V_\theta)$ is naturally considered as the induced $K_\theta$-module
from $(\sigma,V)$,
by defining $K$-homomorphisms
$\iota_s : V\hookrightarrow V_\theta \xrightarrow{\text{projection}} V_s$ ($s=1,\ldots,t$),
we get the isomorphism
\[
\bigoplus_{s=1}^t\Hom_{K_\theta}(V_s,\mathscr H_{K_\theta}(\Liep)) \ni (\Phi^{(1)},\ldots,\Phi^{(t)})
 \mapsto \Phi^{(1)}\circ\iota_1+\cdots+\Phi^{(t)}\circ\iota_t \in \Hom_K(V,\mathscr H_K(\Liep)).
\]
Therefore, if we take a basis $\{\Phi^{(s)}_1,\ldots,\Phi^{(s)}_{m(\sigma_s)}\}$ of
$\Hom_{K_\theta}(V_s,\mathscr H_{K_\theta}(\Liep))$ for each $s=1,\ldots,t$,
then
\[
\left\{\Phi^{(1)}_1\circ\iota_1,\ldots,\Phi^{(1)}_{m(\sigma_1)}\circ\iota_1,
\ldots, \Phi^{(t)}_1\circ\iota_t,\ldots,\Phi^{(t)}_{m(\sigma_t)}\circ\iota_t \right\}
\]
is a basis of $\Hom_K(V,\mathscr H_K(\Liep))$.
Let us consider $P^\sigma$ with respect to this basis and the basis of $V^M$ defined above. 
Note that $M_\theta$ acts on $S(\Liea_\CC)$ trivially
and that $\gamma$ is an $M_\theta$-homomorphism
because $M_\theta$ normalizes $\Liek_\CC$ and $\Lien_\CC$.
Hence,
$\gamma(D)=\gamma\left(\left(\frac1{\#F_1}\sum_{a\in F_1}a\right) D\right)$ for any $D\in U(\Lieg)$.
We thereby get for
$s, s'=1,\ldots,t$, $i=1,\ldots,m(\sigma_{s'})$, and $j=1,\ldots,m(\sigma_s)$,
\[
\gamma\circ\sym\circ\,\Phi^{(s)}_j\circ\iota_s\,[\varv_i^{(s')}]=
 \begin{cases}
  \ \gamma\circ\sym\circ\,\Phi^{(s)}_j\left[\left(\frac1{\#F_1}\sum_{a\in F_1}a\right)\otimes \varv_i^{(s)}\right]
   & \text{ if }s=s', \\
  \ 0 & \text{ if }s\ne s',
 \end{cases}
\]
which implies the equality
\begin{equation}\label{eq:PP}
P^\sigma=\begin{pmatrix}
P^{\sigma_1} & & \bigzerot\\
& \ddots & \\
\bigzerob & & P^{\sigma_t}
\end{pmatrix}.
\end{equation}
Now our claim follows from Kostant's result.
\end{proof}

Next, we consider some basic single-petaled $K$-types.
Clearly the trivial $K$-type is single-petaled.
Other examples are constituents of $(\Ad,\Liep_\CC)$.
\begin{lem}\label{lem:p}
For any $\alpha\in\Sigma$ \norbra{not $\Sigma_1$!},
$X_\alpha\in\Lieg_\alpha$, and $H\in\Liea_\CC$,
\[
\ad(X_\alpha+\theta X_\alpha)\bigl( \ad(X_\alpha+\theta X_\alpha)^2-2|\alpha|^2 B(X_\alpha,\theta X_\alpha)\bigr)H=0.
\]
\end{lem}
\begin{proof}
If $\alpha(H)=0$, then $\ad(X_\alpha+\theta X_\alpha)H=0$.
On the other hand,
\begin{align*}
\bigl( \ad(X_\alpha+\theta X_\alpha)^2&-2|\alpha|^2 B(X_\alpha,\theta X_\alpha)\bigr)H_\alpha \\
 &=|\alpha|^2\ad(X_\alpha+\theta X_\alpha)(-X_\alpha+\theta X_\alpha)
  -2|\alpha|^2 B(X_\alpha,\theta X_\alpha)H_\alpha \\
 &=2|\alpha|^2\bigl( [X_\alpha,\theta X_\alpha] - B(X_\alpha,\theta X_\alpha)H_\alpha \bigr)
  = 0.\qedhere
\end{align*}
\end{proof}
\begin{defn}
In this paper we say $G/K$ is of Hermitian type
if and only if $\Liep$ has a $K$-invariant complex structure.
\end{defn}
\begin{thm}\label{thm:p}
Suppose $\Lieg$ is simple.
\begin{enumerate}[{\normalfont (i)}]
\item\label{item:pHermite}
Suppose $G/K$ is of Hermitian type and
$\Liep$ has a $K$-invariant complex structure $J$.
Extend $J$ to the $\CC$-linear endomorphism on $\Liep_\CC$ and
let $\Liep_\pm\subset\Liep_\CC$ be the eigenspaces of $J$
with eigenvalues $\pm\iu$.
Then $(\Liep_\CC)^M=\Liea_\CC\oplus J\Liea_\CC$.
Moreover, the two $K$-types $(\Ad,\Liep_\pm)$ are single-petaled
and $(\Liep_\pm)^M\simeq\Liea_\CC$ \norbra{the reflection representation}
as $W$-modules.
\item\label{item:pNonHermite}
Suppose $G/K$ is not of Hermitian type.
Then the $K$-type $(\Ad,\Liep_\CC)$ is single-petaled and
$(\Liep_\CC)^M\simeq\Liea_\CC$ \norbra{the reflection representation}
as a $W$-module.
\end{enumerate}
\end{thm}
\begin{proof}
If $G=G_{\ad}$, then
$(\Liep_\CC)^M$ is as stated in the proposition
by \cite[Proposition~4.1]{Joh}.

\noindent(\ref{item:pHermite})
Suppose $G/K$ is of Hermitian type.
Then $G_{\ad}/(K\cap G_{\ad})$ is also of Hermitian type.
Hence by \cite{Joh},
$(\Liep_\CC)^{M\cap G_{\ad}}=\Liea_\CC\oplus J\Liea_\CC$.
Since $J$ and the $K$-action on $\Liep_\CC$ are commutative,
$(\Liep_\CC)^M=\Liea_\CC\oplus J\Liea_\CC$.
The rest is clear from the $K$-isomorphism
$\frac{1\mp\iu J}2 : \Liep \simarrow \Liep_\pm$
and Lemma~\ref{lem:p}.

\noindent(\ref{item:pNonHermite})
Suppose $G/K$ is not of Hermitian type.
If $G_{\ad}/(K\cap G_{\ad})$ is of Hermitian type,
then its complex structure $J$ gives
the $(K\cap G_{\ad})$-module decomposition
$\Liep_\CC=\Liep_+\oplus \Liep_-$.
In this case, from (\ref{item:pHermite}) and Proposition~\ref{prop:KR},
$\dim_\CC\Hom_{K\cap G_{\ad}}(\Liep_+,\mathscr H_K(\Liep))=\dim_\CC \Liea_\CC$.
Since the $K$-action on $\Liep_\CC$ does not commute with $J$,
$\Liep_\CC$ is an irreducible $K$-module.
Hence we get a natural injection 
$\Hom_K(\Liep_\CC,\mathscr H_K(\Liep))\rightarrow\Hom_{K\cap G_{\ad}}(\Liep_+,\mathscr H_K(\Liep))$,
which implies 
$\dim_\CC (\Liep_\CC)^M\le \dim_\CC \Liea_\CC$ in view of Proposition~\ref{prop:KR}.
But since $(\Liep_\CC)^M\supset \Liea_\CC$, we have $(\Liep_\CC)^M=\Liea_\CC$.

On the other hand, if $G_{\ad}/(K\cap G_{\ad})$ is not of Hermitian type,
then $(\Liep_\CC)^M=\Liea_\CC$ since
$(\Liep_\CC)^{(M\cap G_{\ad})}=\Liea_\CC$ by \cite{Joh}.
The rest is clear from Lemma~\ref{lem:p}.
\end{proof}
In the remainder of this section, we assume $\Lieg$ has real rank $1$
and give a close look at its quasi-spherical $K$-types.
Let $\alpha$ be the unique element in $\Sigma_1\cap\Sigma^+$ and choose
$X_\alpha\in\Lieg_\alpha$
so that $B(X_\alpha,\theta X_\alpha)=-\frac1{2|\alpha|^2}$.
If we put $Z=\iu X_\alpha+\iu\theta X_\alpha$,
then the condition~\eqref{eq:single-petaled} for 
a quasi-spherical $(\sigma,V)$ to be single-petaled 
is reduced to $Z(Z^2-1)V^M=0$.
\begin{lem}\label{lem:rank1}
Suppose $(\sigma,V)$ is a quasi-spherical $K$-type.
\begin{enumerate}[{\normalfont (i)}]
\item\label{item:rank1eigen}
All the eigenvalues of $\sigma(Z)$ are integers.
We denote the largest one by $e(\sigma)$.
\item\label{item:rank1one}
$\dim_\CC V^M=1$.
Hence
$\varv^\sigma \in V^M\setminus\{0\}$ and $\Phi^\sigma\in\Hom_K(V,\mathscr H_K(\Liep))\setminus\{0\}$
are uniquely determined up to scalar multiples. 
\item\label{item:rank1intersect}
Let $V^{Z-e(\sigma)}$ be the eigenspace of $\sigma(Z)$ with eigenvalue $e(\sigma)$.
Let $\bigl(V^M\bigr)^\perp$ denote the orthogonal complement of $V^M$ in $V$
with respect to some $K$-invariant Hermitian inner product $(\cdot,\cdot)_V$ on $V$.
Then $V^{Z-e(\sigma)}\not\subset \bigl(V^M\bigr)^\perp$.
\item\label{item:rank1HCimage}
Put $ \delta=\dim\Lieg_{2\alpha}$ and
$h=\frac{\alpha^\vee}2+\frac{\dim\Lieg_\alpha}2\in S(\Liea_\CC)$
\norbra{recall $\alpha^\vee=\frac{2H_\alpha}{|\alpha|^2}$}.
Then we can choose a pair $(i,j)$ of non-negative integers with $2i+j=|e(\sigma)|$
so that $\gamma\circ\sym\circ\,\Phi^\sigma[\varv^\sigma]$ \norbra{$=\det P^\sigma$} equals
\begin{equation}\label{eq:rank1HCimage}
\bigl[ (h+\delta)(h+\delta+2)\cdots(h+\delta+2(i+j)-2) \bigr]
 \cdot \bigl[ (h+1)(h+3)\cdots(h+2i-1) \bigr]
\end{equation}
up to a scalar multiple.
\item\label{item:rank1small}
$(\sigma,V)$ is the trivial $K$-type $ \Leftrightarrow e(\sigma)=0$.
$(\sigma,V)$ is a constituent of $(\Ad,\Liep_\CC) \Leftrightarrow |e(\sigma)|=1$.
\end{enumerate}
\end{lem}
\begin{proof}
Use the same notation as in the proof of Proposition~\ref{prop:Kdet}.
We may assume $\Lieg$ is simple.
If $G=G_\theta$, then all assertions of the lemma are consequences of \cite{Ko:Mat, Ko:Mat2}.

Suppose $G\ne G_\theta$.
Firstly, we consider the case where $\Lieg\not\simeq\Liesl(2,\RR)$.
In this case each quasi-spherical $K_\theta$-type $(\sigma,V)$
is irreducible as a $K$-module and $V^{M_\theta}=V^M$
(\cite[Chapter II, \S2]{Ko:Mat2}).
Also, because a quasi-spherical $K$-type is a constituent of $S(\Liep_\CC)$,
it must be the restriction of some quasi-spherical $K_\theta$-type.
Hence the lemma follows from the case for $G_\theta$.

Secondly, suppose $\Lieg=\Liesl(2,\RR)$ and $\Liek=\Lieso(2,\RR)$.
Then $G=G_{\ad}=\Ad\bigl(\LieSL(2,\RR)\bigr)$,
$F_1=\{1,a\}$ with
$a=\Ad\left(\begin{pmatrix}\iu&0\\0&-\iu\end{pmatrix}\right)$,
and $Z=\pm\begin{pmatrix}0&-\frac12\iu\\\frac12\iu&0\end{pmatrix}$.
For each integer $e$, define the $1$-dimensional 
quasi-spherical $K$-type $(\sigma_e,V_e)$
by $\sigma_e(Z)=e$.
Then a quasi-spherical $K$-type equals some $(\sigma_e,V_e)$.
For each $(\sigma_e,V_e)$,
define the $K_\theta$-module $(V_e)_\theta=1\otimes V_e+a\otimes V_e$
as in the proof of Proposition~\ref{prop:Kdet}.
Then $a\otimes V_e\simeq V_{-e}$ as a $K$-module.
If $e\ne0$, then $(V_e)_\theta$ is an irreducible $K_\theta$-module
and hence $P^{\sigma_e}=P^{(\sigma_e)_\theta}$ by \eqref{eq:PP},
which assures
(\ref{item:rank1eigen})--(\ref{item:rank1HCimage}) for $(\sigma_e,V_e)$.
If $e=0$, then $\sigma_e=\triv$ and clearly
(\ref{item:rank1eigen})--(\ref{item:rank1HCimage}) hold.
It is also clear that 
(\ref{item:rank1small}) follows from the case for $G_\theta$.
\end{proof}
Combining Lemma~\ref{lem:rank1} (\ref{item:rank1intersect}), (\ref{item:rank1small})
and Theorem~\ref{thm:p}, we can conclude
\begin{cor}\label{cor:rank1-single-petaled}
If $\Lieg$ has real rank $1$, then
the trivial $K$-type
and the $K$-types appearing in $(\Ad,\Liep_\CC)$
exhaust all the single-petaled $K$-types.
\end{cor}
\section{The Chevalley restriction theorem}\label{sec:Ch}
The purpose of this section is to prove Theorem~\ref{thm:Ch}.
Although the method is modeled on that of~\cite{Da} in large part,
some points are improved by use of the {\it rational Dunkl operators}.
We note our method is applicable even to the classical case.

Under the setting of \S\ref{sec:intro} suppose $(\sigma,V)$ is a quasi-spherical $K$-type.
Let $\mathscr{F}$ represent one of the following $\CC$-valued function classes:
$\mathscr{C}$ (continuous functions), 
$\mathscr{C}^\infty$ (smooth functions), or $\mathscr{P}$ (polynomial functions).
Define the map
\begin{equation}\label{eq:gChmap3}
\Hom_{K}(V,\mathscr F(\Liep))\ni\Phi
\mapsto \bigl( \varphi: V^M\ni \varv\mapsto \Phi[\varv]|_\Liea \bigr)
\in \Hom_W(V^M,\mathscr F(\Liea)).
\end{equation}
We consider $\Liep$ and $\Liea$ as Euclidean spaces by the Killing form.
Under the natural identifications
$\mathscr{P}(\Liep)\simeq S(\Liep_\CC)$ and $\mathscr{P}(\Liea) \simeq S(\Liea_\CC)$,
$\Gamma_0^\sigma$ in \S\ref{sec:intro} coincides with \eqref{eq:gChmap3}
for $\mathscr F=\mathscr P$.
Hence we use the same symbol $\Gamma_0^\sigma$ for \eqref{eq:gChmap3}
in general cases.
First we shall prove
\begin{prop}\label{prop:Chinj}
	The map $\Gamma_0^\sigma$ for $\mathscr F=\mathscr C, \mathscr C^\infty$ or $\mathscr P$
	is injective.
\end{prop}
\begin{proof}
We may assume $\mathscr F=\mathscr C$.
Let
\begin{equation}\label{eq:V_decomp}
V=V^M \oplus \sum_{\tau\ne\triv} V_\tau
\end{equation}
be the decomposition
into isotypic components of the $M$-module $V=V|_M$
and define the projection map
\begin{equation}\label{eq:Vproj}
p^\sigma:V=V^M\oplus \sum_{\tau\ne\triv} V_\tau
\rightarrow V^M.
\end{equation}
Suppose $\Phi\in\Hom_{K}(V,\mathscr C(\Liep))$.
Then clearly $\Phi[\varv](H)=\Phi[p^\sigma(\varv)](H)$ for any $\varv\in V$ and $H\in\Liea$.
Let $\varphi \in \Hom_W(V^M,\mathscr C(\Liea))$
be the image of $\Phi$.
Since each element $X\in\Liep$ can be written as
$X=\Ad(k)H$ for some $k\in K$ and $H\in\Liea$,
we have for any $\varv\in V$
\begin{equation}\label{eq:lift}
\begin{aligned}
	\Phi[\varv](X)&=\Phi[\varv](\Ad(k)H)=\Phi[\sigma(k^{-1})\varv](H)\\
	&=\Phi[\,p^\sigma\!\left(\sigma(k^{-1})\varv\right)\,](H)
	=\varphi[\,p^\sigma\!\left(\sigma(k^{-1})\varv\right)\,](H).
\end{aligned}
\end{equation}
Thus $\Phi$ can be completely reproduced by $\varphi$.
\end{proof}
To discuss the image of $\Gamma_0^\sigma$
we introduce two $W$-subspaces of $V^M$.
\begin{defn}\label{defn:sd}
Put
\begin{align*}
V_\single^M&=\bigl\{\,\varv\in V^M;\, 
     \sigma(X_\alpha+\theta X_\alpha)
     \bigl( \sigma(X_\alpha+\theta X_\alpha)^2-2|\alpha|^2 B(X_\alpha,\theta X_\alpha)\bigr)
     \varv=0\bigr.\\
     &\qquad\qquad\qquad\qquad\qquad\qquad\qquad\qquad\qquad\qquad
      \bigl.\forall \alpha\in\Sigma_1, \forall X_\alpha\in\Lieg_\alpha \,\bigr\},\\
V_\double^M&=V^M\cap
     \sum\Bigl\{\,
      \sigma(X_\alpha+\theta X_\alpha)
      \bigl( \sigma(X_\alpha+\theta X_\alpha)^2-2|\alpha|^2 B(X_\alpha,\theta X_\alpha)\bigr)
      V;\,\Bigr.\\
     &\qquad\qquad\qquad\qquad\qquad\qquad\qquad\qquad\qquad\qquad
      \Bigl.\alpha\in\Sigma_1, X_\alpha\in\Lieg_\alpha\,
     \Bigr\}.
\end{align*}
\end{defn}
\begin{lem}\label{lem:direct_sum}
$V^M=V_\single^M\oplus V_\double^M$.
\end{lem}
\begin{proof}
Let $(\cdot,\cdot)_V$ be a $K$-invariant Hermitian inner product on $V$.
Then the isotypic components $V^M$ and $V_\tau$ in \eqref{eq:V_decomp}
are orthogonal to one another.
Let $\bigl(V_\double^M\bigr)^\perp$ be
the orthogonal complement of $V_\double^M$ in $V^M$.
Since $\sigma(X_\alpha+\theta X_\alpha)$ is skew-Hermitian
with respect to $(\cdot,\cdot)_V$,
we easily get $V_\single^M \subset \bigl(V_\double^M\bigr)^\perp$.
Conversely, suppose $\varv\in\bigl(V_\double^M\bigr)^\perp$.
Since $V_\double^M$ is the image of the $M$-module
\[
 \sum\Bigl\{\,
      \sigma(X_\alpha+\theta X_\alpha)
      \bigl( \sigma(X_\alpha+\theta X_\alpha)^2-2|\alpha|^2 B(X_\alpha,\theta X_\alpha)\bigr)
      V;\,\alpha\in\Sigma_1,\ X_\alpha\in\Lieg_\alpha\,
     \Bigr\}
\]
under the projection map \eqref{eq:Vproj},
we have for any $\varv'\in V$, $\alpha\in\Sigma_1$, and $X_\alpha\in\Lieg_\alpha$,
\[
\begin{aligned}
\Bigl(\sigma(X_\alpha+\theta X_\alpha)
 \bigl( \sigma(X_\alpha&+\theta X_\alpha)^2-2|\alpha|^2 B(X_\alpha,\theta X_\alpha)\bigr)\varv,\varv'\Bigr)_V\\
&=-\Bigl(\varv, \sigma(X_\alpha+\theta X_\alpha)
    \bigl( \sigma(X_\alpha+\theta X_\alpha)^2-2|\alpha|^2 B(X_\alpha,\theta X_\alpha)\bigr)\varv'\Bigr)_V\\
&=-\Bigl(\varv,p^\sigma \Bigl( \sigma(X_\alpha+\theta X_\alpha)
    \bigl( \sigma(X_\alpha+\theta X_\alpha)^2-2|\alpha|^2 B(X_\alpha,\theta X_\alpha)\bigr)\varv'\Bigr)\Bigr)_V\\
&=0.
\end{aligned}
\]
It shows $\varv\in V_\single^M$.
Thus we get $V_\single^M=\bigl(V_\double^M\bigr)^\perp$.
\end{proof}
From Remark~\ref{rem:single-petaled},
$(\sigma,V)$ is single-petaled if and only if $V^M_\double=0$.
\begin{lem}\label{lem:double_root}
For any $\alpha\in\Sigma$ \norbra{not $\Sigma_1$!}, $\varv\in V_\single^M$, and $X_\alpha\in\Lieg_\alpha$,
\[
\sigma(X_\alpha+\theta X_\alpha)
     \bigl( \sigma(X_\alpha+\theta X_\alpha)^2-2|\alpha|^2 B(X_\alpha,\theta X_\alpha)\bigr)
     \varv=0.
\]
\end{lem}
\begin{proof}
It suffices to show the equality for a root $2\alpha$ with $\alpha\in\Sigma_1$.
Put $\Lieg(\alpha)=\Liem+\RR H_\alpha + \sum_{\beta\in \Sigma\cap \ZZ\alpha}\Lieg_\beta$
and $\Lieg_\ss(\alpha)=[\Lieg(\alpha),\Lieg(\alpha)]$.
Then $\Lieg_\ss(\alpha)$ is a semisimple Lie algebra with real rank $1$.
Let $G_\ss(\alpha)\subset G$ be the analytic subgroup of $\Lieg_\ss(\alpha)$
and put $\Liek_\ss(\alpha)=\Liek\cap\Lieg_\ss(\alpha)$,
$K_\ss(\alpha)=K\cap G_\ss(\alpha)$, and
$M_\ss(\alpha)=Z_{K_\ss(\alpha)}(\RR H_\alpha)$.
Let $U(\Liek_\ss(\alpha)_\CC)\varv=V^{(1)}\oplus\cdots\oplus V^{(t)}$ be an irreducible
decomposition as a $K_\ss(\alpha)$-module and $\varv=\varv^{(1)}+\cdots+\varv^{(t)}$ the corresponding decomposition.
Since $M\cap G_\ss(\alpha)=M_\ss(\alpha)$,
$\varv^{(s)}$ ($s=1,\ldots,t$) is a non-zero $M_\ss(\alpha)$-fixed vector of $V^{(s)}$.
Moreover, since $M_\ss(\alpha)$ includes the center of $G_\ss(\alpha)$,
we can essentially regard each $V^{(s)}$
as a `K-type' of the adjoint group of $\Lieg_\ss(\alpha)$
and apply the results of \S\ref{sec:Kostant} to it.
Choose $X_\alpha\in\Lieg_\alpha$ so that $B(X_\alpha,\theta X_\alpha)=-\frac1{2|\alpha|^2}$
and put $Z=\iu X_\alpha+\iu\theta X_\alpha$.
Since $\varv\in V_\single^M$, $Z(Z^2-1)\varv=0$.
Then $Z(Z^2-1)\varv^{(s)}=0$ for $s=1,\ldots,t$
because $Z\in \Liek_\ss(\alpha)_\CC$.
Hence by Lemma~\ref{lem:rank1} (\ref{item:rank1one})
each $V^{(s)}$ is single-petaled 
as a `K-type' of the adjoint group of $\Lieg_\ss(\alpha)$.
Now it follows from Corollary~\ref{cor:rank1-single-petaled}
that each $V^{(s)}$ is either the trivial $K_\ss(\alpha)$-type
or a $K_\ss(\alpha)$-type appearing in $(\Ad, \Liep_\CC\cap\Lieg_\ss(\alpha)_\CC)$.
Therefore Theorem~\ref{thm:p} and Lemma~\ref{lem:p} imply
\[
\sigma(X_{2\alpha}+\theta X_{2\alpha})
     \bigl( \sigma(X_{2\alpha}+\theta X_{2\alpha})^2-2|2\alpha|^2 B(X_{2\alpha},\theta X_{2\alpha})\bigr)
     \varv^{(s)}=0 \quad\text{for } s=1,\ldots,t.
\]
Thus we get the lemma.
\end{proof}
For any $W$-submodule $V'$ of $V^M$,
we naturally identify
$\Hom_W\bigl(V^M/V',\mathscr F(\Liea)\bigr)$
with the linear space $\bigl\{\,\varphi\in \Hom_W(V^M,\mathscr F(\Liea));\, 
  \varphi[\varv]=0\ \forall \varv\in V'\,\bigr\}$.
Hereafter in this paper,
we repeatedly use similar identifications without notice.
The second assertion of Theorem~\ref{thm:Ch} can be made more precise
as follows:
\begin{thm}\label{thm:Chsurj}
Suppose $\mathscr F=\mathscr C, \mathscr C^\infty$, or $\mathscr P$.
For any $\varphi\in\Hom_W\bigl(V^M/V_\double^M,\mathscr F(\Liea)\bigr)$
there exists a unique $\Phi\in\Hom_K(V,\mathscr F(\Liep))$
such that $\Gamma_0^\sigma(\Phi)=\varphi$.
\end{thm}
The proof is a bit long and
a large part of this section is devoted to it.
Retain the notation in the proof of Proposition~\ref{prop:Chinj}.
We first show the theorem for $\mathscr F=\mathscr C$.
Suppose $\varphi\in\Hom_W\bigl(V^M/V_\double^M,\mathscr C(\Liea)\bigr)$.
For each $\varv\in V$ we define $\Phi_\varv\in \mathscr C(K\times\Liea)$ by
\begin{equation}\label{eq:lift2}
\Phi_\varv(k,H)=\varphi[\,p^\sigma\!\left(\sigma(k^{-1})\varv\right)\,](H)
\qquad\text{ for }(k,H)\in K\times \Liea.
\end{equation}
\begin{lem}\label{lem:Chpf1}
Suppose $k_1, k_2\in K$ and $H_1, H_2\in\Liea$ satisfy
$\Ad(k_1)H_1=\Ad(k_2)H_2$.
Then $\Phi_\varv(k_1,H_1)=\Phi_\varv(k_2,H_2)$ for any $\varv\in V$.
\end{lem}
\begin{proof}
Note that $H_1$ and $H_2$ in the lemma are conjugate
by some element of $N_K(\Liea)$ (\cite[Chapter~VII, Proposition~2.2]{He1}).
By the definition \eqref{eq:lift2}
we see
for any $\varv\in V, k, k_1\in K, \bar w\in N_K(\Liea)$ and $H\in\Liea$,
the following equalities hold:
\begin{align*}
\Phi_\varv(k_1^{-1}k,H)&=\Phi_{\sigma(k_1)\varv}(k,H),\\
\Phi_\varv(k \bar w,H)&=\Phi_\varv(k,wH)
 \qquad\text{with }w\coloneqq \bar w\bmod M \in W.
\end{align*}
Therefore, 
if we  show
\begin{equation}\label{eq:Chpf4}
\Phi_\varv(k,H)=\Phi_\varv(e,H)\quad\text{ for }H\in\Liea,k\in K^H,\text{ and }\varv\in V,
\end{equation}
our claim follows from it.
Here $e$ and $K^H$ in \eqref{eq:Chpf4} are a unit element
and the stabilizer of $H$ in $K$, respectively.
To show \eqref{eq:Chpf4}, fix an arbitrary $H\in\Liea$
and define $\lambda_H\in V^*$ by
\[
	\lambda_H:V\ni \varv\mapsto \Phi_\varv(e,H).
\]
Let $(\sigma^*,V^*)$ be the dual $K$-type of $(\sigma,V)$
and $(\cdot, \cdot)$ the canonical bilinear form on $V^*\times V$.
For $\bar w\in N_K(\Liea)\cap K^H$ and $\varv\in V$,
\[
\begin{aligned}
(\sigma^*(\bar w)\lambda_H, \varv)&=(\lambda_H, \sigma(\bar w^{-1})\varv)=\Phi_{\sigma(\bar w^{-1})\varv}(e,H)\\
 &=\Phi_\varv(\bar w,H)=\Phi_\varv(e,wH)=\Phi_\varv(e,H)=(\lambda_H, \varv).
\end{aligned}
\]
It shows $\sigma^*(\bar w)\lambda_H=\lambda_H$ for $\bar w\in N_K(\Liea)\cap K^H$
and in particular $\lambda_H\in(V^*)^M$.
Furthermore, for $\alpha\in\Sigma_1$, $X_\alpha\in\Lieg_\alpha$, and $\varv\in V$,
\[
\begin{aligned}
\Bigl(\sigma^*(X_\alpha+\theta X_\alpha)
 \bigl( \sigma^*(X_\alpha&+\theta X_\alpha)^2-2|\alpha|^2 B(X_\alpha,\theta X_\alpha)\bigr)\lambda_H,\varv\Bigr)\\
&=-\Bigl(\lambda_H, \sigma(X_\alpha+\theta X_\alpha)
 \bigl( \sigma(X_\alpha+\theta X_\alpha)^2-2|\alpha|^2 B(X_\alpha,\theta X_\alpha)\bigr)\varv\Bigr)\\
&=-\Phi_{\sigma(X_\alpha+\theta X_\alpha)
 \bigl( \sigma(X_\alpha+\theta X_\alpha)^2-2|\alpha|^2 B(X_\alpha,\theta X_\alpha)\bigr)\varv}(e,H)\\
&=-\varphi[\,p^\sigma\!\left(\sigma(X_\alpha+\theta X_\alpha)
 \bigl( \sigma(X_\alpha+\theta X_\alpha)^2-2|\alpha|^2 B(X_\alpha,\theta X_\alpha)\bigr)\varv\right)\,](H)\\
&=0
\end{aligned}
\]
since $p^\sigma\!\left(\sigma(X_\alpha+\theta X_\alpha)
 \bigl( \sigma(X_\alpha+\theta X_\alpha)^2-2|\alpha|^2 B(X_\alpha,\theta X_\alpha)\bigr)\varv\right)
\in V_\double^M$.
Thus $\lambda_H\in(V^*)_\single^M$.

Put $\Sigma^H=\{\alpha\in\Sigma;\,\alpha(H)=0\}$
and take an arbitrary $\alpha\in\Sigma^H\cap\Sigma_1$ and $X_\alpha\in\Lieg_\alpha$.
We shall prove $\sigma^*(X_\alpha+\theta X_\alpha)\lambda_H=0$.
We may assume $B(X_\alpha,\theta X_\alpha)=-\frac1{2|\alpha|^2}$.
Put $Z=\iu X_\alpha+\iu\theta X_\alpha$ and $\bar s_\alpha=\exp(\pi\iu Z)$.
Then $\bar s_\alpha\in N_K(\Liea)\cap K^H$
and hence $\sigma^*(\bar s_\alpha)\lambda_H=\lambda_H$.
Let
\[
\lambda_H=\lambda_H^{(0)}+\lambda_H^{(+)}+\lambda_H^{(-)}
\]
be the decomposition into $\sigma^*(Z)$-eigenvectors
with eigenvalues $0, 1$, and $-1$.
Then we have
\[
\sigma^*(\bar s_\alpha)\lambda_H
=\lambda_H^{(0)}+e^{\pi\iu}\lambda_H^{(+)}
 +e^{-\pi\iu}\lambda_H^{(-)}
=\lambda_H^{(0)}-(\lambda_H^{(+)}+\lambda_H^{(-)}),
\]
which shows $\lambda_H=\lambda_H^{(0)}$
and hence $\sigma^*(Z)\lambda_H=0$.

Note that $\Liek^H\!\coloneqq \Liem
+\sum \bigl\{\, \RR(X_\alpha+\theta X_\alpha);\, \alpha\in\Sigma^H, X_\alpha\in\Lieg_\alpha  \,\bigr\}$
is the Lie algebra corresponding to $K^H$
and is generated by $\Liem$
and $X_\alpha+\theta X_\alpha$ ($\alpha\in\Sigma^H\cap\Sigma_1$, $X_\alpha\in\Lieg_\alpha$).
Hence we get
$\sigma^*(X)\lambda_H=0$ for any $X\in \Liek^H$.
If we define the analytic subgroup $(K^H)_0$ with Lie algebra $\Liek^H$,
a usual argument leads us to
\[
K^H=\bigl(N_K(\Liea)\cap K^H\bigr) \cdot (K^H)_0.
\]
It shows the $K^H$-invariance of $\lambda^H$ and therefore \eqref{eq:Chpf4}.
\end{proof}
\begin{lem}\label{lem:QTop}
The natural topology of $\Liep$ coincides with the quotient topology
of the surjective map $q:K\times\Liea\rightarrow\Liep$ defined by
$q(k,H)=\Ad(k)H$.
\end{lem}
\begin{proof}
Notice that $B(\cdot,\cdot)$ is $K$-invariant.
Hence if we put for any positive number $R$ 
\[
\Liea_R=\{H\in\Liea;\,B(H,H)\leq R\},\qquad
\Liep_R=\{X\in\Liep;\,B(X,X)\leq R\},
\]
then $q\bigl(S\cap (K\times\Liea_R)\bigr)=q(S)\cap \Liep_R$
for any closed subset $S\subset K\times\Liea$.
Here $S\cap (K\times\Liea_R)$ is compact and so is $q(S)\cap \Liep_R$
with respect to the natural topology of $\Liep$.
It implies that $q(S)$ is closed by the natural topology
and hence the lemma. 
\end{proof}
From Lemma~\ref{lem:Chpf1} and Lemma~\ref{lem:QTop}
$\Phi_\varv$ induces the continuous function $\Phi[\varv]$ on $\Liep$ for each $\varv\in V$.
Clearly the correspondence $\Phi:\varv\mapsto\Phi[\varv]$ commutes with the $K$-actions
and satisfies the relation \eqref{eq:lift}. 
Therefore $\Phi$ is a unique element of $\Hom_K(V,\mathscr C(\Liep))$
such that $\Gamma_0^\sigma(\Phi)=\varphi$.
To show Theorem~\ref{thm:Chsurj} for $\mathscr F=\mathscr C^\infty$
we need some preparation.
\begin{defn}\label{defn:Dunkl}
Let $\mathbf k:\Sigma\rightarrow\CC$ be a {\it multiplicity function},
that is, a function which takes the same value on each $W$-orbit of $\Sigma$.
For $\xi\in\Liea$ we define the operator $\mathscr T_{\mathbf k}(\xi)$ acting on
$f\in\mathscr C^\infty(\Liea)$ or $\mathscr D(\Liea)$
(infinitely differentiable functions with compact support) by
\begin{equation}\label{eq:Dunkl}
\mathscr T_{\mathbf k}(\xi)f(H)=\partial(\xi)f(H)
 +\sum_{\alpha\in\Sigma^+}\mathbf k(\alpha)\,\alpha(\xi)\,\frac{f(H)-f(s_\alpha H)}{\alpha(H)},
\end{equation}
where $\partial(\xi)$ is the $\xi$-directional derivative
and $s_\alpha\in W$ is  
the reflection with respect to $\alpha$. 
\end{defn}
\begin{rem}\label{rem:Dunkl}
The result of \eqref{eq:Dunkl} belongs to the original function class
and it holds that $w\mathscr T_{\mathbf k}(\xi)=\mathscr T_{\mathbf k}(w\xi)w$ for any $w\in W$. 
The operator $\mathscr T_{\mathbf k}(\xi)$ is introduced in \cite{Du} and 
is called the rational Dunkl operator.
Is is known that
$\mathscr T_{\mathbf k}(\xi)\mathscr T_{\mathbf k}(\eta)
 =\mathscr T_{\mathbf k}(\eta)\mathscr T_{\mathbf k}(\xi)$
for any $\xi, \eta\in\Liea$. 
In this section we consider only one special case
where ${\mathbf k}(\alpha)=\frac{\dim \Lieg_\alpha}2$.
Hence hereafter we drop the suffix ${\mathbf k}$ in $\mathscr T_{\mathbf k}$.
\end{rem}
\begin{lem}\label{lem:radial}
Let $L_{\Liep}$ is the flat Euclidean Laplacian on ${\Liep}$.
Let $\{\xi_1,\ldots,\xi_\ell\}$ be an orthonormal basis of $\Liea$
and put
$\mathscr L_\Liea=\sum_{i=1}^\ell \mathscr T(\xi_i)^2$.
Suppose $\Phi\in\Hom_K(V,\mathscr C^\infty(\Liep))$
and $\varv\in V_\single^M$.
Then
\[
\bigr.\bigl(L_{\Liep}\Phi[\varv]\bigr)\bigr|_\Liea
 =\mathscr L_\Liea\bigl(\bigl.\Phi[\varv]\bigr|_\Liea\bigr).
\]
\end{lem}
\begin{proof}
Note that for $X\in\Liep$ and $Y\in\Liek$
\begin{equation}\label{eq:difdif}
\begin{aligned}
\Phi[\sigma(Y)\varv](X)&=\left.\frac{d}{dt}\Phi[\sigma(\exp tY)\varv](X)\right|_{t=0}\\
 &=\left.\frac{d}{dt}\Phi[\varv](\Ad(\exp -tY)X)\right|_{t=0}
  =\partial([X,Y])\Phi[\varv](X).
\end{aligned}
\end{equation}
Hence for $H\in\Liea, \alpha\in\Sigma$, and $X_\alpha\in\Lieg_\alpha$ we have
\begin{equation}\label{eq:differential}
\begin{aligned}
\Phi[\sigma(X_\alpha+\theta X_\alpha)^2\varv](H)
 &=\partial([H,X_\alpha+\theta X_\alpha])\Phi[\sigma(X_\alpha+\theta X_\alpha)\varv](H)\\
 &=\alpha(H)\partial(X_\alpha-\theta X_\alpha)\Phi[\sigma(X_\alpha+\theta X_\alpha)\varv](H)\\
 &=\left.\alpha(H)\frac{d}{dt}\Phi[\sigma(X_\alpha+\theta X_\alpha)\varv]\bigl(H+t(X_\alpha-\theta X_\alpha)\bigr)\right|_{t=0}\\
 &=\left.\alpha(H)\frac{d}{dt}\partial\bigl([H+t(X_\alpha-\theta X_\alpha),X_\alpha+\theta X_\alpha]\bigr)
 \Phi[\varv]\bigl(H+t(X_\alpha-\theta X_\alpha)\bigr)\right|_{t=0}\\
 &=\left.\alpha(H)^2\frac{d}{dt}\partial(X_\alpha-\theta X_\alpha)\Phi[\varv]\bigl(H+t(X_\alpha-\theta X_\alpha)\bigr)\right|_{t=0}\\
 &\qquad+\left.\alpha(H)\frac{d}{dt}2t\partial([X_\alpha,\theta X_\alpha])\Phi[\varv]\bigl(H+t(X_\alpha-\theta X_\alpha)\bigr)\right|_{t=0}\\
 &=\alpha(H)^2\partial(X_\alpha-\theta X_\alpha)^2\Phi[\varv](H)
  +2\alpha(H)\partial([X_\alpha,\theta X_\alpha])\Phi[\varv](H)\\
 &=\alpha(H)^2\partial(X_\alpha-\theta X_\alpha)^2\Phi[\varv](H)
  +2\alpha(H)B(X_\alpha,\theta X_\alpha)\partial(H_\alpha)\Phi[\varv](H).
\end{aligned}
\end{equation}
Let $\varv=\varv^{(0)}+\varv^{(+)}+\varv^{(-)}$
be the decomposition into $\sigma(X_\alpha+\theta X_\alpha)$-eigenvectors
with eigenvalues $0, \pm|\alpha|\!\sqrt{2B(X_\alpha,\theta X_\alpha)}$
(Lemma~\ref{lem:double_root} assures there are no other eigenvalues).
Then $\varv - s_\alpha \varv= 2(\varv^{(+)}+\varv^{(-)})$ (cf.~the proof of Lemma~\ref{lem:QTop}) 
and 
\begin{equation}\label{eq:difference}
\begin{aligned}
\sigma(X_\alpha+\theta X_\alpha)^2\varv
 &=\sigma(X_\alpha+\theta X_\alpha)^2(\varv^{(+)}+\varv^{(-)})\\
 &=2|\alpha|^2B(X_\alpha,\theta X_\alpha)(\varv^{(+)}+\varv^{(-)})
 =|\alpha|^2B(X_\alpha,\theta X_\alpha)(1-s_\alpha)\varv.
\end{aligned}
\end{equation}
From \eqref{eq:differential} and \eqref{eq:difference}
we get
\[
\frac{\partial(X_\alpha-\theta X_\alpha)^2}{-2B(X_\alpha,\theta X_\alpha)}\Phi[\varv](H)
 =\frac1{\alpha(H)}\partial(H_\alpha)\Phi[\varv](H)
  -|\alpha|^2\frac{\Phi[\varv](H)-\Phi[\varv](s_\alpha H)}{2\alpha(H)^2}.
\]
Therefore
\[
\begin{aligned}
L_{\Liep}\Phi[\varv](H)&=\sum_{i=1}^\ell \partial(\xi_i)^2\Phi[\varv](H)\\
 &\qquad+ \sum_{\alpha\in\Sigma^+} 
  \frac{\dim \Lieg_\alpha}2 \left( \frac2{\alpha(H)}\partial(H_\alpha)\Phi[\varv](H)
  -|\alpha|^2\frac{\Phi[\varv](H)-\Phi[\varv](s_\alpha H)}{\alpha(H)^2} \right),
\end{aligned}
\]
which equals $\sum_{i=1}^\ell \mathscr T(\xi_i)^2\Phi[\varv](H)$
by~\cite[Theorem~1.10]{Du}.
\end{proof}
Let $dX$ (resp.~$dH$) be the canonical measure of the Euclidean space $\Liep$
(resp.~$\Liea$).
\begin{lem}\label{lem:invint}
There exists a positive constant $C_\Liea$ such that
for any $K$-invariant continuous function $F(X)$ on $\Liep$ with compact support
\[
\int_{\Liep}F(X)\,dX=
C_\Liea\int_\Liea F(H) \prod_{\alpha\in\Sigma^+}\bigl|\alpha(H)\bigr|^{\dim \Lieg_\alpha}dH.
\]
\end{lem}
\begin{proof}
See \cite[Chapter~I, Theorem~5.17]{He2}.
\end{proof}
\begin{lem}\label{lem:Dunkladj}
For any $\varphi \in\mathscr C^\infty(\Liea), f \in\mathscr D(\Liea)$, and $\xi\in\Liea$,
\[
\int_\Liea \bigl(\mathscr T(\xi)\varphi\bigr)(H) f(H) 
 \prod_{\alpha\in\Sigma^+}\bigl|\alpha(H)\bigr|^{\dim \Lieg_\alpha}dH
 = - \int_\Liea \varphi(H) \bigl(\mathscr T(\xi)f\bigr)(H) 
 \prod_{\alpha\in\Sigma^+}\bigl|\alpha(H)\bigr|^{\dim \Lieg_\alpha}dH.  
\]
\end{lem}
\begin{proof}
By a straightforward calculation (cf.~\cite[Lemma~2.9]{Du2}).
\end{proof}
\begin{lem}\label{lem:bases}
Recall the decomposition \eqref{eq:V_decomp}.
Suppose $\{\varv_1,\ldots,\varv_n\}$ is a basis of\/ $V$ such that
$\{\varv_1,\ldots,\varv_{m'}\}, \{\varv_{m'+1},\ldots,\varv_{m}\}$
and $\{\varv_{m+1},\ldots,\varv_n\}$ are bases of\/ $V_\single^M$, $V_\double^M$,
and $\sum_{\tau\ne\triv}V^\tau$, respectively.
Let $\{\varv_1^*,\ldots,\varv_n^*\}$ be the dual basis of $\{\varv_1,\ldots,\varv_n\}$.
Then, 
$\{\varv_1^*,\ldots,\varv_{m'}^*\}, \{\varv_{m'+1}^*,\ldots,\varv_{m}^*\}$
and $\{\varv_{m+1}^*,\ldots,\varv_n^*\}$ are bases of\/ $(V^*)_\single^M$, $(V^*)_\double^M$,
and $\sum_{\tau\ne\triv}(V^*)^\tau$, respectively.
\end{lem}
\begin{proof}
Let $(\cdot, \cdot)$ be the canonical bilinear form on $V^*\times V$.
Then $((V^*)^\tau, V^\mu)=0$ unless $\tau=\mu^*$.
Hence the lemma follows if we show 
\[
\begin{aligned}
V_\single^M&=\{\varv\in V^M;\,(\varv^*,\varv)=0\quad \forall \varv^*\in(V^*)_\double^M\},\\
(V^*)_\single^M&=\{\varv^*\in (V^*)^M;\,(\varv^*, \varv)=0\quad \forall \varv\in V_\double^M\}
\end{aligned}
\]
But this argument is quite similar to the proof of Lemma~\ref{lem:direct_sum}
and we omit it. 
\end{proof}
Now suppose $\varphi\in\Hom_W\bigl(V^M/V_\double^M,\mathscr C^\infty(\Liea)\bigr)$.
Let $\varphi^\sim$ stand for the unique element of $\Hom_K(V,\mathscr C(\Liep))$
such that $\Gamma_0^\sigma(\varphi^\sim)=\varphi$.
It follows from Remark~\ref{rem:Dunkl} that the map
\begin{equation}\label{eq:Lphi}
V^M\ni \varv\mapsto \mathscr L_\Liea\varphi[\varv] \in \mathscr C^\infty(\Liea)
\end{equation}
belongs to $\Hom_W\bigl(V^M/V_\double^M,\mathscr C^\infty(\Liea)\bigr)$.
Denote the map \eqref{eq:Lphi} by $\mathscr L_\Liea\varphi$ and
let us show
\begin{equation}\label{eq:Laplacian}
L_{\Liep}\varphi^\sim[\varv]=(\mathscr L_\Liea \varphi)^\sim[\varv]
\qquad\forall \varv\in V.
\end{equation}
In the left-hand side we consider $\varphi^\sim[\varv]$
an element in $\mathscr D'(\Liep)$
(the space of distributions) on which $L_{\Liep}$ is acting.
On the other hand we know the right-hand side is a continuous function.
Hence by successive use of \eqref{eq:Laplacian} and Weyl's lemma on regularity,
we can conclude $\varphi^\sim[\varv]\in\mathscr C^\infty(\Liep)$.

Let $\{\varv_1,\ldots,\varv_n\}$ be a basis of $V$ as in Lemma~\ref{lem:bases}
and $\{\varv_1^*,\ldots,\varv_n^*\}$ its dual basis.
It suffices to show that
for any $n$ test functions $F_1, \ldots, F_n \in \mathscr D(\Liep)$,
\begin{equation}\label{eq:cuupling}
\sum_{i=1}^n \int_{\Liep} \varphi^\sim[\varv_i] \left(L_{\Liep}F_i\right)\,dX
 = \sum_{i=1}^n \int_{\Liep} (\mathscr L_\Liea \varphi)^\sim[\varv_i] F_i\,dX.
\end{equation}
To do this, using the linear map $F : V^*\rightarrow \mathscr D(\Liep)$ defined by
$\varv_i^*\mapsto F_i$ ($i=1,\ldots,n$), put
\[
\bar F : V^*\ni \varv^*\mapsto \int_K F[\sigma^*(k)\varv^*](\Ad(k)X)\,dk \in \mathscr D(\Liep),
\]
where $dk$ is a normalized measure on $K$.
Then $\bar F\in\Hom_K(V^*,\mathscr C^\infty(\Liep))$
and the left-hand side of \eqref{eq:cuupling} equals
\begin{equation}\label{eq:lhs}
\begin{aligned}
\int_K \sum_{i=1}^n \int_{\Liep} &\varphi^\sim[\sigma(k)\varv_i](X) \left(L_{\Liep}F[\sigma^*(k)\varv_i^*]\right)(X)\,dX\,dk \\
 &=\sum_{i=1}^n \int_K \int_{\Liep} \varphi^\sim[\varv_i](\Ad(k^{-1})X) \left(L_{\Liep}F[\sigma^*(k)\varv_i^*]\right)(X)\,dX\,dk \\
 &=\sum_{i=1}^n \int_{\Liep} \varphi^\sim[\varv_i](X) \int_K \left(L_{\Liep}F[\sigma^*(k)\varv_i^*]\right)(\Ad(k)X)\,dk\,dX \\
 &=\int_{\Liep} \sum_{i=1}^n \varphi^\sim[\varv_i](X) \left(L_{\Liep}\bar F[\varv_i^*]\right)(X)\,dX \\
 &=C_\Liea\int_\Liea \sum_{i=1}^n \varphi^\sim[\varv_i](H) \left(L_{\Liep}\bar F[\varv_i^*]\right)(H) \prod_{\alpha\in\Sigma^+} \bigl|\alpha(H)\bigr|^{\dim \Lieg_\alpha}dH \\
 &=C_\Liea\int_\Liea \sum_{i=1}^{m'} \varphi[\varv_i](H) \left(L_{\Liep}\bar F[\varv_i^*]\right)(H) \prod_{\alpha\in\Sigma^+} \bigl|\alpha(H)\bigr|^{\dim \Lieg_\alpha}dH. 
\end{aligned}
\end{equation}
Here the fourth equality comes from the $K$-invariance of the integrand
and Lemma~\ref{lem:invint}.
The last equality
is based on the fact that $\bigl.\varphi^\sim[\varv_i]\bigr|_\Liea=0$ for $i=m'+1,\ldots,n$
(see the proof of Proposition~\ref{prop:Chinj}).
Similarly the right-hand side of \eqref{eq:cuupling} is changed into the form
\begin{equation}\label{eq:rhs}
C_\Liea\int_\Liea \sum_{i=1}^{m'} (\mathscr L_\Liea \varphi)[\varv_i](H) \bar F[\varv_i^*](H) \prod_{\alpha\in\Sigma^+} \bigl|\alpha(H)\bigr|^{\dim \Lieg_\alpha}dH,
\end{equation}
which equals the final form of \eqref{eq:lhs}
because from Lemma~\ref{lem:Dunkladj} and Lemma~\ref{lem:radial}
we have for $i=1,\ldots,m'$
\[
\begin{aligned}
\int_\Liea ( \mathscr L_\Liea \varphi )[\varv_i](H) \bar F[\varv_i^*] (H) \prod_{\alpha\in\Sigma^+} \bigl|\alpha(H)\bigr|^{\dim \Lieg_\alpha}dH
 &=\int_\Liea \varphi[\varv_i](H) \mathscr L_\Liea \bigl( \bigl.\bar F[\varv_i^*]\bigr|_\Liea \bigr)(H) \prod_{\alpha\in\Sigma^+} \bigl|\alpha(H)\bigr|^{\dim \Lieg_\alpha}dH \\
 &=\int_\Liea \varphi[\varv_i](H) \left(L_{\Liep}\bar F[\varv_i^*]\right)(H) \prod_{\alpha\in\Sigma^+} \bigl|\alpha(H)\bigr|^{\dim \Lieg_\alpha}dH.
\end{aligned}
\]
Thus we get Theorem~\ref{thm:Chsurj} for $\mathscr F=\mathscr C^\infty$.

Finally, to show Theorem~\ref{thm:Chsurj} for $\mathscr F=\mathscr P$,
suppose $\varphi\in\Hom_W\bigl(V^M/V_\double^M,\mathscr P(\Liea)\bigr)$
and put $\Phi=\varphi^\sim$.
Let us prove $\Phi[\varv]\in \mathscr P(\Liep)$ for any $\varv\in V$.
We may assume for any $\varv\in V^M$, 
$\varphi[\varv]$ is a homogenous polynomial of the same degree, say $j$.
Then $\Phi[\varv]$ is also a homogeneous function of degree $j$
for any $\varv\in V$.
It is clear from \eqref{eq:lift2} and the relation
$\Phi[\varv](\Ad(k)H)=\Phi_\varv(k,H)$ for $k\in K$ and $H\in\Liea$.
Since a $\mathscr C^\infty$ homogeneous function defined around $0$
is a polynomial, we get the claim.
This completes the proof of Theorem~\ref{thm:Chsurj}.
\begin{rem}
In our proof of Theorem~\ref{thm:Chsurj},
the results in \S\ref{sec:Kostant} are used to get Lemma~\ref{lem:double_root},
but nowhere else.
Hence if we replace the definitions of $V_\single^M$ and $V_\double^M$ by
\begin{align*}
V_\single^M&=\bigl\{\,\varv\in V^M;\, 
     \sigma(X_\alpha+\theta X_\alpha)
     \bigl( \sigma(X_\alpha+\theta X_\alpha)^2-2|\alpha|^2 B(X_\alpha,\theta X_\alpha)\bigr)
     \varv=0\bigr.\\
     &\qquad\qquad\qquad\qquad\qquad\qquad\qquad\qquad\qquad\qquad
      \bigl.\forall \alpha\in\Sigma, \forall X_\alpha\in\Lieg_\alpha \,\bigr\},\\
V_\double^M&=V^M\cap
     \sum\Bigl\{\,
      \sigma(X_\alpha+\theta X_\alpha)
      \bigl( \sigma(X_\alpha+\theta X_\alpha)^2-2|\alpha|^2 B(X_\alpha,\theta X_\alpha)\bigr)
      V;\,\Bigr.\\
     &\qquad\qquad\qquad\qquad\qquad\qquad\qquad\qquad\qquad\qquad
      \Bigl.\alpha\in\Sigma, X_\alpha\in\Lieg_\alpha\,
     \Bigr\},
\end{align*}
then Theorem~\ref{thm:Chsurj} can be shown without using any result in \S\ref{sec:Kostant}.
By Lemma~\ref{lem:double_root} and a similar argument to the proof of Lemma~\ref{lem:direct_sum},
we can see that the definitions here are equivalent to those in Definition~\ref{defn:sd}.
\end{rem}
To finish the proof of Theorem~\ref{thm:Ch} we shall show
\begin{thm}\label{thm:WhenSurj}
For any $W$-submodule $V'$ of $V^M$,
$\Gamma^\pi_0\left( \Hom_K(V, S(\Liep_\CC)) \right) \supset
\Hom_W(V^M/V', S(\Liea_\CC))$ if and only if
$V'\supset V_\double^M$.
\end{thm}
It suffices to prove the necessity since
Theorem~\ref{thm:Chsurj} implies the sufficiency.

We say an element in $S(\Liea_\CC)\simeq\mathscr P(\Liea)$
is $W$-harmonic if it is killed by $\partial(f)$
for any $f\in S(\Liea_\CC)^W\cap S(\Liea_\CC)\Liea_\CC$.
Let $\mathscr H_W(\Liea)$ denote the set of $W$-harmonics.
\begin{prop}\label{prop:He}
The map
\[
S(\Liea_\CC)^W \otimes \mathscr H_W(\Liea)\rightarrow S(\Liea_\CC)
\]
defined by multiplication is a $W$-module isomorphism.
Moreover $\mathscr H_W(\Liea)\simeq \CC[W]$ as a $W$-module,
where we consider the right-hand side
as the left regular representation of\/ $W$.
Hence
\begin{equation}\label{eq:Ast2}
\Hom_W(V^M,S(\Liea_\CC)) \simeq \mathscr A \otimes \Hom_W(V^M,\mathscr H_W(\Liea))
\simeq \mathscr A^{\oplus m(\sigma)}\quad\text{with }m(\sigma)=\dim_\CC V^M.
\end{equation}
\end{prop}
\begin{proof}
See \cite[Chapter~III,Theorem~3.4]{He2} .
\end{proof}
The necessity of Theorem~\ref{thm:WhenSurj} follows form the next lemma:
\begin{lem}\label{lem:HI}
Suppose $\Phi\in\Hom_K(V,S(\Liep_\CC))$ satisfies
$\Gamma_0^\sigma(\Phi)[\varv]\in \mathscr H_W(\Liea)$ for any $\varv\in V^M$.
Then $\Gamma_0^\sigma(\Phi)[\varv]=0$ for any $\varv\in V_\double^M$.
\end{lem}
\begin{proof}[Proof of Lemma~\ref{lem:HI}]
Recall the notation in Lemma~\ref{lem:radial} and its proof.
Put $L_\Liea=\sum_{i=1}^\ell \partial(\xi_i)^2$.
Notice that \eqref{eq:difdif} and \eqref{eq:differential}
are valid for any $\varv\in V$.
Hence for any $\varv \in V, \alpha\in\Sigma_1$, and $X_\alpha\in\Lieg_\alpha$,
a similar calculation to \eqref{eq:differential} implies
\[
\begin{aligned}
\Phi&[ \sigma(X_\alpha+\theta X_\alpha)
 \bigl( \sigma(X_\alpha+\theta X_\alpha)^2-2|\alpha|^2 B(X_\alpha,\theta X_\alpha)\bigr)\varv ](H)\\
&=\alpha(H)^2\partial(X_\alpha-\theta X_\alpha)^2\Phi[\sigma(X_\alpha+\theta X_\alpha)\varv](H) \\
 &\quad+2\alpha(H)B(X_\alpha,\theta X_\alpha)\partial(H_\alpha)\Phi[\sigma(X_\alpha+\theta X_\alpha)\varv](H)
  -2|\alpha|^2 B(X_\alpha,\theta X_\alpha)\Phi[\sigma(X_\alpha+\theta X_\alpha)\varv](H)\\
&=\alpha(H)^2\partial(X_\alpha-\theta X_\alpha)^2\Phi[\sigma(X_\alpha+\theta X_\alpha)\varv](H) \\
 &\qquad+2\alpha(H)B(X_\alpha,\theta X_\alpha)
  \left.\frac{d}{dt}\Phi[\sigma(X_\alpha+\theta X_\alpha)\varv](H+tH_\alpha)\right|_{t=0}\\
 &\qquad\qquad-2\alpha(H)|\alpha|^2 B(X_\alpha,\theta X_\alpha)\partial(X_\alpha-\theta X_\alpha)\Phi[\varv](H)\\
&=\alpha(H)^2\partial(X_\alpha-\theta X_\alpha)^2\Phi[\sigma(X_\alpha+\theta X_\alpha)\varv](H) \\
 &\qquad+2\alpha(H)B(X_\alpha,\theta X_\alpha)
  \left( |\alpha|^2\partial(X_\alpha-\theta X_\alpha)\Phi[\varv](H) + \alpha(H)\partial(H_\alpha)\partial(X_\alpha-\theta X_\alpha)\Phi[\varv](H) \right)\\
 &\qquad\qquad-2\alpha(H)|\alpha|^2 B(X_\alpha,\theta X_\alpha)\partial(X_\alpha-\theta X_\alpha)\Phi[\varv](H)\\
&=\alpha(H)^2 \left(\partial(X_\alpha-\theta X_\alpha)^2\Phi[\sigma(X_\alpha+\theta X_\alpha)\varv](H)
  +2B(X_\alpha, \theta X_\alpha)\partial(H_\alpha)\partial(X_\alpha-\theta X_\alpha)\Phi[\varv](H)\right).
\end{aligned}
\]
It shows 
\begin{multline*}
\bigl.\Phi[\sigma(X_\alpha+\theta X_\alpha)
 \bigl( \sigma(X_\alpha+\theta X_\alpha)^2-2|\alpha|^2 B(X_\alpha,\theta X_\alpha)\bigr)\varv]\bigr|_\Liea \\
=\bigl.\Phi[p^\sigma\!\left(\sigma(X_\alpha+\theta X_\alpha)
 \bigl( \sigma(X_\alpha+\theta X_\alpha)^2-2|\alpha|^2 B(X_\alpha,\theta X_\alpha)\bigr)\varv\right)]\bigr|_\Liea
 \in \alpha^2 \mathscr P(\Liea).
\end{multline*}
Since no element of $\alpha^2\mathscr P(\Liea)$ other than $0$ is killed by $L_\Liea$,
we have
\[
\Gamma_0^\sigma(\Phi)[p^\sigma\!\left(\sigma(X_\alpha+\theta X_\alpha)
 \bigl( \sigma(X_\alpha+\theta X_\alpha)^2-2|\alpha|^2 B(X_\alpha,\theta X_\alpha)\bigr)\varv\right)]=0,
\]
and the lemma.
\end{proof}
We conclude this section by introducing a new class of $K$-types. 
\begin{defn}
We say a $K$-type $(\sigma, V)$ is quasi-single-petaled when $V_\single^M\ne 0$.
\end{defn}
If $\Lieg$ has real rank $1$, then
Lemma~\ref{lem:rank1} (\ref{item:rank1one}) assures
all the quasi-single-petaled $K$-types are single-petaled
and their number is finite by Corollary~\ref{cor:rank1-single-petaled}.
In general we have
\begin{prop}\label{prop:quasi-single-petaled}
The number of quasi-single-petaled $K$-types is finite.
\end{prop}
\begin{proof}
Suppose $(\sigma, V)$ is quasi-single-petaled.
Then it follows from Theorem~\ref{thm:Chsurj}
there is a nontrivial $\Phi\in\Hom_K(V,S(\Liep_\CC))$
such that $\Gamma_0^\sigma(\Phi)[\varv]\in\mathscr H_W(\Liea)$ for any $\varv\in V^M$.
Recall that the degree of an element of $\mathscr H_W(\Liea)$ 
is not greater than the number of reflections in $W$, say, $r$ (cf. \cite[Chapter~III,Theorem~3.6]{He2}).
Since $\gamma_0$ maps a homogeneous element to a homogeneous element of the same degree,
the degree of $\Phi[\varv]$ for each $\varv\in V$ is at most $r$.
Hence $V$ must be equivalent to an irreducible $K$-subspace 
of $\{F\in S(\Liep_\CC);\,\deg F\leq r\}$.
\end{proof}

\section{The Harish-Chandra homomorphism}\label{sec:HC}

In this section we prove Theorem~\ref{thm:HChom}.
Let us start with the definition of the degenerate affine Hecke algebra,
which is due to \cite{Lu}.
\begin{defn}\label{defn:H}
Let $\mathbf k : \Sigma_1\rightarrow\CC$ be a multiplicity function.
Then there exists uniquely (up to equivalence) an algebra $\mathbf H_{\mathbf k}$ over $\CC$
with the following properties:
\begin{enumerate}[\ (i)]
\item $\mathbf H_{\mathbf k}\simeq S(\Liea_\CC)\otimes\CC[W]$ as a $\CC$-linear space.
\item\label{item:Hsubs}
 The maps $S(\Liea_\CC)\rightarrow \mathbf H_{\mathbf k}, f\mapsto f\otimes1$ and
 $\CC[W]\rightarrow \mathbf H_{\mathbf k}, w\mapsto1\otimes w$ are algebra homomorphisms.
\item $(f\otimes1)\cdot(1\otimes w)=f\otimes w$ for any $f\in S(\Liea_\CC)$ and $w\in W$.
\item\label{item:Hrelorg}
 $(1\otimes s_\alpha)\cdot(\xi\otimes1)
 =s_\alpha(\xi)\otimes s_\alpha-\mathbf k(\alpha)\,\alpha(\xi)$
for any $\alpha\in\Pi$ and $\xi\in\Liea_\CC$.
Here $s_\alpha\in W$ is the reflection corresponding to $\alpha$.
\end{enumerate}
We call $\mathbf H_{\mathbf k}$ the degenerate affine Hecke algebra
associated to the data $(\Liea_\CC, \Pi, \mathbf k)$.
\end{defn}

\begin{rem}
By (\ref{item:Hsubs}) we identify $S(\Liea_\CC)$ and $\CC[W]$
with subalgebras of $\mathbf H_{\mathbf k}$.
Then (\ref{item:Hrelorg}) is simply written as
\begin{equation}\label{eq:Hrel}
s_\alpha\cdot\xi
 =s_\alpha(\xi)\cdot s_\alpha
 - \mathbf k(\alpha)\,\alpha(\xi)\qquad\forall\alpha\in\Pi\ \forall\xi\in\Liea_\CC.
\end{equation}
The center of $\mathbf H_{\mathbf k}$
equals $S(\Liea_\CC)^W$ (\cite[Theorem 6.5]{Lu})
as we stated in \S\ref{sec:intro}.
In this section we fix
\begin{equation}\label{eq:param}
\mathbf k(\alpha)=
\dim \Lieg_\alpha + 2\dim \Lieg_{2\alpha}
\end{equation}
and drop the suffix ${\mathbf k}$ in $\mathbf H_{\mathbf k}$.
Note that $\mathbf H$ is fully determined by the data $(\Lien,\Liea)$.
\end{rem}
As in \S\ref{sec:intro}, we define the left $\mathbf H$-module
$S_{\mathbf H}(\Liea_\CC)$ by \eqref{eq:Sh}.
\begin{lem}\label{lem:Hecke}
Suppose $\alpha\in\Pi$ and put $\Liea(\alpha)=\{H\in\Liea;\,\alpha(H)=0\}$.
Then 
\[
S_{\mathbf H}(\Liea_\CC)=S(\Liea(\alpha)_\CC)\cdot \CC[(\alpha^\vee)^2]
\,\oplus\, S(\Liea(\alpha)_\CC)\cdot \CC[(\alpha^\vee)^2](\alpha^\vee+\dim \Lieg_\alpha + 2\dim \Lieg_{2\alpha})
\]
is the decomposition into the eigenspaces of $s_\alpha\in\mathbf H$
with eigenvalues $1,-1$.
\end{lem}
\begin{proof}
Using \eqref{eq:Hrel}, we have
\[
s_\alpha(\alpha^\vee+\mathbf k(\alpha))=
-\alpha^\vee s_\alpha - \mathbf k(\alpha)\cdot2 + \mathbf k(\alpha)\,s_\alpha
 \equiv -(\alpha^\vee+\mathbf k(\alpha))
\mod \sum_{w\in W\setminus\{1\}}\mathbf H(w-1).
\]
Likewise, $s_\alpha\cdot(\alpha^\vee)^2=(\alpha^\vee)^2\cdot s_\alpha$,
and $s_\alpha\cdot\xi=\xi\cdot s_\alpha$ for $\xi\in\Liea(\alpha)$.
Now the lemma is clear.
\end{proof}
\begin{cor}\label{cor:ShW}
Under the natural identification 
$S_{\mathbf H}(\Liea_\CC)\simeq S(\Liea_\CC)$ (see \S\ref{sec:intro}),
The space of $W$-fixed elements in $S_{\mathbf H}(\Liea_\CC)$
equals $S(\Liea_\CC)^W$.
\end{cor}
For a quasi-spherical $K$-type $(\sigma,V)$
we shall investigate the map $\Gamma^\sigma$ in \S\ref{sec:intro}.
By virtue of \eqref{eq:K-decomp} we have
\[
\Hom_K(V,U(\Lieg_\CC))=
 \Hom_K(V,\sym(S(\Liep_\CC))) \oplus \Hom_K(V,U(\Lieg_\CC)\Liek_\CC).
\]
Let $\{S^d(\Liep_\CC)\}_{d=0}^\infty$ (resp.~$\{S^d(\Liea_\CC)\}_{d=0}^\infty$) be
the standard grading of $S(\Liep_\CC)$ (resp.~$S(\Liea_\CC)$).
It is easy to see
\begin{equation}\label{eq:SymmEst}
\gamma\circ\sym(F)-\gamma_0(F)\in\sum_{i=0}^{d-1}S^i(\Liea_\CC)
\qquad \forall F\in S^d(\Liep_\CC).
\end{equation}
Therefore, as a corollary of Proposition~\ref{prop:Chinj},
we get the following exact sequence:
\begin{equation}\label{eq:left_half}
 0 \rightarrow \Hom_K(V, U(\Lieg_\CC)\Liek_\CC)
   \rightarrow \Hom_K(V, U(\Lieg_\CC))
   \xrightarrow{\Gamma^\sigma} \Hom_\CC(V^M, S_{\mathbf H}(\Liea_\CC)).
\end{equation}
By the decomposition
$\mathscr A=\bigoplus_{d=0}^\infty S(\Liep_\CC)^K\cap S^d(\Liep_\CC)
=\bigoplus_{d=0}^\infty S(\Liea_\CC)^W\cap S^d(\Liea_\CC)$,
$\mathscr A$ is a graded algebra.
Also, by the decompositions
$\Hom_K(V,S(\Liep_\CC))=\bigoplus_{d=0}^\infty\Hom_K(V,S^d(\Liep_\CC))$ and
$\Hom_W(V^M,S(\Liea_\CC))=\bigoplus_{d=0}^\infty\Hom_W(V^M,S^d(\Liea_\CC))$,
$\Hom_K(V,S(\Liep_\CC))$ and $\Hom_W(V^M,S(\Liea_\CC))$ are graded $\mathscr A$-modules.
Homogeneity of an element of these modules
is defined in the usual way.
\begin{lem}\label{lem:higheronto}
There is a non-zero homogeneous element $b\in\mathscr A$
such that
$b\cdot \Hom_W(V^M,S(\Liea_\CC)) \subset \Gamma^\sigma_0\left(\Hom_K(V,S(\Liep_\CC))\right)$.
\end{lem}
\begin{proof}
Note that $\mathscr A$ is an integral domain.
In view of \eqref{eq:Ast1} and \eqref{eq:Ast2},
both $\Hom_K(V,S(\Liep_\CC))$ and $\Hom_W(V^M,S(\Liea_\CC))$
are free $\mathscr A$-modules of the same rank
and admit bases consisting of homogeneous elements.
By Proposition~\ref{prop:Chinj},
$\Gamma^\sigma_0 : \Hom_K(V,S(\Liep_\CC))\rightarrow \Hom_W(V^M,S(\Liea_\CC))$
is an injective $\mathscr A$-homomorphism.
Moreover, clearly $\Gamma^\sigma_0$ maps a homogeneous element to a homogeneous element of the same degree.
From these facts the lemma follows easily.
\end{proof}
Let the map $p^\sigma : V\rightarrow V^M$ be as in the proof of Proposition~\ref{prop:Chinj}.
\begin{lem}\label{lem:HCprojM}
For any $\Psi\in\Hom_K(V, U(\Lieg_\CC))$ and $\varv\in V$,
$\gamma\left(\Psi[\varv]\right)=\Gamma^\sigma(\Psi)[p^\sigma(\varv)]$.
\end{lem}
\begin{proof}
Since $M$ normalizes $\Liek_\CC$ and $\Lien_\CC$,
$\gamma$ is an $M$-homomorphism from $U(\Lieg_\CC)$
to a trivial $M$-module $S(\Liea_\CC)$.
Hence the map $V\ni \varv\mapsto \gamma\left(\Psi[\varv]\right) \in S(\Liea_\CC)$
is an $M$-homomorphism and the lemma follows.
\end{proof}
The sufficiency of the second statement of Theorem~\ref{thm:HChom}
comes from
\begin{thm}\label{thm:W-inv}
Suppose $\varv\in V_\single^M$.
Then for any $\Psi\in\Hom_K(V, U(\Lieg_\CC))$ and $w\in W$,
\begin{equation}\label{eq:W-inv}
\Gamma^\sigma(\Psi)[w\varv]=w\,\Gamma^\sigma(\Psi)[\varv],
\end{equation}
where the action of $w$ in the right-hand side is that on $S_{\mathbf H}(\Liea_\CC)$.
\end{thm}
\begin{proof}
Suppose $\Psi\in\Hom_K(V, U(\Lieg_\CC))$.
For each $\alpha\in\Pi$,
we define
$\Lieg(\alpha)$,
$\Lieg_\ss(\alpha)$,
$\Liek_\ss(\alpha)$,
$G_\ss(\alpha)$, $K_\ss(\alpha)$, and
$M_\ss(\alpha)$
as in the proof of Lemma~\ref{lem:double_root}.
Recall $\Liea(\alpha)=\{H\in\Liea;\,\alpha(H)=0\}$.
Moreover, put
\[
\Liez(\alpha)=\text{ the center of }\Lieg(\alpha),\quad
\Lien_\alpha =\sum_{\beta\in \Sigma^+\setminus \ZZ\alpha}\Lieg_\beta, \quad
\Liek_\alpha =\sum \{ \RR(X_\beta+\theta X_\beta) ;\,
   \beta\in \Sigma \setminus \ZZ\alpha, X_\beta\in\Lieg_\beta\}.
\]
If we define the projection map
\[
\gamma_\alpha : 
\begin{aligned}
U(\Lieg_\CC)&=
 \bigl((\Lien_\alpha)_\CC U(\Lieg_\CC) 
  + U(\Lieg_\CC) (\Liek_\alpha)_\CC
  + U(\Lieg_\CC) \Liez(\alpha)_\CC \bigr)
\,\oplus\,
 U(\Liea(\alpha)_\CC+\Lieg_\ss(\alpha)_\CC) \\
 &\qquad\longrightarrow U(\Liea(\alpha)_\CC+\Lieg_\ss(\alpha)_\CC),
\end{aligned}
\]
then $\gamma_\alpha$ is a $K_\ss(\alpha)$-homomorphism and $\gamma\circ\gamma_\alpha=\gamma$.
Let $U(\Liek_\ss(\alpha)_\CC)\varv=V^{(1)}\oplus\cdots\oplus V^{(t)}$ be an irreducible
decomposition as a $K_\ss(\alpha)$-module and $\varv=\varv^{(1)}+\cdots+\varv^{(t)}$ the corresponding decomposition.
Then by the same argument as in the proof of Lemma~\ref{lem:double_root},
for each $s=1,\ldots,t$,
$V^{(s)}$ is considered as a single-petaled `$K$-type' of the adjoint group of
$\Lieg_\ss(\alpha)$ and
$\varv^{(s)}$ is a non-zero $M_\ss(\alpha)$-fixed vector of $V^{(s)}$.
Let $\mathscr H_\alpha$ be the set of $K_\ss(\alpha)$-harmonic elements in
$S(\Liep_\CC\cap \Lieg_\ss(\alpha)_\CC)$ and fix an arbitrary
$\Psi^{(s)}\in\Hom_{K_\ss(\alpha)}(V^{(s)},\sym(\mathscr H_\alpha))\setminus\{0\}$.
Put $S_\alpha=\sym(S(\Liep_\CC\cap \Lieg_\ss(\alpha)_\CC)^{K_\ss(\alpha)})$.
Define $\Psi_s \in \Hom_{K_\ss(\alpha)}(V^{(s)},U(\Liea(\alpha)_\CC+\Lieg_\ss(\alpha)_\CC))$
by
\[
V^{(s)}\hookrightarrow V \xrightarrow{\Psi} U(\Lieg_\CC) \xrightarrow{\gamma_\alpha} 
U(\Liea(\alpha)_\CC+\Lieg_\ss(\alpha)_\CC).
\]
Since
\[\begin{aligned}
\Hom_{K_\ss(\alpha)}(V^{(s)},U(\Liea(\alpha)_\CC+\Lieg_\ss(\alpha)_\CC))
 &\simeq \Hom_{K_\ss(\alpha)}(V^{(s)},U(\Lieg_\ss(\alpha)_\CC))\otimes S(\Liea(\alpha)_\CC),\\
\Hom_{K_\ss(\alpha)}(V^{(s)},U(\Lieg_\ss(\alpha)_\CC))
 &\simeq \Hom_{K_\ss(\alpha)}(V^{(s)},U(\Lieg_\ss(\alpha)_\CC)\Liek_\ss(\alpha)_\CC) \\
 &\qquad\qquad  \oplus \Hom_{K_\ss(\alpha)}(V^{(s)},\sym(\mathscr H_\alpha))\otimes S_\alpha,\\
\Hom_{K_\ss(\alpha)}(V^{(s)},\sym(\mathscr H_\alpha))
 &=\CC\Psi^{(s)},
\end{aligned}
\]
we can choose $f_1,\ldots,f_\mu\in S(\Liea(\alpha)_\CC)$
and $D_1,\ldots,D_\mu\in S_\alpha$ so that
$\Psi_s-\Psi^{(s)}(D_1f_1+\cdots+D_\mu f_\mu)
\in\Hom_{K_\ss(\alpha)}(V^{(s)},U(\Liea(\alpha)_\CC+\Lieg_\ss(\alpha)_\CC)\Liek_\ss(\alpha)_\CC)$.
Then by Lemma~\ref{lem:HCprojM},
\[\begin{aligned}
&\Gamma^\sigma(\Psi)[p^\sigma\!\left(\varv^{(s)}\right)]=
\gamma\left(\Psi_s[\varv^{(s)}]\right)=
\gamma\left(\Psi^{(s)}[\varv^{(s)}]\right)\cdot 
(\gamma(D_1)f_1+\cdots+\gamma(D_\mu)f_\mu),\\
&\gamma(D_1)f_1+\cdots+\gamma(D_\mu)f_\mu \in S(\Liea(\alpha)_\CC)\cdot\CC[(\alpha^\vee)^2].
\end{aligned}
\]
Now by Corollary~\ref{cor:rank1-single-petaled},
each $V^{(s)}$ is either the trivial $K_\ss(\alpha)$-type
or a $K_\ss(\alpha)$-type appearing in $(\Ad, \Liep_\CC\cap\Lieg_\ss(\alpha)_\CC)$.
Suppose $V^{(s)}$ is the trivial $K_\ss(\alpha)$-type.
Then it follows from
Lemma~\ref{lem:rank1} (\ref{item:rank1HCimage}), (\ref{item:rank1small})
that $\gamma\left(\Psi^{(s)}[\varv^{(s)}]\right)$ is a scalar.
Hence by Lemma~\ref{lem:Hecke},
$s_\alpha\Gamma^\sigma(\Psi)[p^\sigma\!\left(\varv^{(s)}\right)]
=\Gamma^\sigma(\Psi)[p^\sigma\!\left(\varv^{(s)}\right)]$.
On the other hand, suppose
$(\Ad, \Liep_\CC\cap\Lieg_\ss(\alpha)_\CC)$ is a $K_\ss(\alpha)$-type
appearing in $(\Ad, \Liep_\CC\cap\Lieg_\ss(\alpha)_\CC)$.
This time it follows from Lemma~\ref{lem:rank1} (\ref{item:rank1HCimage}), (\ref{item:rank1small})
that $\gamma\left(\Psi^{(s)}[\varv^{(s)}]\right)$ equals
$\alpha^\vee+\dim \Lieg_\alpha + 2\dim \Lieg_{2\alpha}$
up to a scalar multiple.
Hence Lemma~\ref{lem:Hecke} implies
$s_\alpha\Gamma^\sigma(\Psi)[p^\sigma\!\left(\varv^{(s)}\right)]
=-\Gamma^\sigma(\Psi)[p^\sigma\!\left(\varv^{(s)}\right)]$.
Also, in this case $s_\alpha \varv^{(s)}=-\varv^{(s)}$ by Theorem~\ref{thm:p}.
Thus, in either case we get 
$\Gamma^\sigma(\Psi)[p^\sigma\!\left(s_\alpha \varv^{(s)}\right)]
=s_\alpha\Gamma^\sigma(\Psi)[p^\sigma\!\left(\varv^{(s)}\right)]$
for each $s=1,\ldots,t$.
Hence
$\Gamma^\sigma(\Psi)[s_\alpha \varv]
=\Gamma^\sigma(\Psi)[p^\sigma\!\left(s_\alpha \varv^{(1)}+\cdots+s_\alpha \varv^{(t)}\right)]=
s_\alpha \Gamma^\sigma(\Psi)[\varv]$
for each $\alpha\in\Pi$,
which assures \eqref{eq:W-inv} for any $w\in W$.
\end{proof}
Suppose $\alpha\in\Pi$ and retain the notation of the proof of Theorem~\ref{thm:W-inv}.
We say a $K_\ss(\alpha)$-type $(\sigma',V')$ is quasi-spherical
if $V'$ has a non-zero $M_\ss(\alpha)$-fixed vector.
A quasi-spherical $K_\ss(\alpha)$-type is naturally identified with a
quasi-spherical `$K$-type' of the adjoint group of $\Lieg_\ss(\alpha)$.
Choose $X_\alpha\in\Lieg_\alpha$ so that $B(X_\alpha,\theta X_\alpha)=-\frac1{2|\alpha|^2}$
and put $Z=\iu X_\alpha+\iu\theta X_\alpha$.
According to Lemma~\ref{lem:rank1} (\ref{item:rank1eigen}),
for each quasi-spherical $K_\ss(\alpha)$-type $(\sigma',V')$
we define the integer $e(\sigma')$ as the largest eigenvalue of $\sigma'(Z)$.
If we decompose the $K_\ss(\alpha)$-module $U(\Liek_\ss(\alpha)_\CC)V^M$
into irreducible submodules, then each submodule is quasi-spherical.
In fact, if the corresponding decomposition of a given $\varv \in V^M$
is $\varv=\varv^{(1)}+\varv^{(2)}+\cdots$, then each $\varv^{(s)}$ ($s=1,2,\ldots$)
is an $M_\ss(\alpha)$-fixed vector.
Let us consider the direct sum decomposition
\[
U(\Liek_\ss(\alpha)_\CC)V^M=V_{[0]}\oplus V_{[1]} \oplus\cdots\oplus V_{[k]},
\]
where $V_{[s]}$ ($s=1,\ldots,k$) is the sum of
all irreducible $K_\ss(\alpha)$-submodules which are isomorphic to
some $K_\ss(\alpha)$-type $(\sigma',V')$ with $|e(\sigma')|=s$.
\begin{lem}\label{lem:d_str}
$p^\sigma\!\left(V_{[0]}^{M_\ss(\alpha)}\right)
 \cap
 p^\sigma\!\left(V_{[1]}^{M_\ss(\alpha)}\right) = 0$ and
$p^\sigma\!\left(V_{[2]}^{M_\ss(\alpha)}+\cdots+V_{[k]}^{M_\ss(\alpha)}\right) \subset V^M_\double$.
Moreover,
\begin{equation}\label{eq:d_str}
V^M_\double=\left(
 V^M_\double \cap p^\sigma\!\left(V_{[0]}^{M_\ss(\alpha)}\right)
 \oplus
 V^M_\double \cap p^\sigma\!\left(V_{[1]}^{M_\ss(\alpha)}\right)
\right) + p^\sigma\!\left(V_{[2]}^{M_\ss(\alpha)}+\cdots+V_{[k]}^{M_\ss(\alpha)}\right).
\end{equation}
\end{lem}
\begin{proof}
From Lemma~\ref{lem:rank1} (\ref{item:rank1one}), (\ref{item:rank1intersect})
and the inclusion relation
\[
V^M\subset V_{[0]}^{M_\ss(\alpha)} \oplus V_{[1]}^{M_\ss(\alpha)} 
 \oplus\cdots\oplus V_{[k]}^{M_\ss(\alpha)},
\]
we have $V^M_\single\subset V_{[0]}^{M_\ss(\alpha)} \oplus V_{[1]}^{M_\ss(\alpha)}$.
Let $(\cdot,\cdot)_V$ be a $K$-invariant Hermitian inner product on $V$.
Then the proof of Lemma~\ref{lem:direct_sum} 
says the orthogonal complement $\bigl(V^M_\single\bigr)^\perp$
of $V^M_\single$ in $V^M$ equals $V^M_\double$.
Since
$\bigl(V^M_\single , p^\sigma\!\left(V_{[2]}+\cdots+V_{[k]}\right)\bigr)_V
=\bigl(V^M_\single , V_{[2]}+\cdots+V_{[k]}\bigr)_V
\subset(V_{[0]}+V_{[1]}, V_{[2]}+\cdots+V_{[k]}\bigr)_V=0$,
$p^\sigma\!\left(V_{[2]}+\cdots+V_{[k]}\right) \subset \bigl(V^M_\single\bigr)^\perp=V^M_\double$.
On the other hand, it follows from Lemma~\ref{lem:rank1} (\ref{item:rank1small}) and Theorem~\ref{thm:p}
that $s_\alpha$ acts on $V_{[0]}^{M_\ss(\alpha)}$ and $V_{[1]}^{M_\ss(\alpha)}$
by $+1$ and  $-1$, respectively.
But since $p^\sigma : V \rightarrow V^M$ is an $N_K(\Liea)$-homomorphism,
$s_\alpha$ acts on $p^\sigma\!\left(V_{[0]}^{M_\ss(\alpha)}\right)$
and $p^\sigma\!\left(V_{[1]}^{M_\ss(\alpha)}\right)$ by $+1$ and $-1$, respectively.
Hence the decomposition of $V^M_\double \cap p^\sigma\!\left(V_{[0]}+V_{[1]}\right)
 =V^M_\double \cap p^\sigma\!\left(V_{[0]}^{M_\ss(\alpha)}+V_{[1]}^{M_\ss(\alpha)}\right)$
into the eigenspaces of $s_\alpha$ is
$V^M_\double \cap p^\sigma\!\left(V_{[0]}^{M_\ss(\alpha)}\right) \oplus
 V^M_\double \cap p^\sigma\!\left(V_{[1]}^{M_\ss(\alpha)}\right)$.
We thus get
\[\begin{aligned}
V^M_\double&=V^M_\double \cap p^\sigma\!\left(V_{[0]}+\cdots+V_{[k]}\right) \\
 &= V^M_\double \cap p^\sigma\!\left(V_{[0]}+V_{[1]}\right)
   + p^\sigma\!\left(V_{[2]}+\cdots+V_{[k]}\right) \\
 &= \left(
 V^M_\double \cap p^\sigma\!\left(V_{[0]}^{M_\ss(\alpha)}\right)
 \oplus
 V^M_\double \cap p^\sigma\!\left(V_{[1]}^{M_\ss(\alpha)}\right)
\right) + p^\sigma\!\left(V_{[2]}^{M_\ss(\alpha)}+\cdots+V_{[k]}^{M_\ss(\alpha)}\right).\qedhere
\end{aligned}\]
\end{proof}
The necessity of the second statement of Theorem~\ref{thm:HChom}
follows from the next proposition:
\begin{prop}\label{prop:W-inv2}
Suppose $V^M_\double\ne0$.
Then \eqref{eq:W-inv} does not hold
for a suitable combination of $\varv\in V^M_\double$, $w\in W$ and $\Psi\in\Hom_K(V,U(\Lieg_\CC))$.
\end{prop}
\begin{proof}
By the assumption of the proposition, there exists $\alpha\in\Pi$ for which
$V_{[2]}+\cdots+V_{[k]}$ in the above argument is not $0$.
Take $s$ ($2\le s\le k$) so that $V_{[s]}\ne 0$.
Fix  an irreducible $K_\ss(\alpha)$-submodule $V_0$ of $V_{[s]}$ and
$\varv_0 \in V_0^{M_\ss(\alpha)}\setminus\{0\}$.

First, we shall show $p^\sigma(\varv_0)\ne0$.
For this, let $U(\Liek_\ss(\alpha))V^M=V^{(1)}\oplus V^{(2)}\oplus \cdots$
be an irreducible decomposition as a $K_\ss(\alpha)$-module such that
$V^{(1)}=V_0$ and each component is orthogonal to the other components
with respect to a $K$-invariant Hermitian inner product $(\cdot,\cdot)_V$ on $V$.
If $\varv=\varv^{(1)}+\varv^{(2)}+\cdots$ is the corresponding decomposition of any $\varv\in V^M$,
then each $\varv^{(s)}$ ($s=1,2,\ldots$) is an $M_\ss(\alpha)$-fixed vector.
Since $V^M$ generates $U(\Liek_\ss(\alpha)_\CC)V^M$,
there exists $\varv\in V^M$ such that $\varv^{(1)}\ne 0$.
Since $\dim_\CC V_0^{M_\ss(\alpha)}=1$,
$\varv^{(1)}=c\varv_0$ for some constant $c\ne0$.
Now $(\varv,p^\sigma(\varv_0))_V=(\varv,\varv_0)_V=(\varv^{(1)},\varv_0)_V=c(\varv_0, \varv_0)_V\ne 0$,
which proves $p^\sigma(\varv_0)\ne0$.

Choose a homogeneous $\varphi \in \Hom_W(V^M,S(\Liea_\CC))$ so that $\varphi[p^\sigma(\varv_0)]\ne 0$.
Let $b\in\mathscr A$ be the homogeneous element in Lemma~\ref{lem:higheronto}.
Then $\Gamma^\sigma_0(\Phi)=b\cdot\varphi$ for some $\Phi \in \Hom_K(V,S(\Liep_\CC))$.
Put $\Psi=\sym\circ\Phi \in \Hom_K(V,U(\Lieg_\CC))$.
Since $\Phi$ is homogeneous and 
$\gamma_0(\Phi[\varv_0])=b\cdot\varphi[p^\sigma(\varv_0)]\ne 0$,
\eqref{eq:SymmEst} implies $\gamma\left(\Psi[\varv_0]\right)\ne0$.

Let $\gamma_\alpha$, $S_\alpha$, $\mathscr H_\alpha$ be the same as in the proof of Theorem~\ref{thm:W-inv}.
Fix an arbitrary $\Psi^0\in \Hom_{K_\ss(\alpha)}(V_0, \sym(\mathscr H_\alpha))\setminus\{0\}$ and
define $\Psi_0 \in \Hom_{K_\ss(\alpha)}(V_0,U(\Liea(\alpha)_\CC+\Lieg_\ss(\alpha)_\CC))$
by
\[
V_0 \hookrightarrow V \xrightarrow{\Psi} U(\Lieg_\CC) \xrightarrow{\gamma_\alpha} 
U(\Liea(\alpha)_\CC+\Lieg_\ss(\alpha)_\CC).
\]
Then, as in the proof of Theorem~\ref{thm:W-inv},
we can choose $f_1,\ldots,f_\mu\in S(\Liea(\alpha)_\CC)$
and $D_1,\ldots,D_\mu\in S_\alpha$ so that
$\Psi_0-\Psi^0(D_1f_1+\cdots+D_\mu f_\mu)
\in\Hom_{K_\ss(\alpha)}(V_0,U(\Liea(\alpha)_\CC+\Lieg_\ss(\alpha)_\CC)\Liek_\ss(\alpha)_\CC)$.
Hence
\[\begin{aligned}
&\gamma\left(\Psi[\varv_0]\right)=
\gamma\left(\Psi_0[\varv_0]\right)=
\gamma\left(\Psi^0[\varv_0]\right)\cdot 
(\gamma(D_1)f_1+\cdots+\gamma(D_\mu)f_\mu),\\
&\gamma(D_1)f_1+\cdots+\gamma(D_\mu)f_\mu \in S(\Liea(\alpha)_\CC)\cdot\CC[(\alpha^\vee)^2].
\end{aligned}
\]
From Lemma~\ref{lem:rank1} (\ref{item:rank1HCimage}),
there is a pair $(i, j)$ of non-negative integers with $2i+j=s\ (\ge 2)$
such that $\gamma\left(\Psi^0[\varv_0]\right)$ equals
\begin{equation}\label{eq:rank1HCimage2}
\bigl[ (h+\delta)(h+\delta+2)\cdots(h+\delta+2(i+j)-2) \bigr]
 \cdot \bigl[ (h+1)(h+3)\cdots(h+2i-1) \bigr]
\end{equation}
up to a scalar multiple (here $h=\frac{\alpha^\vee}2+\frac{\dim\Lieg_\alpha}2\in S(\Liea_\CC)$ and
$ \delta=\dim\Lieg_{2\alpha}$).
Since it is clear that $\gamma\left(\Psi^0[\varv_0]\right), 
 \frac1{(h+\delta)}\gamma\left(\Psi^0[\varv_0]\right)
 \notin
 \CC[(\alpha^\vee)^2]$, 
we have $\gamma\left(\Psi[\varv_0]\right) = \Gamma^\sigma(\Psi)[p^\sigma(\varv_0)] 
  \,\notin\, 
 S(\Liea(\alpha)_\CC)\cdot \CC[(\alpha^\vee)^2] \,\cup\,
 S(\Liea(\alpha)_\CC)\cdot \CC[(\alpha^\vee)^2](\alpha^\vee+\dim \Lieg_\alpha + 2\dim \Lieg_{2\alpha})$.
Hence from Lemma~\ref{lem:Hecke},
$s_\alpha \Gamma^\sigma(\Psi)[p^\sigma(\varv_0)] \ne \pm \Gamma^\sigma(\Psi)[p^\sigma(\varv_0)]$.
On the other hand, since $\dim_\CC V_0^{M_\ss(\alpha)}=1$,
we have $s_\alpha \varv_0= \pm \varv_0$ and therefore $s_\alpha p^\sigma(\varv_0)= \pm p^\sigma(\varv_0)$.
\end{proof}
\begin{rem}\label{rem:alpha_sign}
In the above proof we can say $s_\alpha p^\sigma(\varv_0)= (-1)^{s} p^\sigma(\varv_0)$.
Indeed, if we choose $\Phi^0\in\Hom_{K_\ss(\alpha)}(V_0,\mathscr H_\alpha)$
so that $\Psi^0=\sym\circ\Phi^0$,
then it is a homogeneous element by Lemma~\ref{lem:rank1} (\ref{item:rank1one}).
Hence from \eqref{eq:SymmEst} and \eqref{eq:rank1HCimage2},
$\gamma_0\left(\Phi^0[\varv_0]\right)=C\left(\frac{\alpha^\vee}2\right)^s$
for some non-zero constant $C$.
With respect to the ordinary action of $s_\alpha$ on $\CC[\alpha^\vee]$,
\[
 \gamma_0\left(\Phi^0[s_\alpha \varv_0]\right)
 =s_\alpha \gamma_0\left(\Phi^0[\varv_0]\right)
 =s_\alpha C\left(\frac{\alpha^\vee}2\right)^s
 =(-1)^s C\left(\frac{\alpha^\vee}2\right)^s
 =\gamma_0\left(\Phi^0[(-1)^s \varv_0]\right).
\]
Thus we have $s_\alpha \varv_0= (-1)^{s} \varv_0$, which shows our claim.
\end{rem}
In the rest of this section
we shall prove the following theorem,
which is considered as a non-commutative counterpart of Theorem~\ref{thm:Chsurj}.
\begin{thm}\label{thm:HCsurj}
For any $\psi\in\Hom_W\bigl(V^M/V_\double^M, S_{\mathbf H}(\Liea_\CC)\bigr)$
there exists $\Psi\in\Hom_K(V, U(\Lieg_\CC))$ such that $\Gamma^\sigma(\Psi)=\psi$.
\end{thm}
Suppose $V'\subset V^M$ is an arbitrary $W$-submodule.
Let $\iota_{V'}^\#$ denote the map
\[
 \Hom_\CC(V^M,S_{\mathbf H}(\Liea_\CC)) \ni \psi 
   \mapsto \psi|_{V'} \in \Hom_\CC({V'},S_{\mathbf H}(\Liea_\CC))
\]
or the map
\[
 \Hom_W(V^M,S(\Liea_\CC)) \ni \varphi 
    \mapsto \varphi|_{V'} \in \Hom_W({V'},S(\Liea_\CC)).
\]
Under the natural identification $S_{\mathbf H}(\Liea_\CC)\simeq S(\Liea_\CC)$,
put $S_{\mathbf H}^d(\Liea_\CC)=\sum_{i=0}^d S^i(\Liea_\CC)$ ($d=0,1,\ldots$).
Then $\Hom_\CC({V'},S_{\mathbf H}(\Liea_\CC))$ has
the natural filtration $\bigl\{\Hom_\CC({V'},S_{\mathbf H}^d(\Liea_\CC))\bigr\}_{d=0}^\infty$
by which it is considered as a filtered $\mathscr A$-module.
For each $\psi \in \Hom_\CC({V'},S_{\mathbf H}(\Liea_\CC))$ put
\[
\deg \psi = \begin{cases}
\ d & \text{if }\psi \in \Hom_\CC({V'},S_{\mathbf H}^d(\Liea_\CC)) \setminus
\Hom_\CC({V'},S_{\mathbf H}^{d-1}(\Liea_\CC)), \\
\ -\infty & \text{if }\psi=0.
\end{cases}
\]
Also, for $d=-\infty,0,1,\ldots$, define the natural map
\[
q_d : \Hom_\CC({V'},S_{\mathbf H}^d(\Liea_\CC)) \longrightarrow 
 \Hom_\CC({V'}, S^d(\Liea_\CC)).
\]
\begin{lem}\label{lem:single-W-inv}
For any $\Psi\in\Hom_K(V,U(\Lieg_\CC))$,
$\iota_{V^M_\single}^{\#} \!\!\circ\Gamma^\pi(\Psi) \in
\Hom_W(V^M_\single,S_{\mathbf H}(\Liea_\CC))$.
For any $\psi\in \Hom_W(V^M_\single,S_{\mathbf H}(\Liea_\CC))$
with $d=\deg\psi$,
$q_d(\psi) \in \Hom_W(V^M_\single, S^d(\Liea_\CC))$.
\end{lem}
\begin{proof}
The first assertion is due to Theorem~\ref{thm:W-inv}.
The second assertion follows from the fact that the map
$S_{\mathbf H}^d(\Liea_\CC)=\sum_{i=0}^d S^i(\Liea_\CC) \xrightarrow{\text{projection}} S^d(\Liea_\CC)$
is a $W$-homomorphism, which is easily checked by use of \eqref{eq:Hrel} and \eqref{eq:Sh}.
\end{proof}
\begin{lem}\label{lem:double-W-inv}
For any $\Psi\in\Hom_K(V,U(\Lieg_\CC))$, put
$\psi\coloneqq\iota_{V^M_\double}^{\#} \!\!\circ\Gamma^\pi(\Psi) \in \Hom_\CC(V^M_\double,S_{\mathbf H}(\Liea_\CC))$
and $d=\deg\psi$. Then $q_d(\psi) \in \Hom_W(V^M_\double, S^d(\Liea_\CC))$.
\end{lem}
\begin{proof}
We shall check for any $\alpha \in \Pi$ and $\varv\in V^M_\double$,
\begin{equation}\label{eq:sa_kakan}
q_d(\psi)[s_\alpha \varv]=s_\alpha q_d(\psi)[\varv].
\end{equation}
For this, fix $\alpha\in\Pi$ and consider the decomposition \eqref{eq:d_str}.

Take an arbitrary irreducible $K_\ss(\alpha)$-submodule
$V_0\subset V_{[s]}$ ($s=2,3,\ldots$) and $\varv_0 \in V_0^{M_\ss(\alpha)}$.
Then for the same reason as the proofs of Theorem~\ref{thm:W-inv} and Proposition~\ref{prop:W-inv2},
there are some non-negative integers $i$,$j$ with $2i+j=s$ and some
$D_0 \in S(\Liea(\alpha)_\CC)\cdot\CC[(\alpha^\vee)^2]$ such that
\[
\Gamma^\sigma(\Psi)[p^\sigma(\varv_0)]=
D_0 \cdot \bigl[ (h+\delta)(h+\delta+2)\cdots(h+\delta+2(i+j)-2) \bigr]
 \cdot \bigl[ (h+1)(h+3)\cdots(h+2i-1) \bigr].
\]
(Recall $h=\frac{\alpha^\vee}2+\frac{\dim\Lieg_\alpha}2\in S(\Liea_\CC)$ and
$ \delta=\dim\Lieg_{2\alpha}$.)
Therefore there exists a homogeneous element $\bar{D}_0$ in 
$S(\Liea(\alpha)_\CC)\cdot\CC[(\alpha^\vee)^2]$
such that
\[
q_d(\psi)[p^\sigma(\varv_0)]=\bar D_0\left(\frac{\alpha^\vee}2\right)^s.
\]
Hence $s_\alpha q_d(\psi)[p^\sigma(\varv_0)]=(-1)^s q_d(\psi)[p^\sigma(\varv_0)]$.
But since $s_\alpha p^\sigma(\varv_0)=(-1)^s p^\sigma(\varv_0)$ by Remark~\ref{rem:alpha_sign},
\eqref{eq:sa_kakan} is valid for $\varv=p^\sigma(\varv_0)$.

Secondly, let $V^{(1)} \oplus \cdots \oplus V^{(t)}$ be an irreducible decomposition
of the $K_\ss(\alpha)$-module $V_{[1]}$.
Suppose $\varv_{[1]} \in V_{[1]}^{M_\ss(\alpha)}$ satisfies
$p^\sigma(\varv_{[1]}) \in V^M_\double$ and
let $\varv_{[1]}=\varv^{(1)}+\cdots+\varv^{(t)}$
be the decomposition according to the above decomposition.
Since $\varv^{(s)}\in \left(V^{(s)}\right)^{M_\ss(\alpha)}$ ($s=1,\ldots,t$), 
as in the proof of Theorem~\ref{thm:W-inv}, there exist
$D_1,\ldots,D_t\in S(\Liea(\alpha)_\CC)\cdot\CC[(\alpha^\vee)^2]$ such that
\[
\Gamma^\pi(\Psi)[p^\sigma\!\left(\varv^{(s)}\right)]=D_s (\alpha^\vee+\dim \Lieg_\alpha + 2\dim \Lieg_{2\alpha})
\qquad s=1,\ldots,t.
\]
Hence there exists a homogeneous element $\bar D$ in
$S(\Liea(\alpha)_\CC)\cdot\CC[(\alpha^\vee)^2]$ such that
\[
q_d(\psi)[p^\sigma(\varv_{[1]})]=\bar D\, \alpha^\vee.
\]
Since $s_\alpha p^\sigma(\varv_{[1]})=-p^\sigma(\varv_{[1]})$ and
$s_\alpha \bar D\, \alpha^\vee = -\bar D\, \alpha^\vee$,
\eqref{eq:sa_kakan} is valid for $\varv=p^\sigma(\varv_{[1]})$.

Similarly, we can show \eqref{eq:sa_kakan}
for any $\varv \in V^M_\double \cap p^\sigma\!\left(V_{[0]}^{M_\ss(\alpha)}\right)$.
Hence from \eqref{eq:d_str}, \eqref{eq:sa_kakan} is valid for any $\varv \in V^M_\double$.
\end{proof}
\begin{lem}\label{lem:single_appro}
For any $\psi \in \Hom_W(V^M_\single, S_{\mathbf H}(\Liea_\CC))\setminus\{0\}$
there exists $\Psi \in \Hom_K(V, U(\Lieg_\CC))$ such that
\[
\iota^\#_{V^M_\single} \!\!\circ\Gamma^\sigma(\Psi) = 
 \psi,\qquad
 \deg \iota^\#_{V^M_\double} \!\!\circ\Gamma^\sigma(\Psi)
 < \deg \psi.
\]
\end{lem}
\begin{proof}
Put $d=\deg \psi$.
Assume that for some $i\in\{d+1,d,d-1,\ldots,1\}$
we already have
$\Psi_i \in \Hom_K(V, U(\Lieg_\CC))$ such that
\[
\deg \left(\psi-\iota^\#_{V^M_\single} \!\!\circ\Gamma^\sigma(\Psi_i)\right) < i,\qquad
 \deg \iota^\#_{V^M_\double} \!\!\circ\Gamma^\sigma(\Psi_i) < d.
\]
Then from Lemma~\ref{lem:single-W-inv}, we get
$\varphi_{i-1} \coloneqq q_{i-1}\left(\psi-\iota^\#_{V^M_\single} \!\!\circ\Gamma^\sigma(\Psi_i)\right)
\in \Hom_W(V^M_\single,S^{i-1}(\Liea_\CC))$.
Since $V^M_\single\simeq V^M/V^M_\double$,
we identify $\varphi_{i-1}$ with an element of $\Hom_W(V^M/V^M_\double,S^{i-1}(\Liea_\CC))$.
Then by Theorem~\ref{thm:Chsurj},
there exists a unique $\Phi_{i-1} \in \Hom_K(V, S^{i-1}(\Liep_\CC))$ such that
$\Gamma^\sigma_0(\Phi_{i-1})=\varphi_{i-1}$.
In view of \eqref{eq:SymmEst} we see
$\deg \Gamma^\sigma(\sym\circ\,\Phi_{i-1}) \le i-1$ and
\[
\begin{aligned}
&q_{i-1}\circ \iota^\#_{V^M_\single} \!\!\circ\Gamma^\sigma(\sym\circ\,\Phi_{i-1})
 =\iota^\#_{V^M_\single} (\varphi_{i-1}) = q_{i-1}\left(\psi-\iota^\#_{V^M_\single} \!\!\circ\Gamma^\sigma(\Psi_i)\right),\\
&q_{i-1}\circ \iota^\#_{V^M_\double} \!\!\circ\Gamma^\sigma(\sym\circ\,\Phi_{i-1})
 =\iota^\#_{V^M_\double} (\varphi_{i-1}) = 0.
\end{aligned}
\]
Hence, $\deg \iota^\#_{V^M_\double} \!\!\circ\Gamma^\sigma(\sym\circ\,\Phi_{i-1})<i-1$
and if we put $\Psi_{i-1}=\Psi_i+\sym\circ\,\Phi_{i-1}$,
then $\deg \left(\psi-\iota^\#_{V^M_\single} \!\!\circ\Gamma^\sigma(\Psi_{i-1})\right) < i-1$.

Thus, if we start with $\Psi_{d+1}\coloneqq0$
and define $\Psi_d,\Psi_{d-1},\ldots$ as above,
then $\Psi\coloneqq\Psi_0$ satisfies the desired properties.
\end{proof}
\begin{proof}[Proof of Theorem~\ref{thm:HCsurj}]
Put $m=\dim_\CC V^M$ and $m'=\dim_\CC V^M_\single$.
Take a basis $\bigl\{\varphi_{m'+1},\ldots, \varphi_{m}\bigr\}$ of
$\Hom_W(V^M/V^M_\single, \mathscr H_W(\Liea))$ so that
each $\varphi_i$ is homogeneous
(note that $\mathscr H_W(\Liea)\simeq\CC[W]$).
Let $b\in\mathscr A$ be the homogeneous element in Lemma~\ref{lem:higheronto}.
Then there exist $\Phi_{m'+1},\ldots,\Phi_{m} \in \Hom_K(V,S(\Liep_\CC))$ such that
$\Gamma^\sigma(\Phi_i)=b\cdot\varphi_i$.
Put $d_i=\deg\varphi_i$ ($i=m'+1,\ldots,m$) and $d_0=\deg b$.
Owing to \eqref{eq:SymmEst}, Lemma~\ref{lem:single-W-inv}, and Lemma~\ref{lem:single_appro},
by modifying $\sym\circ\Phi_{m'+1},\ldots,\sym\circ\Phi_{m}$
in lower-order terms,
we can get
$\Psi_{m'+1},\ldots,\Psi_{m} \in \Hom_K(V,U(\Lieg_\CC))$
which satisfy for each $i=m'+1,\ldots,m$,
\begin{equation}\label{eq:double_basis}
 \deg \Gamma^\sigma(\Psi_i) = \deg \Phi_i = d_i + d_0, \quad
 q_{d_i+d_0}\circ\iota^\#_{V^M_\double}\!\!\circ\Gamma^\sigma(\Psi_i)
  = b\cdot \iota^\#_{V^M_\double}(\varphi_i), \quad
 \iota^\#_{V^M_\single}\!\!\circ\Gamma^\sigma(\Psi_i)=0.
\end{equation}

Put $\mathscr M=\iota^\#_{V^M_\double}\!\!\circ\Gamma^\sigma
 \left(\Hom_K(V,U(\Lieg_\CC))\right)$.
Then by \eqref{eq:Uk-inv},
it is a submodule of the filtered $\mathscr A$-module $\Hom_\CC(V^M_\double, S_{\mathbf H}(\Liea_\CC))$.
Also from Lemma~\ref{lem:double-W-inv},
$\gr \mathscr M \subset \Hom_W(V^M_\double, S(\Liea_\CC))$.
Since $\gr \mathscr M$ is finitely generated over $\mathscr A$,
we can take $\tilde\Psi_1,\ldots,\tilde\Psi_k \in \Hom_K(V,U(\Lieg_\CC))$ so that
$\bigl\{ q_{\tilde d_i}\circ\iota^\#_{V^M_\double}\!\!\circ\Gamma^\sigma(\tilde\Psi_i);\, i=1,\ldots k\bigr\}$
generates $\gr \mathscr M$ over $\mathscr A$
(here $\tilde d_i\coloneqq\deg \iota^\#_{V^M_\double}\!\!\circ\Gamma^\sigma(\tilde\Psi_i)$).
Now, from \eqref{eq:K-decomp} we have
$U(\Lieg_\CC)^K = \sym(S(\Liep_\CC)^K) \oplus U(\Lieg_\CC)^K\cap U(\Lieg_\CC)\Liek_\CC$.
Hence by the exactness of \eqref{eq:HChom},
$\gamma$ gives the isomorphism $\sym(S(\Liep_\CC)^K) \simarrow S(\Liea_\CC)^W$.
For each $a \in \mathscr A=S(\Liea_\CC)^W$,
we denote by $\Hat a$
the unique element of $\sym\bigl(S(\Liep_\CC)^K\bigr)$
such that $\gamma(\Hat{a})=a$.
Then for any $\Psi \in \Hom_K(V,U(\Lieg_\CC))$,
there exist $a_1,\ldots,a_k\in\mathscr A$ with
$\deg a_i \le \deg \iota^\#_{V^M_\double}\!\!\circ\Gamma^\sigma(\Psi)-\tilde d_i$
such that
$\iota^\#_{V^M_\double}\!\!\circ\Gamma^\sigma(\Psi-\tilde\Psi_1 \hat a_1-\cdots-\tilde\Psi_k \hat a_k)=0$.
Now since $\bigl\{\iota^\#_{V^M_\double}(\varphi_i);\,i=m'+1,\ldots,m\bigr\}$ is
a basis of $\Hom_W(V^M_\double, S(\Liea_\CC))$ over $\mathscr A$,
we can take homogeneous elements $b_{is}\in\mathscr A$ ($i=1,\ldots,k,\ s=m'+1,\ldots,m$)
so that
\[
q_{\tilde d_i}\circ\iota^\#_{V^M_\double}\!\!\circ\Gamma^\sigma(\tilde\Psi_i)
 = \sum_{s=m'+1}^{m} b_{is}\cdot \iota^\#_{V^M_\double}(\varphi_s),\qquad
\tilde d_i = \deg b_{is} + d_s\text{ or }b_{is}=0.
\]
Put $\Tilde{\Tilde\Psi}_i \coloneqq \tilde \Psi_i\,\Hat b  - \sum_{s=m'+1}^{m} \Psi_s\,\Hat b_{is}$.
Then from \eqref{eq:double_basis},
$\deg \iota^\#_{V^M_\double}\!\!\circ\Gamma^\sigma(\Tilde{\Tilde\Psi}_i)<\deg b_0+\tilde d_i$.
Hence there exist $a_{ij}\in\mathscr A$ ($i,j=1,\ldots,k$) with $\deg a_{ij} < \deg d_0+\tilde d_i-\tilde d_j$
such that
\[
\iota^\#_{V^M_\double}\!\!\circ\Gamma^\sigma\left(
\Tilde{\Tilde\Psi}_i-\sum_{j=1}^k \tilde\Psi_j\,\hat a_{ij}
\right)=0.
\]
Let us define the $\mathscr A$-valued $k \times k$-matrix
$A\coloneqq \diag(b,\ldots,b)-(a_{ij})_{1\leq i,j\leq k}$.
By estimating the degree of each coefficient of $A$,
we can easily see $\det A \ne 0$.
Let $\tilde A=(\tilde a_{ij})_{1\leq i,j\leq k}$ be the cofactor matrix of $A$.
Observe that \begin{equation}\label{eq:base_change}
\iota^\#_{V^M_\double}\!\!\circ\Gamma^\sigma\left(
 \tilde\Psi_i\cdot\widehat{\det A}
 - \sum_{j=1}^k\sum_{s=m'+1}^m  \Psi_s\,\Hat{\Tilde a}_{ij}\,\hat b_{js}
\right)=0
\qquad
i=1,\ldots,k.
\end{equation}

Now let $\psi\in\Hom_W(V^M/V^M_\double,S_{\mathbf H}(\Liea_\CC))$.
By Lemma~\ref{lem:single_appro}, there is $\Psi \in\Hom_K(V,U(\Lieg_\CC))$ such that
$\iota^\#_{V^M_\single}\!\!\circ\Gamma^\sigma(\Psi)=\iota^\#_{V^M_\single}(\psi)$.
Then we can take $a_1,\ldots,a_k\in\mathscr A$ so that
$\iota^\#_{V^M_\double}\!\!\circ\Gamma^\sigma(\Psi-\tilde\Psi_1 \hat a_1-\cdots-\tilde\Psi_k \hat a_k)=0$.
Hence if we put
\[
\Psi'=\Psi \cdot \widehat{\det A} - \sum_{i=1}^k
\sum_{j=1}^k\sum_{s=m'+1}^m \Psi_s\, \hat a_i\,\Hat{\Tilde a}_{ij}\,\hat b_{js},
\]
then from \eqref{eq:double_basis} and \eqref{eq:base_change} we have
$\iota^\#_{V^M_\single}\!\!\circ\Gamma^\sigma(\Psi')=\det A\cdot \iota^\#_{V^M_\single}(\psi)$
and $\iota^\#_{V^M_\double}\!\!\circ\Gamma^\sigma(\Psi')=0$, 
namely,
$\Gamma^\sigma(\Psi')=\det A\cdot\psi$.
Hence if we put
$\mathscr I=\{c\in\mathscr A;\, c\cdot\psi \in \Gamma^\sigma\left(\Hom_K(V,U(\Lieg_\CC))\right) \}$,
then $\mathscr I\ni\det a\ne0$.
Note that $\mathscr I$ is an ideal of $\mathscr A$.

In order to complete the proof, it suffices to show $\mathscr I=\mathscr A$.
Assume $\mathscr I \subsetneq \mathscr A$.
From inside $\mathscr A\setminus\{0\}$, take an element $c$ 
so that it has the lowest degree.
Then by assumption, $c$ is not a constant.
Let $\Psi''\in\Hom_K(V,U(\Lieg_\CC))$ be such that
\[
\Gamma^\sigma(\Psi'')=c\cdot\psi.
\]
With respect to a basis $\bigl\{\bar\Psi_1,\ldots,\bar\Psi_m\bigr\}$
of $\Hom_K(V,\sym(\mathscr H_K(\Liep)))$ and a basis
$\bigl\{\varv_1,\ldots,\varv_m\bigr\}$ of $V^M$, we define
the matrix $P^\sigma=(\gamma\circ\bar\Psi_j[\varv_i])_{1\leq i,j\leq m}$ as in \S\ref{sec:Kostant}.
By virtue of Corollary~\ref{cor:UgDecomp} and the exactness of \eqref{eq:left_half},
we can take $e_1,\ldots,e_m\in\mathscr A$ so that
$\Psi''-\bar\Psi_1 \hat e_1-\cdots-\bar\Psi_m\hat  e_m
\in \Hom_K(V,U(\Lieg_\CC)\Liek_\CC)$.
Then we have
$\Gamma^\sigma(\bar\Psi_1)\,e_1+\cdots+\Gamma^\sigma(\bar\Psi_m)\,e_m
 = c\cdot\psi$, which, using $P^\sigma$, are rewritten as
\begin{equation}\label{eq:mat_eq}
P^\sigma\begin{pmatrix}
e_1 \\
\vdots \\
e_m
\end{pmatrix}
=c\begin{pmatrix}
\psi[\varv_1] \\
\vdots \\
\psi[\varv_m]
\end{pmatrix}.
\end{equation}

We assert that if $\lambda\in\Liea_\CC^*$ satisfies $c(\lambda)=0$,
then $e_1(\lambda)=\cdots=e_m(\lambda)=0$.
To show this, suppose $\lambda\in\Liea_\CC^*$ satisfies $c(\lambda)=0$.
Then there exists $w\in W$ such that
$\Real \ang{w\lambda}{\alpha}\ge 0$ for any $\alpha\in\Sigma^+$.
Since $c \in \mathscr A=S(\Liea_\CC)^W$,
$c(\lambda)=0$ implies $c(w\lambda)=0$.
Evaluating both sides of \eqref{eq:mat_eq} at $w\lambda$, we have
\[
P^\sigma(w\lambda)\begin{pmatrix}
e_1(w\lambda) \\
\vdots \\
e_m(w\lambda)
\end{pmatrix}
=\begin{pmatrix}
0 \\
\vdots \\
0
\end{pmatrix}.
\]
Then $e_1(w\lambda)=\cdots= e_m(w\lambda)=0$
since $P^\sigma(w\lambda)$ is a regular matrix by Proposition~\ref{prop:Kdet}.
Because $e_1,\ldots,e_m\in \mathscr A=S(\Liea_\CC)^W$,
$e_1(\lambda)=\cdots=e_m(\lambda)=0$. Thus we get the assertion.

Now, $\mathscr A$ is isomorphic to a polynomial ring
(\cite[Ch.~III, Theorem~3.1]{He2}) and
a maximal ideal of $\mathscr A$ equals $\{f\in\mathscr A;\,f(\lambda)=0\}$
for some $\lambda\in\Liea_\CC^*$ (\ibid[Ch.~III, Lemma~3.11]).
Hence by the fact shown above,
$e_j$ ($j=1,\ldots,m$) are divisible by any irreducible factor $c_0$ of $c$.
Let $c', e'_1,\ldots,e'_m\in\mathscr A$ be such that
$c=c'c_0$, $e_j=e'_j\,c_0$ and put
$\Psi'''=\bar\Psi_1 \hat e'_1+\cdots+\bar\Psi_m\hat  e'_m$.
Then we have
\[
\Gamma^\sigma(\Psi''')=c'\cdot\psi,\qquad\deg c'<\deg c,\qquad c'\ne0.
\]
It contradicts the minimality of the degree of $c$.
Thus we get $\mathscr I=\mathscr A$.
\end{proof}
\section{Complex semisimple Lie algebras}\label{sec:complex}
Suppose $\Lieg$ is a complex semisimple Lie algebra with complex structure $J$.
In this case one has $G=G_{\ad}=G_\theta$.
Throughout this section
we use the symbols $U$, $\Lieu$, $\Lieh_\RR$, and $\Lieh$ in place of
$K$, $\Liek$, $\Liea$, and $\Liem+\Liea$, respectively.
Then $\Lieh$ is a Cartan subalgebra.
Extend each $\alpha\in\Sigma$ to a complex linear form on $\Lieh$ and
put $\Tilde\rho=\frac12\sum_{\alpha\in\Sigma^+}\alpha\in\Lieh^*$.
By the unitary trick, a $U$-type $(\sigma,V)$ is naturally identified with
a finite-dimensional irreducible holomorphic representation of $G$.
Since $M=\exp (J\Lieh_\RR)$, $V^M$ equals $V^\Lieh$,
the space of $0$-weight vectors.
Hence in this section,
we denote $V^M_\single$ and $V^M_\double$ by $V^\Lieh_\single$ and $V^\Lieh_\double$, respectively.
One knows each finite-dimensional irreducible holomorphic representation $(\sigma, V)$  of $G$
satisfies $V^\Lieh\ne0$.
This means all $U$-types are quasi-spherical.
From now on we always assume a representation of $G$ is holomorphic. 
\begin{prop}\label{prop:comp_subs}
For any finite-dimensional irreducible representation $(\sigma, V)$ of $G$,
\begin{align}
V^\Lieh_\single &= \Bigl\{\, \varv\in V^\Lieh;\, \sigma(X_\alpha)^2\varv=0\quad \forall \alpha\in\Sigma,\
 \forall X_\alpha\in\Lieg_\alpha  \,\Bigr\},\label{eq:new_single}\\
V^\Lieh_\double &=V^\Lieh \cap \sum\Bigl\{\, \sigma(X_\alpha)^2 V;\, \alpha\in\Sigma,\
 X_\alpha\in\Lieg_\alpha \,\Bigr\}.\label{eq:new_double}
\end{align}
\end{prop}
\begin{proof}
Suppose $\varv\in V^\Lieh$ and $\alpha\in\Sigma$.
Choose $X_\alpha\in\Lieg_\alpha$ so that $B(X_\alpha,\theta X_\alpha)=-\frac1{2|\alpha|^2}$.
Put $Z=\iu X_\alpha+\iu\theta X_\alpha \in \Lieu_\CC$ and
$Z'=J X_\alpha+J\theta X_\alpha \in J\Lieu\subset \Lieg$.
Then $\sigma(Z)=\sigma(Z')$.
Denote by $\Liesl_\alpha(2,\CC)$ 
the three-dimensional simple subalgebra spanned by
$\{X_\alpha, \alpha^\vee, \theta X_\alpha\}$ over $\CC_J\coloneqq\RR\oplus\RR J$.
We identify $\Liesl_\alpha(2,\CC)$ with $\Liesl(2,\CC)$ by the following correspondence:
\[
X_\alpha \leftrightarrow \bmat0{\frac12}00,\qquad
\alpha^\vee \leftrightarrow \bmat100{-1},\qquad
\theta X_\alpha \leftrightarrow \bmat00{-\frac12}0.
\]
Let $U(\Liesl_\alpha(2,\CC))\varv=V^{(1)}+\ldots+V^{(t)}$ be an irreducible decomposition
as an $\Liesl_\alpha(2,\CC)$-module and
$\varv=\varv^{(1)}+\cdots+\varv^{(t)}$ the corresponding decomposition.
Then for each $s=1,\ldots,t$, we have
$\varv^{(s)}\ne0$ and $\sigma(\alpha^\vee)\varv^{(s)}=0$.
If we put $d_s = \dim_{\CC_J}V^{(s)}$,
then by the representation theory of $\Liesl(2,\CC)$,
we can make the following identification:
\[
\begin{aligned}
&V^{(s)}=\sum_{i=0}^{d_s-1}\CC x^{d_s-1-i}y^i \subset \CC[x,y],\\
&\sigma(X_\alpha+\theta X_\alpha)|_{V^{(s)}}=
 -\frac12\left(y\frac{\partial}{\partial x} - x\frac{\partial}{\partial y} \right), \qquad 
\sigma(\alpha^\vee)|_{V^{(s)}}=-x\frac{\partial}{\partial x}+y\frac{\partial}{\partial y}.
\end{aligned}
\]
We see $x^{d_s-1-i}y^i$ is a $\sigma(\alpha^\vee)$-eigenvector 
with eigenvalue $2i+1-d_s$ and because $\sigma(\alpha^\vee)\varv^{(s)}=0$,
$d_s$ is necessarily odd.
On the other hand, if we put $z=x+\iu y$, $\bar z=x-\iu y$, then
$V^{(s)}=\sum_{i=0}^{d_s-1}\CC z^{d_s-1-i}\bar z^i$ and
$z^{d_s-1-i}\bar z^i$ is a $\sigma(X_\alpha+\theta X_\alpha)$-eigenvector
with eigenvalue $\frac{\iu}2(d_s-1-2i)$.
Also, $\varv^{(s)}$ equals
\[
(z+\bar z)^{\frac{d_s-1}2}(z-\bar z)^{\frac{d_s-1}2}
=(z^2-\bar z^2)^{\frac{d_s-1}2}
\]
up to a scalar multiple.
Since $\sigma(Z')\,z^{d_s-1}=-\frac{d_s-1}2\,z^{d_s-1}$ and
$\sigma(Z')\,\bar z^{d_s-1}=\frac{d_s-1}2\,\bar z^{d_s-1}$,
it follows that
\[
\sigma(Z')(\sigma(Z')^2-1)\varv^{(s)}=0
\Leftrightarrow
d_s=1\text{ or }d_s=3
\Leftrightarrow
\sigma(X_\alpha)^2\varv^{(s)}=0.
\]
Thus we get
$\sigma(Z)(\sigma(Z)^2-1)\varv \Leftrightarrow \sigma(X_\alpha)^2\varv=0$,
which proves \eqref{eq:new_single}.

To prove \eqref{eq:new_double}
it suffices to show with respect to a $U$-invariant Hermitian inner product $(\cdot,\cdot)_V$ on $V$,
$V^\Lieh_\single$ equals the orthogonal complement of $V^\Lieh_\double$ in $V^\Lieh$.
But it can be checked in a quite similar way to the proof of Lemma~\ref{lem:direct_sum}.
\end{proof}
As a consequence of this proposition,
the condition that $V^\Lieh=V^\Lieh_\single$ is equivalent to 
the condition that twice a root is not a weight of $(\sigma,V)$.
Hence we get
\begin{cor}\label{cor:small}
A single-petaled $U$-type $(\sigma, V)$ is nothing but
an irreducible small representation of $G$ in the sense of \cite{Br}.
\end{cor}
\begin{defn}
We say an irreducible representation $(\sigma,V)$ of $G$
is {\it quasi-small\/} when $V^\Lieh_\single\ne0$, 
that is,
$(\sigma,V)$ is quasi-single-petaled as a $U$-type.
\end{defn}
Since $\theta$ is the conjugation map of $\Lieg$ with respect to 
the real form $\Lieu$,
the $U$-homomorphism
\[
\Lieg \oplus \Lieg
\xrightarrow{\id \oplus \, \theta}
\Lieg \oplus \Lieg
\xrightarrow{\frac{1-\iu J}2\cdot\, \oplus \frac{1+\iu J}2\cdot}
\frac{1-\iu J}2\Lieg \oplus \frac{1+\iu J}2\Lieg
= \Lieg_\CC
\]
gives an isomorphism
$\Lieg\oplus\Lieg\simeq \Lieg_\CC$
of complex Lie algebras.
If we identify these two complex Lie algebras,
then their subspaces correspond in the following way:
\begin{equation}\label{eq:ids}
\begin{aligned}
 \{(X,\theta X);\, X \in\Lieg \} &\leftrightarrow \Lieg, \quad
  &\Delta\Lieg \coloneqq \{(X,X);\, X \in\Lieg \} &\leftrightarrow \Lieu_\CC, \\
 \nabla\Lieg \coloneqq \{(X,-X);\, X \in\Lieg \} &\leftrightarrow (J\Lieu)_\CC, \quad
  &\{(H, -H);\, H \in \Lieh \} &\leftrightarrow (\Lieh_\RR)_\CC, \\
 \{(X,-X);\, X\in \bar\Lien+\Lien\} &\leftrightarrow \bigl(\Lieh_\RR\bigr)^\perp_\CC, \quad
  &\Lien \oplus \bar\Lien &\leftrightarrow \Lien_\CC.
\end{aligned}
\end{equation}
Here $\bigl(\Lieh_\RR\bigr)^\perp$ is the orthogonal complement of
$\Lieh_\RR$ in $J\Lieu$ and 
$\bar\Lien=\theta\Lien$.
Extend the $U$-isomorphism $\eta_0 : \Lieg\ni X\mapsto (-X,X)\in \nabla\Lieg \simeq (J\Lieu)_\CC$
to the algebra isomorphism $\eta_0 : S(\Lieg) \simarrow S(\nabla\Lieg) \simeq S((J\Lieu)_\CC)$.
Then the restriction of $\eta_0$ to $S(\Lieh)$ gives an isomorphism 
$S(\Lieh)\simarrow S((\Lieh_\RR)_\CC)$.
We denote its inverse by $\chi_0$.
Clearly $\chi_0$ commutes with the $W$-actions.
As a variation of the map
$\gamma_0 : S((J\Lieu)_\CC) \rightarrow S((\Lieh_\RR)_\CC)$,
define the map
\[
\Tilde\gamma_0 \coloneqq \chi_0 \circ\gamma_0\circ\eta_0 :
 S(\Lieg) \rightarrow S(\Lieh).
\]
Then \eqref{eq:Ch} induces the algebra isomorphism $\Tilde\gamma_0 : S(\Lieg)^G\simarrow S(\Lieh)^W$,
by which we identify the two algebras and denote both of them by the same symbol $\Tilde{\mathscr A}$.
Note that by \eqref{eq:ids}
$\gamma_0$ is reduced to the projection map
\[
S(\Lieg)=S(\Lieh) \oplus S(\Lieg) (\bar\Lien+\Lien)
\rightarrow S(\Lieh).
\]
Now the result of \cite{Br} is generalized to the case of a quasi-small representation as follows:
\begin{thm}\label{thm:exBroer}
For a finite-dimensional irreducible representation $(\sigma,V)$ of $G$, define the map
\[
 \Tilde\Gamma_0^\sigma : \Hom_G(V,S(\Lieg))\ni\Phi\mapsto \varphi \in \Hom_W(V^\Lieh,S(\Lieh))
\]
so that the image $\varphi$ is given by the composition
\[
	\varphi:V^\Lieh\hookrightarrow V\xrightarrow{\Phi}
	S(\Lieg)\xrightarrow{\Tilde\gamma_0}S(\Lieh).
\]
Then $\Tilde\Gamma_0^\sigma$ is an injective
$\Tilde{\mathscr A}$-homomorphism
\norbra{clearly $\Hom_G(V,S(\Lieg))$ and $\Hom_W(V,S(\Lieh))$ have natural $\Tilde{\mathscr A}$-module structures}.
On the other hand, $\Tilde\Gamma_0^\sigma$ is surjective if and only if
$(\sigma,V)$ is small.
Furthermore, for any $\varphi \in \Hom_W(V^\Lieh/V^\Lieh_\double,S(\Lieh))$
there exists $\Phi \in \Hom_G(V,S(\Lieg))$
such that $\Tilde\Gamma_0^\sigma(\Phi)=\varphi$.
\end{thm}
\begin{proof}
The theorem follows immediately from 
the results of \S\ref{sec:Ch}
and the fact that
$\Tilde\Gamma_0^\sigma$ coincides with the composition
\begin{align*}
\Hom_G(V,S(\Lieg))&=\Hom_U(V,S(\Lieg))\xrightarrow{\eta_0\circ\cdot}
\Hom_U(V,S((J\Lieu)_\CC)) \\
&\xrightarrow{\Gamma^\sigma_0}\Hom_W(V^\Lieh,S((\Lieh_\RR)_\CC))
\xrightarrow{\chi_0\circ\cdot}\Hom_W(V^\Lieh,S(\Lieh)).\qedhere
\end{align*}
\end{proof}
\begin{defn}
We define the map $\Tilde\gamma$ of $U(\Lieg)$ into $S(\Lieh)$
by the projection
\begin{equation}\label{eq:cpxHCproj}
U(\Lieg)=\left(\bar\Lien U(\Lieg)+U(\Lieg)\Lien\right)\oplus U(\Lieh)\rightarrow U(\Lieh)\simeq S(\Lieh)
\end{equation}
followed by the translation
\[
S(\Lieh) \ni f(\lambda) \mapsto f(\lambda-\Tilde\rho) \in S(\Lieh).
\]
Here we identified $S(\Lieh)$ with the space of holomorphic polynomials
on the (complex) dual space $\Lieh^*$ of $\Lieh$.
\end{defn}
\begin{lem}\label{lem:HCcomute}
For any $D_1\in U(\Lieg)$ and $D_2\in U(\Lieg)^\Lieh$,
$\Tilde\gamma(D_1D_2)=\Tilde\gamma(D_2D_1)=\Tilde\gamma(D_1)\Tilde\gamma(D_2)$.
\end{lem}
\begin{proof}
Let $\bar D_1$ and $\bar D_2$ be
the images of $D_1$ and $D_2$ under the projection \eqref{eq:cpxHCproj},
respectively.
Since $U(\Lieg)=\bar\Lien U(\bar\Lien+\Lieh)\oplus\left(U(\Lieh)+U(\Lieg)\Lien\right)
 =\left(\bar\Lien U(\Lieg) + U(\Lieh)\right)\oplus U(\Lieh+\Lien)\Lien$
as an $\ad(\Lieh)$-module,
we have $D_2 \in \left(U(\Lieh)+U(\Lieg)\Lien\right) \cap \left(\bar\Lien U(\Lieg) + U(\Lieh)\right)$.
Hence $D_1D_2\equiv D_1\bar D_2\pmod{U(\Lieg)\Lien}$ and
$D_2D_1\equiv \bar D_2D_1\pmod{\bar \Lien U(\Lieg)}$.
Since $[\Lien,\Lieh]=\Lien$, $[\Lieh,\bar\Lien]=\bar\Lien$, we get
$D_1\bar D_2 \equiv \bar D_1\bar D_2 \pmod{\bar \Lien U(\Lieg)+U(\Lieg)\Lien}$ and 
$\bar D_2 D_1 \equiv \bar D_2 \bar D_1\pmod{\bar \Lien U(\Lieg)+U(\Lieg)\Lien}$.
\end{proof}
The isomorphism $\Lieg \oplus \Lieg \simeq \Lieg_\CC$ induces the algebra isomorphism
$U(\Lieg)\otimes U(\Lieg) = U(\Lieg \oplus \Lieg) \simeq U(\Lieg_\CC)$,
which clearly commutes with the $U$-actions.
Define the map $\eta : U(\Lieg) \ni D \mapsto 1\otimes D \in U(\Lieg)\otimes U(\Lieg) \simeq U(\Lieg_\CC)$.
Then we obtain the direct sum decomposition 
\begin{equation}\label{eq:compdeco}
U(\Lieg_\CC) =
 U(\Lieg_\CC) \Lieu_\CC
  \oplus \eta\left( U(\Lieg) \right)
 \simeq \left( U(\Lieg \oplus \Lieg) \Delta\Lieg \right)
  \oplus \left( 1\otimes U(\Lieg) \right)
\end{equation}
as a $U$-module.
\begin{lem}\label{lem:HCtrans}
Under the composition map
\begin{equation}\label{eq:compcompo}
U(\Lieg) \xrightarrow{\eta} U(\Lieg_\CC) \xrightarrow{\gamma}
 S((\Lieh_\RR)_\CC) \xrightarrow{\chi_0} S(\Lieh),
\end{equation}
the image of $D\in U(\Lieg)$ equals $\Tilde\gamma(D)\left(\frac\lambda2\right)$.
\end{lem}
\begin{proof}
For $D\in U(\Lieg)$ and $X\in\Lien$
\[
\eta(DX)=1\otimes DX = (1\otimes D)\cdot(1\otimes X+X\otimes1)-(X\otimes1)\cdot(1\otimes D).
\]
Hence by \eqref{eq:ids},
$\eta\left(U(\Lieg)\Lien\right) 
 \subset \left( 1\otimes U(\Lieg) \right)\Delta\Lieg + \left( \Lien \otimes 1 \right)\left( 1\otimes U(\Lieg) \right)
 \subset U(\Lieg_\CC)\Lieu_\CC + \Lien_\CC U(\Lieg_\CC)$.
Also, since $\eta\left(\bar\Lien U(\Lieg)\right) =
 \left( 1 \otimes \bar\Lien \right)\left( 1\otimes U(\Lieg) \right)
 \subset \Lien_\CC U(\Lieg_\CC)$, we have
 $\gamma\circ\eta\left( \bar\Lien U(\Lieg)+U(\Lieg)\Lien\right)=0$.
On the other hand, if $f\in S(\Lieh) \simeq U(\Lieh)$ and $H \in \Lieh$,
then $\eta(f\cdot H)=(1\otimes f)\cdot (1\otimes \frac{H}2 + \frac{H}2\otimes 1)
 +(1\otimes \frac{H}2 - \frac{H}2\otimes 1)\cdot(1\otimes f)
 \equiv \eta_0(\frac{H}2)\,\eta(f) \pmod{ U(\Lieg_\CC)\Lieu_\CC}$.
Hence for any $f(\lambda)\in S(\Lieh)$,
$\eta(f(\lambda))\equiv\eta_0(f(\frac{\lambda}2)) \pmod{ U(\Lieg_\CC)\Lieu_\CC}$.
In addition, by the correspondence 
\begin{equation}\label{eq:eta0}
 \eta_0(H_1+JH_2)=(-H_1-JH_2, H_1+JH_2)\leftrightarrow
 -H_1-\iu H_2 \in (\Lieh_\RR)_\CC
 \qquad \forall H_1, H_2\in\Lieh_\RR,
\end{equation}
we have $\eta_0(f(\cdot-2\Tilde\rho))=(\eta_0(f))(\cdot+\rho)$
for $f\in S(\Lieh)$.
Therefore \eqref{eq:compcompo} maps
$f\in S(\Lieh)$ in the following way:
\[
f(\lambda) \xrightarrow{\eta} \eta(f(\lambda))
 \xrightarrow{\text{projection to }S((\Lieh_\RR)_\CC)}
 \eta_0\left(f\left(\frac{\lambda}2\right)\right)
 \xrightarrow{\rho\text{-shift}} 
  \eta_0\left(f\left(\frac{\lambda-2\Tilde\rho}2\right)\right)
 \xrightarrow{\chi_0} f\left(\frac{\lambda}2-\Tilde\rho\right).
\]
Since this image is $\Tilde\gamma(f)\left(\frac\lambda2\right)$,
we get the lemma.
\end{proof}
By the multiplicity function $\Tilde{\mathbf k} : \Sigma \ni \alpha \mapsto -1 \in \CC$,
we define the degenerate affine Hecke algebra $\mathbf H_{\Tilde{\mathbf k}}$
associated to the data $(\Lieh,\Pi,\Tilde{\mathbf k})$ as in Definition~\ref{defn:H}
and denote it simply by $\Tilde{\mathbf H}$.
The key relations in $\Tilde{\mathbf H}$ are
\begin{equation*}%
s_\alpha\cdot\xi
 =s_\alpha(\xi)\cdot s_\alpha
  + \alpha(\xi)\qquad\forall\alpha\in\Pi\ \forall\xi\in\Lieh.
\end{equation*}
As in the case of $S_{\mathbf H}((\Lieh_\RR)_\CC)\simeq S((\Lieh_\RR)_\CC)$,
the left $\Tilde{\mathbf H}$-module 
\[
S_{\Tilde{\mathbf H}}(\Lieh)\coloneqq \Tilde{\mathbf H}\ \Bigm/\!\!\sum_{w\in W\setminus\{1\}}\! \Tilde{\mathbf H}(w-1)
\]
is naturally identified with $S(\Lieh)$ as a left $S(\Lieh)$-module
and under this identification
the space of $W$-fixed elements in $S_{\Tilde{\mathbf H}}(\Lieh)$
equals $S(\Lieh)^W$.
\begin{lem}\label{lem:compHecke}
The map
\[
 S_{\mathbf H}((\Lieh_\RR)_\CC)\simeq S((\Lieh_\RR)_\CC)
  \xrightarrow{\chi_0} S(\Lieh) \ni f(\lambda)
  \mapsto f(2\lambda) \in S(\Lieh)
  \simeq S_{\Tilde{\mathbf H}}(\Lieh)
\]
commutes with the $W$-actions.
\end{lem}
\begin{proof}
Let $\alpha \in \Pi$ and put
$\Lieh_\RR(\alpha)=\{H\in\Lieh_\RR;\,\alpha(H)=0\}$. 
Then from Lemma~\ref{lem:Hecke},
the eigenspaces of $s_\alpha$
in $S_{\mathbf H}((\Lieh_\RR)_\CC)$ are
\[
S(\Lieh_\RR(\alpha)_\CC)\cdot\CC[(\alpha^\vee)^2],
\qquad S(\Lieh_\RR(\alpha)_\CC)\cdot\CC[(\alpha^\vee)^2](\alpha^\vee+2),
\]
which have eigenvalues $1,-1$, respectively.
If we apply the map in the lemma to them,
then by \eqref{eq:eta0} their images are respectively
\[
S(\Lieh_\RR(\alpha)\otimes \CC_J)\cdot\CC_J[(\alpha^\vee)^2],
\qquad S(\Lieh_\RR(\alpha)\otimes \CC_J)\cdot\CC_J[(\alpha^\vee)^2](\alpha^\vee-1),
\]
where $\CC_J=\RR\oplus\RR J$.
They, in turn, are shown to be the eigenspaces of $s_\alpha$
in $S_{\Tilde{\mathbf H}}(\Lieh)$ with eigenvalues $1,-1$ 
by the same argument as the proof of Lemma~\ref{lem:Hecke}.
\end{proof}
Now in view of \eqref{eq:compdeco}, Lemma~\ref{lem:HCtrans}, and Lemma~\ref{lem:compHecke},
the results of \S\ref{sec:HC} give the following generalization of the Harish-Chandra isomorphism:
\begin{thm}\label{thm:HCisomo}
For a finite-dimensional irreducible representation $(\sigma,V)$ of $G$, define the map
\[
 \Tilde\Gamma^\sigma : \Hom_G(V,U(\Lieg))\ni\Psi\mapsto \psi \in \Hom_\CC(V^\Lieh,S_{\Tilde{\mathbf H}}(\Lieh))
\]
so that the image $\psi$ is given by the composition
\[
	\varphi:V^\Lieh\hookrightarrow V\xrightarrow{\Psi}
	S(\Lieg)\xrightarrow{\Tilde\gamma}S(\Lieh)\simeq S_{\Tilde{\mathbf H}}(\Lieh).
\]
Then we have
\begin{enumerate}[{\normalfont (i)}]
\item\label{item:injHC}
$\Tilde\Gamma^\sigma$ is injective.

\item\label{item:doubleHC}
For any $\psi \in \Hom_W(V^\Lieh/V^\Lieh_\double,S_{\Tilde{\mathbf H}}(\Lieh))$
there exists $\Psi \in \Hom_G(V,U(\Lieg))$ such that $\Tilde\Gamma^\sigma(\Psi)=\psi$.

\item
$\Tilde\Gamma^\sigma\left( \Hom_G(V,U(\Lieg))\right) \subset 
 \Hom_W(V^\Lieh,S_{\Tilde{\mathbf H}}(\Lieh))$ if and only if
$(\sigma,V)$ is small.
If this condition is satisfied, then
from {\normalfont(\ref{item:injHC})}, {\normalfont(\ref{item:doubleHC})} we have the isomorphism
\begin{equation}\label{eq:gHCisomo}
\Tilde\Gamma^\sigma : \Hom_G(V,U(\Lieg)) \simarrow \Hom_W(V^\Lieh,S_{\Tilde{\mathbf H}}(\Lieh)).
\end{equation}

\item
In particular, for the trivial representation $(\triv,\CC)$, the map
$\Tilde\Gamma^\triv : \Hom_G(\triv,U(\Lieg)) \simarrow \Hom_W(\triv,S_{\Tilde{\mathbf H}}(\Lieh))$
is essentially equal to the classical Harish-Chandra isomorphism
\[
\Tilde\gamma : U(\Lieg)^G\simarrow S(\Lieh)^W.
\]
Since it is an algebra isomorphism by Lemma~\ref{lem:HCcomute},
in addition to $S(\Lieg)^G$ and $S(\Lieh)^W$,
we denote $U(\Lieg)^G$ also by $\Tilde{\mathscr A}$.

\item
For a general $(\sigma,V)$, $\Hom_G(V,U(\Lieg))$ and $\Hom_\CC(V^\Lieh,S_{\Tilde{\mathbf H}}(\Lieh))$
have the natural $\Tilde{\mathscr A}$-module structures which are intertwined by
$\Tilde\Gamma^\sigma$.
Especially, if $(\sigma,V)$ is small, then \eqref{eq:gHCisomo} is an
$\Tilde{\mathscr A}$-module isomorphism.
\end{enumerate}
\end{thm}
In the rest of this section,
we never refer to $\Lieg_\CC$ and
consider $\Lieg$ itself to be defined over $\CC$
by letting $J=\iu$.
Let $\Tilde B(\cdot,\cdot)$ be the Killing form for the complex Lie algebra $\Lieg$
and $\cang{\cdot}{\cdot}$ the bilinear form on $\Lieh^*\times\Lieh^*$
induced by $\Tilde B(\cdot,\cdot)$.
Note that $B(\cdot,\cdot)=2\Real \Tilde B(\cdot,\cdot)$.
Clearly, each irreducible constituent of the adjoint representation $(\Ad,\Lieg)$ is small and
$\Lieg^\Lieh=\Lieh$ is the reflection representation $(\refl,\Lieh)$ of $W$.
More generally, all the irreducible small representations of each type of complex simple Lie algebra
are listed in \cite{Ssl, Sml} (the classification is also given in \cite{{gHC}}).
As for quasi-small representations we have
\begin{prop}\label{prop:quasi-small}
Let $(\sigma_{\Tilde\rho},V_{\Tilde\rho})$ be the finite-dimensional irreducible
representation with highest weight $\Tilde\rho$
and $(\sigma, V)$ an arbitrary irreducible quasi-small representation.
Then $(\sigma, V)$ is isomorphic to
an irreducible constituent of the $G$-module $\End V_{\Tilde\rho}\simeq V_{\Tilde\rho}\otimes V_{\Tilde\rho}$.
Moreover, the multiplicity of $(\sigma,V)$ in $\End V_{\Tilde\rho}$ is
\[\dim_\CC\bigl\{\varv \in V^\Lieh;\,\sigma(X_\alpha)^2\varv=0\quad \forall\alpha\in\Pi, \forall X_\alpha\in\Lieg_\alpha\bigr\}.\]
\end{prop}
\begin{proof}
For any finite-dimensional irreducible representation $(\sigma,V)$ of $G$,
put $V^\Lieh(\Lien)\coloneqq\bigl\{\varv \in V^\Lieh;\,\sigma(X_\alpha)^2\varv=0\ \forall\alpha\in\Pi, \forall X_\alpha\in\Lieg_\alpha\bigr\}$.
Then the multiplicity of $(\sigma,V)$ in $\End V_{\Tilde\rho}$
equals $\dim_\CC V^\Lieh(\Lien)$ (\cite[4.3]{Adg}, \cite[Theorem 47]{Ko:Ad}).
In particular, if $(\sigma,V)$ is quasi-small,
then $V^\Lieh(\Lien)\supset V^\Lieh_\single\ne 0$.
\end{proof}
Let $\ell=\dim_\CC\Lieh$.
Then one has
$\bigwedge \Lieg \simeq \Bigl(\End V_{\Tilde\rho}\Bigr)^{\,\oplus\, 2^{\ell}}$ as $G$-modules
(see \cite{Ko:Ad,Adg}).
Hence each irreducible quasi-small representation
appears also in $\bigwedge \Lieg$.
Related to this, we have the following:
For $k=0,\ldots,\ell$, we consider the $G$-module $(\sigma,V)=(\bigwedge^k\Ad,\bigwedge^k\Lieg)$.
Although it may be reducible, we define $V^\Lieh$, $V^\Lieh_\single$, $V^\Lieh_\double$, and $\Tilde\Gamma^\sigma$
as in the irreducible case.
Observe that $\bigwedge^k\Lieh \subset V^\Lieh_\single$.
Let $\Tilde B_k(\cdot,\cdot)$ be the unique $G$-invariant non-degenerate symmetric bilinear form 
on $V\times V$ such that
\[
\Tilde B_k(u_1 \wedge\cdots\wedge u_k,\varv_1 \wedge\cdots\wedge \varv_k)=\det\left(\Tilde B(u_i,\varv_j)\right)_{1\le i,j\le k}.
\]
Since $\Tilde B_k(\cdot,\cdot)$ is non-degenerate both on 
$V^\Lieh\times V^\Lieh$ and on $\bigwedge^k\Lieh\times \bigwedge^k\Lieh$, 
we can define the orthogonal complement $\bigl(\bigwedge^k\Lieh\bigr)^\perp$
of $\bigwedge^k\Lieh$ in $V^\Lieh$ with respect to $\Tilde B_k(\cdot,\cdot)$.
It is easy to see
$V^\Lieh_\double \subset \bigl(\bigwedge^k\Lieh\bigr)^\perp$.
Therefore it follows from Theorem~\ref{thm:HCisomo} (\ref{item:injHC}), (\ref{item:doubleHC})
that for any $\psi \in \Hom_W\left(V^\Lieh/\bigl(\bigwedge^k\Lieh\bigr)^\perp,S_{\Tilde{\mathbf H}}(\Lieh)\right)$
there exists a unique $\Psi \in \Hom_G(V,U(\Lieg))$ such that $\Tilde\Gamma^\sigma(\Psi)=\psi$.
We denote the set of all such $\Psi$ by $\mathscr M_k$, in other words, we put
\[
\mathscr M_k=\Bigl\{\,\Psi\in\Hom_G(V,U(\Lieg));\,
 \Tilde\Gamma^\sigma(\Psi)[\varv]=0\quad\forall\varv\in\Bigl(\bigwedge^k\Lieh\Bigr)^\perp\,\Bigl\}.
\]
Then $\mathscr M_k$ is an $\Tilde{\mathscr A}$-submodule of the
$\Tilde{\mathscr A}$-module $\Hom_G(V,U(\Lieg))$.
\begin{thm}\label{thm:adprod}
Define the $\Tilde{\mathscr A}$-homomorphism $\omega : \bigwedge^k\Hom_G(\Lieg,U(\Lieg))
 \rightarrow \Hom_G\left(\bigotimes^k\Lieg,U(\Lieg)\right)$ so that
the image of\/ $\Psi_1\wedge\cdots\wedge\Psi_k \in \bigwedge^k\Hom_G(\Lieg,U(\Lieg))$ is given by
\begin{equation}\label{eq:adprodmap}
\omega\bigl( \Psi_1\wedge\cdots\wedge\Psi_k \bigr) :
\bigotimes^k\Lieg \ni X_1\otimes\cdots\otimes X_k\longmapsto
 \coldet\left(\Psi_i[X_j]\right)_{1\le i,j\le k} \in U(\Lieg).
\end{equation}
Here the symbol $\coldet$ in \eqref{eq:adprodmap} stands for a so-called `column-determinant',
that is,
\[
 \coldet\left(\Psi_i[X_j]\right)_{1\le i,j\le k}
 =\sum_{\mu \in\mathfrak S_k}(\sgn\mu)\,
  \Psi_{\mu(1)}[X_1]\cdots\Psi_{\mu(k)}[X_k].
\]
\norbra{$\mathfrak S_k$ denotes the $k$-th symmetric group.}
Then, for any $\Psi_1\wedge\cdots\wedge\Psi_k \in \bigwedge^k\Hom_G(\Lieg,U(\Lieg))$, 
$X_1,\ldots,X_k\in\Lieg$, and
$\tau\in\mathfrak S_k$,
\begin{equation}\label{eq:adcolcmute}
\omega\left(\Psi_1\wedge\cdots\wedge\Psi_k\right)
 [X_{\tau(1)}\otimes\cdots\otimes X_{\tau(k)}]
=(\sgn \tau)\,\omega\left(\Psi_1\wedge\cdots\wedge\Psi_k\right)
 [X_1\otimes\cdots\otimes X_k].
\end{equation}
By this, we consider $\omega\left(\Psi_1\wedge\cdots\wedge\Psi_k\right)\in\Hom_G(V, U(\Lieg))$.
Then 
$\omega\left(\Psi_1\wedge\cdots\wedge\Psi_k\right)
\in \mathscr M_k$ and moreover,
\begin{equation}\label{eq:adprodIm}
\Tilde\Gamma^\sigma\bigl(\omega\left(\Psi_1\wedge\cdots\wedge\Psi_k\right))
 [H_1\wedge\cdots\wedge H_k]
 = \det\left(\Tilde\Gamma^{\Ad}(\Psi_i)[H_j]\right)_{1\le i,j\le k}
 \qquad \forall H_1,\ldots,H_k \in \Lieh.
\end{equation}
Here the right-hand side is the determinant of an $S(\Lieh)$-valued matrix.
Furthermore, $\omega : \bigwedge^k\Hom_G(\Lieg,U(\Lieg))
 \rightarrow \mathscr M_k$
is an $\Tilde{\mathscr A}$-module isomorphism.
\end{thm}
\begin{proof}
By an elementary argument,
we can see \eqref{eq:adcolcmute} follows
if we prove it for the special case where $k=2$ and $\tau=(1,2)$.
Hence for a while we assume $k=2$ and $\tau=(1,2)$.
In order to show \eqref{eq:adcolcmute},
for any $\Psi_1, \Psi_2\in \Hom_G(\Lieg,U(\Lieg))$ define
$\Psi \in \Hom_G(\Lieg\otimes\Lieg,U(\Lieg))$ by
\[\begin{aligned}
\Psi[X_1\otimes X_2]
 &=\omega(\Psi_1\wedge \Psi_2)[X_1\otimes X_2+X_2\otimes X_1]\\
 &=\coldet\bmat{\Psi_1[X_1]}{\Psi_1[X_2]}{\Psi_2[X_1]}{\Psi_2[X_2]}
 +\coldet\bmat{\Psi_1[X_2]}{\Psi_1[X_1]}{\Psi_2[X_2]}{\Psi_2[X_1]}\\
 &= \Psi_1[X_1]\Psi_2[X_2]-\Psi_2[X_1]\Psi_1[X_2]
  + \Psi_1[X_2]\Psi_2[X_1]-\Psi_2[X_2]\Psi_1[X_1],
\end{aligned}\]
and let us prove $\Psi=0$.
By Theorem~\ref{thm:HCisomo} (\ref{item:injHC}),
it suffices to show $\Tilde\Gamma^{\Ad^{\otimes2}}(\Psi)=0$.
Assume $\Tilde\Gamma^{\Ad^{\otimes2}}(\Psi)\ne0$.
First we note
$(\Lieg\otimes\Lieg)^\Lieh=\Lieh\otimes\Lieh + \sum_{\alpha\in\Sigma}\Lieg_{-\alpha}\otimes\Lieg_{\alpha}$.
For $H_1\otimes H_2\in\Lieh\otimes\Lieh$, Lemma~\ref{lem:HCcomute} implies
$\Tilde\gamma\left( \Psi_1[H_1]\Psi_2[H_2] \right) = \Tilde\gamma\left( \Psi_2[H_2]\Psi_1[H_1] \right)$, 
$\Tilde\gamma\left( \Psi_1[H_2]\Psi_2[H_1] \right) = \Tilde\gamma\left( \Psi_2[H_1]\Psi_1[H_2] \right)$,
and hence $\Tilde\Gamma^{\Ad^{\otimes2}}(\Psi)[H_1\otimes H_2]=0$.
For $d\in\ZZ_{\ge0}$, put $S^d_{\Tilde{\mathbf H}}(\Lieh)=\bigoplus_{i=0}^d S^i(\Lieh)$
and define the projection map $q_d : S^d_{\Tilde{\mathbf H}}(\Lieh)\rightarrow S^d(\Lieh)$.
Let $d'\in\ZZ_{\ge0}$ be such that
$\Tilde\Gamma^{\Ad^{\otimes2}}(\Psi) \in \Hom_\CC(V^\Lieh, S^{d'}_{\Tilde{\mathbf H}}(\Lieh))
\setminus \Hom_\CC(V^\Lieh, S^{d'-1}_{\Tilde{\mathbf H}}(\Lieh))$.
Then we can easily observe
$q_{d'}\circ \Gamma^{\Ad^{\otimes2}}(\Psi) \in \Hom_W(V^\Lieh, S^{d'}(\Lieh))\setminus\{0\}$.
Hence there exist $\alpha\in\Pi$, $X_{-\alpha}\in \Lieg_{-\alpha}$, 
and $X_{\alpha}\in\Lieg_\alpha$ such that
$\Tilde\Gamma^{\Ad^{\otimes2}}(\Psi)[X_{-\alpha}\otimes X_\alpha]\ne 0$.
Put
\[
\Liesl_\alpha(2,\CC) =\Lieg_{-\alpha}+\CC\alpha^\vee+\Lieg_\alpha,\quad
\Lieh(\alpha) =\{H\in\Lieh;\,\alpha(H)=0\},\quad
\Lien_\alpha =\sum_{\beta\in\Sigma^+\setminus\{\alpha\}}\Lieg_\beta,\quad
\bar\Lien_\alpha =\theta\Lien_\alpha
\]
and define the projection map
\[
\Tilde\gamma_\alpha :
U(\Lieg)=\left(\bar\Lien_\alpha U(\Lieg) + U(\Lieg) \Lien_\alpha\right)
 \oplus U(\Lieh(\alpha)+\Liesl_\alpha(2,\CC)) 
\rightarrow
U(\Lieh(\alpha)+\Liesl_\alpha(2,\CC)). 
\]
Then $\gamma_\alpha$ is an $\Liesl_\alpha(2,\CC)$-homomorphism and
$\Tilde\gamma\circ\Tilde\gamma_\alpha=\Tilde\gamma$.
Also, in a similar way to the proof of Lemma~\ref{lem:HCcomute},
we can show for any $D_1\in U(\Lieg)$ and $D_2\in U(\Lieg)^{\Lieh(\alpha)}$,
\[
\Tilde\gamma_\alpha(D_1D_2)=\Tilde\gamma_\alpha(D_1)\Tilde\gamma_\alpha(D_2),
\qquad 
\Tilde\gamma_\alpha(D_2D_1)=\Tilde\gamma_\alpha(D_2)\Tilde\gamma_\alpha(D_1).
\]
Define the $\Liesl_\alpha(2,\CC)$-homomorphism
$\Psi^\alpha : \Liesl_\alpha(2,\CC)\otimes\Liesl_\alpha(2,\CC) \rightarrow 
U(\Lieh(\alpha)+\Liesl_\alpha(2,\CC))$ by
\[
\Liesl_\alpha(2,\CC)\otimes\Liesl_\alpha(2,\CC)
\hookrightarrow \Lieg\otimes\Lieg \xrightarrow{\Psi} U(\Lieg)
 \xrightarrow{\Tilde\gamma_\alpha} U(\Lieh(\alpha)+\Liesl_\alpha(2,\CC)),
\]
and the $\Liesl_\alpha(2,\CC)$-homomorphism
$\Psi_i^\alpha : \Liesl_\alpha(2,\CC) \rightarrow 
U(\Lieh(\alpha)+\Liesl_\alpha(2,\CC))$ ($i=1,2$) by
\[
\Liesl_\alpha(2,\CC)
\hookrightarrow \Lieg \xrightarrow{\Psi_i} U(\Lieg)
 \xrightarrow{\Tilde\gamma_\alpha} U(\Lieh(\alpha)+\Liesl_\alpha(2,\CC)).
\]
Since $\Psi_i[X_j] \in U(\Lieg)^{\Lieh(\alpha)}$
for any $X_1, X_2\in \Liesl_\alpha(2,\CC)$ ($i,j=1,2$),
we have
\begin{equation}\label{eq:Psi_alpha}
\Psi^\alpha[X_1\otimes X_2]
 = \Psi_1^\alpha[X_1]\Psi_2^\alpha[X_2]-\Psi_2^\alpha[X_1]\Psi_1^\alpha[X_2]
  + \Psi_1^\alpha[X_2]\Psi_2^\alpha[X_1]-\Psi_2^\alpha[X_2]\Psi_1^\alpha[X_1].
\end{equation}
Now from Kostant's theorem (\cite{Ko1}),
the inclusion map
$\iota^\alpha : \Liesl_\alpha(2,\CC) \hookrightarrow U(\Lieh(\alpha)+\Liesl_\alpha(2,\CC))$
satisfies
\[
\Hom_{\Liesl_\alpha(2,\CC)}(\Liesl_\alpha(2,\CC),U(\Lieh(\alpha)+\Liesl_\alpha(2,\CC)))
 = U(\Liesl_\alpha(2,\CC))^{\Liesl_\alpha(2,\CC)}\cdot U(\Lieh(\alpha))\cdot \iota^\alpha.
\]
Hence there exist 
$Z_1,Z_2\in U(\Liesl_\alpha(2,\CC))^{\Liesl_\alpha(2,\CC)}\cdot U(\Lieh(\alpha))$
such that $\Psi_i^\alpha=Z_i\cdot\iota^\alpha$ ($i,j=1,2$).
Accordingly, the right-hand side of \eqref{eq:Psi_alpha} becomes
\[\begin{aligned}
Z_1\iota^\alpha[X_1]\,&Z_2\iota^\alpha[X_2]-Z_2\iota^\alpha[X_1]\,Z_1\iota^\alpha[X_2]
  + Z_1\iota^\alpha[X_2]\,Z_2\iota^\alpha[X_1]-Z_2\iota^\alpha[X_2]\,Z_1\iota^\alpha[X_1] \\
&=Z_1Z_2\left(\iota^\alpha[X_1]\,\iota^\alpha[X_2]-\iota^\alpha[X_1]\,\iota^\alpha[X_2]
  + \iota^\alpha[X_2]\,\iota^\alpha[X_1]-\iota^\alpha[X_2]\,\iota^\alpha[X_1]\right) =0.
\end{aligned}\]
In particular we have
$\Tilde\Gamma^{\Ad^{\otimes2}}(\Psi)[X_{-\alpha}\otimes X_\alpha]
=\Tilde\gamma\circ\Psi^\alpha[X_{-\alpha}\otimes X_\alpha]=0$,
a contradiction.
Thus we get $\Tilde\Gamma^{\Ad^{\otimes2}}(\Psi)=0$ and hence \eqref{eq:adcolcmute}.

Suppose $k$ is general and
$\Psi_1\wedge\cdots\wedge\Psi_k \in \bigwedge^k\Hom_G(\Lieg,U(\Lieg))$.
Put $\Lieg_0\coloneqq\Lieh$.
Take $\alpha_1,\ldots,\alpha_k\in\Sigma\cup\{0\}$
so that $\alpha_1+\cdots+\alpha_k=0$
and at least one $\alpha_j$ is not $0$.
Also, take $X_{\alpha_j} \in \Lieg_{\alpha_j}$
for each $\alpha_j$ and consider
the element $X_{\alpha_1}\wedge\cdots \wedge X_{\alpha_k}\in V^\Lieh$.
Then there exists at least one $j_0=1,\ldots,k$
such that $\alpha_{j_0}\in\Sigma^+$. 
Since $\Psi_i[X_{\alpha_{j_0}}]\in U(\Lieg)\Lien$ ($i=1,\ldots,k$), 
\eqref{eq:adcolcmute} implies
$\omega\left(\Psi_1\wedge\cdots\wedge\Psi_k\right)[X_{\alpha_1}\wedge\cdots \wedge X_{\alpha_k}] \in U(\Lieg)\Lien$
and hence $\Gamma^{\sigma}\bigl(\omega\left(\Psi_1\wedge\cdots\wedge\Psi_k\right)\bigr)[X_{\alpha_1}\wedge\cdots \wedge X_{\alpha_k}]=0$.
Since such $X_{\alpha_1}\wedge\cdots \wedge X_{\alpha_k}$ span
$\bigl(\bigwedge^k\Lieh\bigr)^\perp$,
we get $\omega\left(\Psi_1\wedge\cdots\wedge\Psi_k\right) \in \mathscr M_k$.
Obviously \eqref{eq:adprodIm} follows from Lemma~\ref{lem:HCcomute}.

Lastly the next lemma assures
$\omega : \bigwedge^k\Hom_G(\Lieg,U(\Lieg))
 \rightarrow \mathscr M_k$ is an isomorphism.
\end{proof}
\begin{lem}\label{lem:}
Define the $\Tilde{\mathscr A}$-homomorphism
$\omega_0 : \bigwedge^k \Hom_W(\Lieh,S(\Lieh))
 \rightarrow \Hom_W\left( \bigwedge^k\Lieh,S(\Lieh)\right)$ 
so that the image of $\varphi_1\wedge\cdots\wedge\varphi_k$ is given by the map
\[
  \bigwedge^k\Lieh\ni H_1\wedge\cdots\wedge H_k
  \mapsto \det\left(\varphi_i[H_j] \right)_{1\le i,j\le k} \in S(\Lieh).
\]
Then it is an isomorphism.
\end{lem}
\begin{proof}
$\Tilde B_k(\cdot,\cdot)$
induces the $W$-module isomorphism
$\left(\bigwedge^k\Lieh\right)^*\simeq \bigwedge^k\Lieh$.
Hence we have the following natural $\Tilde{\mathscr{A}}$-module isomorphisms:
$\Hom_W\left( \bigwedge^k\Lieh,S(\Lieh)\right)\simeq 
\left(\left(\bigwedge^k\Lieh\right)^* \otimes S(\Lieh) \right)^W
\simeq \left(\bigwedge^k\Lieh \otimes S(\Lieh) \right)^W
\simeq \{\text{$W$-invariant polynomial coefficient $p$-form on $\Lieh_\RR^*$}\}$.
Suppose $I_1,\ldots,I_\ell$ are algebraically independent homogeneous elements of $\Tilde{\mathscr A}=S(\Lieh)^W$
and they constitute a generator system of $\Tilde{\mathscr A}$. 
Then $\{ dI_{i_1}\wedge\cdots\wedge dI_{i_k};\,1\le i_1<\cdots<i_k\le\ell\}$
forms a basis of $\left(\bigwedge^k\Lieh \otimes S(\Lieh) \right)^W$ over $\Tilde{\mathscr A}$
(\cite{So}).
In particular, $\{dI_i;\,i=1,\ldots,\ell\}$ is a basis of
$(\Lieh \otimes S(\Lieh))^W\simeq \Hom_W(\Lieh,S(\Lieh))$
over $\Tilde{\mathscr A}$.
It is easy to check under the identification $\Hom_W\left( \bigwedge^k\Lieh,S(\Lieh)\right)\simeq \left(\bigwedge^k\Lieh \otimes S(\Lieh) \right)^W$,
the image $\omega_0(dI_{i_1}\wedge\cdots\wedge dI_{i_k}) \in \Hom_W\left( \bigwedge^k\Lieh,S(\Lieh)\right)$
of $dI_{i_1}\wedge\cdots\wedge dI_{i_k} \in \bigwedge^k \Hom_W(\Lieh,S(\Lieh))$
equals $dI_{i_1}\wedge\cdots\wedge dI_{i_k} \in \left(\bigwedge^k\Lieh \otimes S(\Lieh) \right)^W$.
\end{proof}
To find all the equivalence classes of irreducible quasi-small representations
and to determine the $W$-module structure of $V^\Lieh_\single$ for each irreducible quasi-small representation
$(\sigma,V)$ are both important problems.
The next lemma seems to be a help to solving them:
\begin{lem}\label{lem:Casimir}
Suppose $\Omega_\Lieg \in U(\Lieg)$ is the Casimir element of $\Lieg$.
That is, if we choose $X_\alpha\in\Lieg_\alpha$ for each $\alpha\in\Sigma$
so that $\Tilde B(X_\alpha,X_{-\alpha})=1$ and 
if we take a basis $\{H_1,\ldots,H_\ell\}$ of $\Lieh$
so that $\Tilde B(H_i,H_j)=\delta_{ij}$,
then
\[
\Omega_\Lieg=\sum_{i-1}^\ell H_i^2 
 + \sum_{\alpha\in\Sigma^+}(X_\alpha X_{-\alpha}+X_{-\alpha}X_\alpha).
\]
Define the following central element of \/$\CC[W]$\,:
\[
\Omega_W=\sum_{\alpha\in\Sigma^+}\cang{\alpha}{\alpha}(1-s_\alpha).
\]
Suppose $(\sigma,V)$ is an irreducible quasi-small representation of $G$
and its highest weight is $\lambda$.
Then for any $\varv\in V^\Lieh_\single$,
\begin{equation}\label{eq:Casimir}
\Omega_\Lieg \varv = \cang{\lambda}{\lambda+2\Tilde\rho}\varv = \Omega_W \varv.
\end{equation}
\end{lem}
\begin{proof}
Since $\Omega_\Lieg$ acts on $V$ by a scalar,
we get the first equality of \eqref{eq:Casimir}
by calculating the action of $\Omega_\Lieg$ on a highest weight vector of $V$.
To show the second equality,
take an arbitrary $\alpha\in\Sigma^+$
and define $\Liesl_\alpha(2,\CC)$ to be
the three-dimensional simple subalgebra spanned by
$\{X_\alpha, \alpha^\vee, X_{-\alpha}\}$.
If we consider $U(\Liesl_\alpha(2,\CC))\varv$ as an $\Liesl_\alpha(2,\CC)$-module,
then from \eqref{eq:new_single}
each irreducible constituent of $U(\Liesl_\alpha(2,\CC))\varv$
is isomorphic either to the trivial representation or to the adjoint representation.
Hence if we put $\varv_0\coloneqq\frac{1+s_\alpha}2\varv$,
$\varv_1\coloneqq\frac{1-s_\alpha}2\varv$,
then $\Liesl_\alpha(2,\CC)$ acts trivially on $\varv_0$. 
Thus $(X_\alpha X_{-\alpha}+X_{-\alpha}X_\alpha)\varv_0=0$.
On the other hand, if $\varv_1\ne0$, then
there exists an isomorphism from  
$\Liesl_\alpha(2,\CC)\varv_1$ to $\Liesl_\alpha(2,\CC)$ such that
$\varv_1\mapsto\alpha^\vee$.
Since
\[
\ad(X_\alpha X_{-\alpha}+X_{-\alpha}X_\alpha)\alpha^\vee
=2\ad(X_\alpha)X_{-\alpha}-2\ad(X_{-\alpha})X_{\alpha}
=4[X_\alpha, X_{-\alpha}]=2\cang{\alpha}{\alpha}\alpha^\vee,
\]
we get $(X_\alpha X_{-\alpha}+X_{-\alpha}X_\alpha)\varv_1=2\cang{\alpha}{\alpha} \varv_1$.
Therefore the second equality of \eqref{eq:Casimir} holds.
\end{proof}
\begin{exm}\label{exm:B2}
Suppose $\Lieg$ is the complex simple Lie algebra of type $(B_2)$.
As usual, take a basis $\{e_1,e_2\}$ of $\Lieh^*$ so that
$\Sigma^+=\{e_i\pm e_2, e_1, e_2\}$ and
$\cang{e_i}{e_j}=\frac16\delta_{ij}$.

In addition to the equivalence class of the trivial representation `$\triv$', 
that of the reflection representation `$\refl$', and that of the sign representation `$\sgn$',
we have two other equivalence classes of irreducible representations of $W$:
One is the class of the one-dimensional representation $\tau$
which takes the value $\tau(w)=1$ or $-1$ ($w\in W$)
according as 
the number of appearances of the reflections associated to short roots $\{\pm e_1,\pm e_2\}$
is even or odd
when we express $w$ as a product of reflections.
The other is the class of $\tau\otimes\sgn$,
which behaves similarly for the long roots.
On each irreducible representation of $W$, $\Omega_W$ acts by the scalar
whose value is indicated in Table~\ref{tbl:W}.
\begin{table}[!h]
\begin{center}
\begin{minipage}{.47\linewidth}
\begin{center}
\caption{values of $\Omega_W$}\label{tbl:W}
\begin{tabular}{|c||c|c|c|c|c|}
\hline%
rep. & $\triv$ & $\refl$ & $\sgn$ & $\tau$ & $\tau\otimes\sgn$ \\
\hline
$\Omega_W$ & $0$ & $1$& $2$& $\Bigl.\frac23\Bigr.$& $\frac43$ \\
\hline
\end{tabular}
\end{center}
\end{minipage}
\begin{minipage}{.47\linewidth}
\begin{center}
\caption{values of $\Omega_\Lieg$}\label{tbl:g}
\begin{tabular}{|c||c|c|c|c|}
\hline
rep. & $\sigma_{(0,0)}$ & $\sigma_{(1,0)}$ & $\sigma_{(1,1)}$ & $\sigma_{(2,1)}$ \\
\hline
$\Omega_\Lieg$ & $0$ & $\Bigl.\frac23\Bigr.$ & $1$& $2$  \\
\hline
\end{tabular}
\end{center}
\end{minipage}
\end{center}
\end{table}
The set of the highest weights of all finite-dimensional irreducible representations of $G$
is $\{(i+j)e_1+je_2;\, i,j \in \ZZ_{\ge 0}\}$.
Let $\sigma_{(i,j)}$ denote the finite-dimensional irreducible representation of $G$
with highest weight $ie_1+je_2$.
Note that $2\Tilde\rho=3e_1+e_2$ and that
$\cang{\lambda}{\lambda+2\Tilde\rho}
=\cang{\lambda+\Tilde\rho}{\lambda+\Tilde\rho}-\cang{\Tilde\rho}{\Tilde\rho}$.
Then we easily recognize
a representation $\sigma_{(i,j)}$
for which the value of $\sigma_{(i,j)}(\Omega_\Lieg)$ coincides with
one of the values of $\Omega_W$ in Table~\ref{tbl:W}
is isomorphic to one of those representations in Table~\ref{tbl:g}.

Among the irreducible representations in Table~\ref{tbl:g}, 
$\sigma_{(0,0)}, \sigma_{(1,0)}$, and $\sigma_{(1,1)}$ are small (cf.~Figure~\ref{fig:B2qs}).
Hence by Lemma~\ref{lem:Casimir},
their $0$-weight spaces have neither
$\sgn$ nor $\tau\otimes\sgn$ as a constituent.
On the other hand, $\sigma_{(2,1)}$ is not small because it has the weight $2e_1$.
But since the $W$-module $\bigwedge^2\Lieh$ is isomorphic to $\sgn$,
the argument before Theorem~\ref{thm:adprod} implies
that there is an irreducible quasi-small submodule $V$ of $\bigwedge^2\Lieg$
such that $\sgn$ appears in $V^\Lieh_\single$ as a constituent.
From Lemma~\ref{lem:Casimir}, Table~\ref{tbl:W}, and Table~\ref{tbl:g}, 
we conclude this $V$ must be isomorphic to $\sigma_{(2,1)}$.
\begin{figure}[!h]
\begin{center}
\caption{quasi-small representations for type $(B_2)$}
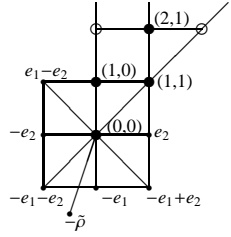
\label{fig:B2qs}
\begin{picture}(85,95)
 \put(15,20){\line(1,1){70}}
 \put(15,60){\line(1,-1){40}}
 \put(15,20){\line(1,0){40}}
 \put(15,40){\line(1,0){40}}
 \put(15,60){\line(1,0){40}}
 \put(35,80){\line(1,0){40}}
 \put(15,20){\line(0,1){40}}
 \put(35,20){\line(0,1){70}}
 \put(55,20){\line(0,1){70}}
 \put(35,40){\circle*{4}} 
 \put(35,60){\circle*{4}}
 \put(35,80){\circle{4}}
 \put(55,60){\circle*{4}}
 \put(55,80){\circle*{4}}
 \put(75,80){\circle{4}}
 \put(15,20){\circle*{2}}
 \put(15,40){\circle*{2}}
 \put(15,60){\circle*{2}}
 \put(35,20){\circle*{2}}
 \put(55,20){\circle*{2}}
 \put(55,40){\circle*{2}}
 \put(25,10){\circle*{2}}
 \put(39,41){$\scriptstyle (0,0)$}
 \put(37,62){$\scriptstyle (1,0)$}
 \put(8,62){$\scriptstyle e_1-e_2$}
 \put(2,39){$\scriptstyle -e_2$}
 \put(2,14){$\scriptstyle -e_1-e_2$}
 \put(37,14){$\scriptstyle -e_1$}
 \put(54,14){$\scriptstyle -e_1+e_2$}
 \put(57,39){$\scriptstyle e_2$}
 \put(58,58){$\scriptstyle (1,1)$}
 \put(57,82){$\scriptstyle (2,1)$}
 \put(23,5){$\scriptstyle -\Tilde\rho$}
 \thinlines
 \put(25,10){\line(1,3){10}}
\end{picture}
\end{center}
\end{figure}
$\sigma_{(2,1)}$ is an example of an irreducible representation of $G$
which is quasi-small but not small.
$\tau\otimes\sgn$ is an example of an irreducible representation of $W$
which does not appear in $V^\Lieh_\single$ as a constituent for any
irreducible quasi-small representation $(\sigma, V)$ of $G$.
\end{exm}
Suppose $\lambda\in\Lieh^*$.
In the rest of this section, we apply Theorem~\ref{thm:HCisomo}
to the study of a new relation we shall establish between
a basic $\Tilde{\mathbf H}$-module 
\begin{equation}\label{eq:Al}
A(\lambda)\coloneqq S_{\Tilde{\mathbf H}}(\Lieh) / \sum_{f\in S(\Lieh)^W} (f-f(\lambda))S_{\Tilde{\mathbf H}}(\Lieh)
\end{equation}
and the Verma module 
\[
M(\lambda)\coloneqq U(\Lieg)/\sum_{H\in\Lieh}U(\Lieg)(H-(\lambda-\Tilde\rho)(H))+U(\Lieg)\Lien
\]
of $U(\Lieg)$.
Hereafter we assume $\lambda\in\Lieh^*$ is dominant,
that is,
$\lambda(\alpha^\vee)\notin\{-1,-2,\ldots\}$
for any $\alpha\in\Sigma^+$.
In this case, it is known that
$A(\lambda)\simeq \Tilde{\mathbf H} / \sum_{\xi\in\Lieh}\Tilde{\mathbf H}(\xi-\lambda(\xi))$
(cf.~\cite{ChJ}).
\begin{defn}\label{defi:V-A}
We introduce the correspondence
$\Xi_\lambda : \bigl\{\, \text{$\Tilde{\mathbf H}$-submodule of }A(\lambda) \,\bigr\} \rightarrow$\newline
$\bigl\{\, \text{$U(\Lieg)$-submodule of }M(\lambda) \,\bigr\}
$ defined by
\[
S\mapsto \sum \Bigl\{\,V M(\lambda);\, 
\text{$V\subset U(\Lieg)$ is $\ad(\Lieg)$-stable
and satisfies $\Tilde\gamma(V) A(\lambda) \subset S$}\,\Bigr\},
\]
and the correspondence
$Y_\lambda : \bigl\{\, \text{$U(\Lieg)$-submodule of }M(\lambda) \,\bigr\} \rightarrow
\bigl\{\, \text{$\Tilde{\mathbf H}$-submodule of }A(\lambda) \,\bigr\}
$
defined by
\[
T\mapsto \sum \Bigl\{\,\Tilde{\mathbf H}\,\Tilde\gamma(V) A(\lambda);\, 
\text{$V\subset U(\Lieg)$ is $\ad(\Lieg)$-stable and satisfies $V M(\lambda) \subset T$}\,\Bigr\}.
\]
Obviously they are well-defined and preserve any inclusion relation.
\end{defn}
\begin{prop}\label{prop:V-A}
$\Xi_\lambda(A(\lambda))=M(\lambda)$ and
$\Xi_\lambda(0)=0$.
Let $S$ be an arbitrary
$\Tilde{\mathbf H}$-submodule of $A(\lambda)$.
Then $Y_\lambda\circ\Xi_\lambda(S) \subset S$.
Moreover, suppose $S$ admits a $W$-stable subspace $E$ with the following properties:
{\normalfont(a)}
$S=\Tilde{\mathbf H}E$;
{\normalfont(b)}
as a $W$-module, 
each irreducible constituent of $E$ belongs to $\widehat{W}_\single$,
where $\widehat{W}_\single$ denotes the set of equivalence classes
of those irreducible representations of $W$ which appear in $V^\Lieh_\single$
for some irreducible quasi-small representation $(\sigma,V)$ of $G$.
Then $Y_\lambda\circ\Xi_\lambda(S)=S$.
\end{prop}
\begin{proof}
It is clear that $\Xi_\lambda(A(\lambda))=M(\lambda)$.
To show $\Xi_\lambda(0)=0$,
let $T_0$ be the unique irreducible $U(\Lieg)$-submodule of $M(\lambda)$.
Then $T_0$ is isomorphic to $M(w_0\lambda)$ for some $w_0\in W$.
Fix a highest weight vector $\varv_0$ of $T_0$.

Suppose an $\ad(\Lieg)$-stable subspace $V$
satisfies $\Tilde\gamma(V) A(\lambda) = 0$.
Then by the definition \eqref{eq:Al} of $A(\lambda)$,
$\Tilde\gamma(V) \subset \sum_{f\in S(\Lieh)^W}S(\Lieh)\left(f-f(\lambda)\right)
=\sum_{f\in S(\Lieh)^W}S(\Lieh)\left(f-f(w_0\lambda)\right)$.
By using the direct sum decomposition $U(\Lieg)=\bar\Lien U(\bar\Lien+\Lieh) \oplus U(\Lieh) \oplus U(\Lieg)\Lien$,
we have $V \varv_0\subset \bigl(\Tilde\gamma(V)(w_0\lambda) + \bar \Lien U(\bar \Lien) \bigr)\varv_0
= \bar \Lien U(\bar \Lien)\varv_0$ and hence
$VT_0=VU(\bar\Lien)\varv_0=U(\bar\Lien)V\varv_0\subset\bar \Lien U(\bar \Lien)\varv_0\subsetneq T_0$.
Since $VT_0$ is a $U(\Lieg)$-submodule of $T_0$,
we get $V T_0 = 0$.
But from Duflo's theorem (\cite{Duf}),
$\Ann M(\lambda)=\Ann M(w_0\lambda)=\sum_{D\in U(\Lieg)^G}U(\Lieg)\bigl(D-\Tilde\gamma(D)(\lambda)\bigr)$.
Hence we get $V\subset \Ann M(w_0\lambda)=\Ann M(\lambda)$,
or equivalently $V M(\lambda)=0$.
It proves $\Xi_\lambda(0)=0$.

Secondly, let $S$ be an arbitrary
$\Tilde{\mathbf H}$-submodule of $A(\lambda)$.
Put \[
I\coloneqq\sum \Bigl\{\, V\subset U(\Lieg);\, 
\text{$V$ is $\ad(\Lieg)$-stable and satisfies $\Tilde\gamma(V) A(\lambda) \subset S$}\,\Bigr\}.\]
Then
\[\Xi_\lambda(S)=IM(\lambda).\]
We assert $I$ is a two-sided ideal of $U(\Lieg)$.
Indeed, $\Lieg I$ is an $\ad(\Lieg)$-stable subspace of $U(\Lieg)$
and satisfies
$\Lieg I \subset (\bar\Lien I + I\Lien) + \ad(\Lien)I + \Lieh I
\subset (\bar\Lien I + I\Lien) + I + \Lieh I$,
from which we easily deduce
$\Tilde\gamma(\Lieg I) \subset \Tilde\gamma(I) + \Lieh\,\Tilde\gamma(I)$.
Since $\left(\Tilde\gamma(I) + \Lieh\,\Tilde\gamma(I)\right) A(\lambda) \subset S$,
$\Lieg I\subset I$.
Similarly we can show $I\Lieg \subset I$. Thus $I$ is a two-sided ideal.
Moreover, since
$\Tilde\gamma\left(\Ann M(\lambda)\right)
 =\sum_{f\in S(\Lieh)^W}S(\Lieh)\bigl(f-f(\lambda)\bigr)$,
$I\supset\Ann M(\lambda)$.
Hence from \cite{Jos,BG}, we get
$I=\Ann \left(M(\lambda)/I M(\lambda)\right)$ and therefore
\begin{equation}\label{eq:YXi}
Y_\lambda(IM(\lambda))
=\Tilde{\mathbf H}\,\Tilde\gamma\left(\Ann \left(M(\lambda)/I M(\lambda)\right)\right)A(\lambda)
=\Tilde{\mathbf H}\,\Tilde\gamma(I)A(\lambda)\subset S.
\end{equation}

Lastly, we assume the above $S$ admits a
$W$-stable subspace $E$ which satisfies the conditions (a) and (b)
in the proposition.
Let $E_1\subset E$ be an irreducible $W$-submodule.
Then there exists a $W$-stable subspace $\tilde E_1$ in $S_{\mathbf H}(\Lieh)$
which is isomorphic to $E_1$
via the natural surjective $W$-homomorphism $S_{\mathbf H}(\Lieh) \rightarrow A(\lambda)$.
By the condition (b), the equivalence class of $\tilde E_1$ belongs to $\widehat{W}_\single$.
Hence it follows from
Theorem~\ref{thm:HCisomo}~(\ref{item:doubleHC}) that
there exists an $\ad(\Lieg)$-submodule $V$ of $U(\Lieg)$
which is isomorphic to an irreducible quasi-small representation of $G$
and satisfies $\Tilde\gamma(V)=\tilde E_1$.
Since $\Tilde\gamma(V)A(\lambda)=S(\Lieh)E_1\subset S$,
the above $I$ satisfies $I\supset V$ and therefore $\Tilde\gamma(I)A(\lambda)\supset E_1$.
Hence by the condition (a), we can replace
the last inclusion relation in \eqref{eq:YXi} with `$=$'.
Thus we get $Y_\lambda\circ\Xi_\lambda(S)=S$.
\end{proof}
\begin{cor}\label{cor:V-A_type_A}
Suppose $\Lieg$ is a complex simple Lie algebra of type $(A)$.
Then $\Xi_\lambda$ is injective and it holds that
$Y_\lambda\circ\Xi_\lambda(S)=S$
for any $\Tilde{\mathbf H}$-submodule $S$ of $A(\lambda)$.
\end{cor}
\begin{proof}
Under the assumption that $\Lieg$ is of type $(A)$,
all the equivalence classes of irreducible representations of $W$
belong to $\widehat{W}_\single$ (\cite{Br}).
\end{proof}
\begin{rem}
If $\lambda(\alpha^\vee)\ne1$ for all $\alpha \in \Sigma$,
then $A(\lambda)$ is irreducible (cf.~\cite{ChJ}).
Hence $\Xi_\lambda$ is not necessarily surjective.
\end{rem}

\ifx\undefined\bysame 
\newcommand{\bysame}{\leavevmode\hbox to3em{\hrulefill}\,} 
\fi 


\begin{thebibliography}{Dun2}


\bibitem[BG]{BG} J. N. Bernstein and S. I. Gelfand,
{\it Tensor products of finite and infinite dimensional
representations of semisimple Lie algebras\/},
Comp. Math. {\bf 41} (1980), 245--285.

%
%
%

\bibitem[Br]{Br} A. Broer, {\it The sum of generalized exponents and Chevalley's 
restriction theorem for modules of covariants\/}, 
Indag. Math. N. S. {\bf 6} (1995), no.4, 385--396.

%
%
%

\bibitem[Ch]{ChJ} I. Cherednik,
{\it Lectures on Knizhnik-Zamolodchikov
equations and Hecke algebras\/}, 
In:~Quantum Many-Body Problems And Representation Theory,
MSJ Mem. {\bf 1}, Math. Soc. Japan, Tokyo, 1998, 1--96.

\bibitem[Da]{Da} J. Dadok, {\it On the $C^\infty$ Chevalley's theorem\/}, 
Adv. in Math. {\bf 44} (1982), 121--131.

%
%
%
%

%
%
%
%

\bibitem[Duf]{Duf} M. Duflo,
{\it Construction of primitive ideals in an enveloping algebra\/},
In:~Lie Groups and Their Representations, edited by I. M. Gelfand,
Summer School Conference, Budapest, 1971, Halsted press, New York (1975), 77--93.

\bibitem[Dun1]{Du} C. F. Dunkl,
{\it Differential-difference operators associated to reflection groups}, 
Trans. Amer. Math. Soc. {\bf 311} (1989), 167--183.

\bibitem[Dun2]{Du2} \bysame, %
{\it Hankel transforms associated to finite groups}, 
In:~Hypergeometric Functions on Domains of Positivity,
Jack Polynomials and Applications, 
edited by Donald St. P. Richards, Proceedings, Tampa 1991,
Contemp. Math. {\bf 138} (1992), 123--138.

\bibitem[HC]{HC} Harish-Chandra, {\it Spherical functions on a semisimple Lie group, I\/},
Amer. J. Math. {\bf 80} (1958), 241--310.

\bibitem[He1]{He1} S. Helgason, {\it Differential Geometry, Lie Groups, and Symmetric Spaces\/}, 
American Mathematical Society, 2001, c1978.

\bibitem[He2]{He2} \bysame, {\it Groups and Geometric Analysis\/}, 
American Mathematical Society, 2000, c1984.

\bibitem[Joh]{Joh} K. D. Johnson, {\it Generalized Hua operators and parabolic subgroups},
Ann. of Math. {\bf 120} (1984), Springer, 477--495.

\bibitem[Jos]{Jos} A. Joseph, {\it Dixmier's problem for Verma and principal
series submodules}, J. London Math. Soc. {\bf 20} (1979), 193--204.

%
%

\bibitem[Ko1]{Ko1} \bysame,  {\it Lie group representations on polynomial rings},
Amer. J. Math. {\bf 85} (1963), 327--404.

%
%
%
%
\bibitem[Ko2]{Ko:Mat} \bysame, %
{\it On the existence and irreducibility of certain series of representations}, 
Bull. A. M. S. {\bf 75} (1969), 627--642.

\bibitem[Ko3]{Ko:Mat2} \bysame, %
{\it On the existence and irreducibility of certain series of representations}, 
In:~Lie Groups and Their Representations, edited by I. M. Gelfand,
Summer School Conference, Budapest, 1971, Halsted press, New York (1975), 231--329.

\bibitem[Ko4]{Ko:Ad} \bysame, %
{\it Clifford algebra analogue of the Hopf-Koszul-Samelson theorem,
 the $\rho$-decomposition $C(\Lieg)=\End V_\rho\otimes C(P)$,
 and the $\Lieg$-module structure of $\bigwedge\Lieg$}, 
Adv. in Math. {\bf 135} (1997), 275--350.

%
%
%
%
%

\bibitem[KR]{KR} B. Kostant and S. Rallis,  {\it Orbits and representations
associated with symmetric spaces}, 
Amer. J. Math. {\bf 93} (1971), 753--809.

\bibitem[Lu]{Lu} G. Lusztig, {\it Cuspidal local systems and graded Hecke algebras, I}, 
Publ. Math. de IHES {\bf 67} (1988), 145--202.

%
%
%
%
%
%
%
%
%

\bibitem[Od]{gHC} H. Oda, {\it Generalization of Harish-Chandra isomorphism}, 
RIMS K\^oky\^uroku, Kyoto Univ. {\bf 1294} (2002), 141--151,
in Japanese.

%
%
%

\bibitem[PRV]{PRV} K. R. Parthasarathy, R. Ranga Rao and V. S. Varadarajan,
{\it Representations of complex semi-simple Lie groups and Lie algebras},
Ann. of Math. {\bf 85} (1967), 383--429.

\bibitem[Re1]{Adg} M. Reeder, {\it Exterior powers of the adjoint representation},
Can. J. Math. {\bf 49} (1997), 133--159.


\bibitem[Re2]{Ssl} \bysame, {\it Zero weight spaces and the Springer correspondence}, 
Indag. Math. N. S. {\bf 9} (1998), 431--441.

%
%

\bibitem[Re3]{Sml} \bysame, {\it Small representations and minuscule Richardson orbits}, 
Int. Math. Res. Not. (2002), 257--275.

\bibitem[So]{So} L. Solomon, {\it Invariants of finite reflection groups}, 
Nagoya Math. J. {\bf 22} (1963), 57--64.

%
%
%
%
%
%
%


\end{thebibliography}
\end{document}